%% file: A_submission.tex
\documentclass[a4paper,9pt,parskip=half,english]{article}
\usepackage[top=60pt,bottom=60pt,left=78pt,right=76pt]{geometry}
\usepackage{babel}
\usepackage[T1]{fontenc}
\usepackage{lmodern,textcomp,latexsym}
\usepackage{amsthm,amssymb,amsmath,amsfonts,dsfont,mathrsfs,bbm,bm,mathtools}
\usepackage{graphicx,psfrag,color,rotate}
\usepackage{paralist,threeparttable}
\usepackage[footnotesize]{caption}
\usepackage{dsfont,booktabs,threeparttable}
\usepackage{longtable,booktabs}
\usepackage{cite}
\usepackage{tikz}
\usepackage{pgfplotstable}
\usepackage{pgfplots}
\usetikzlibrary{shapes.geometric}
\usetikzlibrary{matrix}
\usepackage{nicematrix}
\usepackage[title]{appendix}
\usetikzlibrary{calc,spy,external}
\usepackage{caption}
\usepackage{subcaption}
\usepackage{float}

\usepackage{titlesec}

\usepackage{booktabs}

\usepackage{gettitlestring}
\usepackage{hyperref}

\titleformat*{\section}{\large\bfseries}
\titleformat*{\subsection}{\large\bfseries}
\titleformat*{\subsubsection}{\large\bfseries}
\titleformat*{\paragraph}{\large\bfseries}
\titleformat*{\subparagraph}{\large\bfseries}

\setdefaultenum{1)}{\theenumi.1)}{}{}

\theoremstyle{plain}
\newtheorem{theorem}{Theorem}
\newtheorem{proposition}{Proposition}
\newtheorem{remark}[theorem]{Remark}
\theoremstyle{definition}
\newtheorem{definition}{Definition}

\newcommand{\diff}[1][]{\mathrm{d}#1}

\newcommand{\pd}[2]{\frac{\partial #1}{\partial #2}}

\newcommand{\T}{^{\!\top}} % \top with better spacing

\newcommand{\mR}{\mathbb{R}}

\newcommand{\veta}{\bm{\eta}}

\newcommand{\vxi}{\bm{\xi}}

\newcommand{\vph}{\bm{\varphi}}

\newcommand{\ve}{\mathbf e}
\newcommand{\vf}{\mathbf f}
\newcommand{\vg}{\mathbf g}

\newcommand{\vn}{\mathbf n}

\newcommand{\vp}{\mathbf p}
\newcommand{\vq}{\mathbf q}

\newcommand{\vw}{\mathbf w}
\newcommand{\vx}{\mathbf x}
\newcommand{\vy}{\mathbf y}

\newcommand{\vA}{\mathbf A}

\newcommand{\vC}{\mathbf C}

\newcommand{\vF}{\mathbf F}
\newcommand{\vG}{\mathbf G}

\newcommand{\vI}{\mathbf I}

\newcommand{\vK}{\mathbf K}

\newcommand{\vM}{\mathbf M}

\newcommand{\vO}{\mathbf O}

\newcommand{\vT}{\mathbf T}

\newcommand{\vX}{\mathbf X}

\newcommand{\cG}{\mathcal{G}}

\newcommand{\cP}{\mathcal{P}}

\newcommand{\cV}{\mathcal{V}}
\newcommand{\cVplus}{\mathcal{V}^+}
\newcommand{\cVminus}{\mathcal{V}^-}

\newcommand{\vzero}{\bm{0}}
\def\addlegendimage{\csname pgfplots@addlegendimage\endcsname}

\pgfplotsset{
	cycle list={%
		{draw=black,mark=star,solid},
		{draw=black, mark=square,solid},%densely dashed},
		{draw=black,mark=+,solid},%dashdotted}, %every mark/.append style={rotate=90},
		{black,mark=o},}}

\pgfdeclareplotmark{mystar}{
	\node[
	star,
	star point ratio=2.25,
	minimum size=6pt,
	inner sep=0pt,
	draw=red,% Sets the border to red
	fill=yellow,% Sets the fill to yellow
	line width=0.5pt% Define a thin line width for the border
	] {};
}

\begin{document}
	\begin{center}
		\textbf{\large Model order reduction of piecewise linear mechanical systems using invariant cones }
		
		\qquad

		\renewcommand{\thefootnote}{\fnsymbol{footnote}}
		A.Yassine Karoui$^{1}$\footnote{karoui@inm.uni-stuttgart.de},
		and
		Remco I. Leine$^{1}$\footnote{leine@inm.uni-stuttgart.de}
		
		\qquad
		
		$^{1}$Institute for Nonlinear Mechanics
		
		University of Stuttgart
		
		Pfaffenwaldring 9, 70569 Stuttgart, Germany

	\end{center}
	\vspace{5pt}
	\hrule
	\vspace{8pt}
	\textbf{Abstract}: 
	We present a methodology that extends invariant manifold theory to a class of autonomous piecewise linear systems with nonsmoothness at the equilibrium, providing a framework for model order reduction in mechanical structures with compliant contact laws. The key idea is to make the absence of a local linearization around the equilibrium tractable by leveraging the positive homogeneity property. This property simplifies the invariance equations defining the geometry of the invariant cones, from a set of partial differential equations to a system of ordinary differential equations, enabling their effective solution. We introduce two techniques to compute these invariant cones. First, an intuitive graph-style parametrization is proposed that utilizes Fourier expansions and Chebyshev polynomials to derive explicit reduced-order models in closed form. Second, an arc-length parametrization is introduced to robustly compute invariant cones with complex folding geometries, which are intractable with a standard graph-style technique. The approach is demonstrated on mechanical oscillators with unilateral visco-elastic supports, showcasing its applicability for systems with both continuous (unilateral elastic) and discontinuous (unilateral visco-elastic) unilateral force laws. \\[8pt]
	\vspace{8pt}
	\hrule
\setcounter{footnote}{0}
\renewcommand{\thefootnote}{\roman{footnote}}
\section{Introduction}\label{section:intro}
This paper introduces a full geometric parametrization of invariant cones, which serve as the nonlinear normal modes (NNMs) for a class of piecewise linear (PWL) mechanical systems. Our approach extends the classical graph-style parametrization method of invariant manifolds, pioneered by Shaw and Pierre \cite{Shaw1993} for smooth nonlinear systems, to the realm of PWL dynamics. This parametrization reveals the geometric shape of the invariant cone and enables the derivation of closed-form, mathematically exact reduced-order models (ROMs) by constraining the system's dynamics to its invariant cone. To the best of our knowledge, this work presents the first explicit invariant cone parametrization and the first derivation of such two-dimensional ROMs for nonlinearizable PWL systems.

The analysis of high-dimensional dynamical systems, particularly those arising from complex finite element models, often necessitates the use of model order reduction methods (MOR) to achieve qualitative understanding while reducing the computational burden~\cite{Kerschen2009,Touze2021}. For systems exhibiting significant nonlinear behavior, traditional linear analysis methods are frequently insufficient, as linear subspaces are inherently not invariant under the flow of the nonlinear full system~\cite{Haller2016}. This has motivated the development of nonlinear MOR methods that reduce the full system's trajectories onto attracting, forward-invariant manifolds, thereby providing a mathematically justifiable and highly accurate ROMs.
In the context of MOR for smooth nonlinear systems, NNMs emerged as a powerful tool. The concept offers a robust framework for analyzing nonlinear structures by drawing a clear conceptual relation to their linear counterparts, the Linear Normal Modes (LNMs). While fundamental properties of LNMs such as orthogonality, modal superposition and frequency-energy independence are lost in the nonlinear setting, NNMs inherit the most crucial property: invariance. This means that trajectories initiated on an NNM remain on it indefinitely. This invariance property makes NNMs very suitable for MOR, as they provide a more effective basis than LNMs for capturing the system's dynamics, especially for large-amplitude vibrations and in the presence of strong nonlinearities. This modern, manifold-based understanding of NNMs evolved from a rich history of research aimed at extending linear modal concepts. The term "nonlinear mode" can be traced back to the pioneering work of Rosenberg in the 1960s \cite{Rosenberg1960, Rosenberg1961, Rosenberg1962, Rosenberg1964, Rosenberg1964a, Rosenberg1966}, who defined NNMs as synchronous periodic oscillations in conservative nonlinear systems, viewing them as a natural generalization of LNMs. This classical definition, which required that all material points of the system oscillate in unison, was limited to undamped mechanical systems. Further foundational work on NNMs was conducted in the 1970s by researchers such as Rand \cite{Rand1971,Rand1971a,Rand1974} and Manevitch and Mikhlin \cite{Manevich1972}. The 1990s witnessed a resurgence of interest in NNM research, spearheaded by the work of Vakakis et al. \cite{Caughey1990,Vakakis1992,Vakakis1996,Vakakis1997} and, most notably, Shaw and Pierre \cite{Shaw1991,Shaw1993,Shaw1994}. By introducing the invariant manifold definition, Shaw and Pierre broadened the scope of NNM analysis to include damping \cite{Shaw1993} and external forcing \cite{Jiang2005}, laying the ground work for the modern computational tools and extensive research that has established NNMs as a cornerstone of nonlinear vibration analysis and model reduction \cite{Touze2006,Amabili2007,Touze2021}. 

The modern, widely accepted definition of an NNM, pioneered by Shaw and Pierre, characterizes it as an invariant manifold in the system's phase space that is tangent to a linear modal subspace at the equilibrium. While this geometric view was a pivotal generalization, it was later shown that for dissipative systems, this definition is not sufficient to ensure uniqueness; infinitely many such invariant manifolds can exist. This ambiguity was formally resolved by the theory of Spectral Submanifolds (SSMs), introduced by Haller and Ponsioen \cite{Haller2016, Ponsioen2020}. An SSM is defined as the unique, smoothest invariant manifold that continues a given spectral subspace into the nonlinear regime. Its existence and uniqueness are guaranteed under specific non-resonance conditions related to the system's spectrum. This established SSMs as the definitive target for model reduction in smooth damped nonlinear systems.

The computational framework for constructing these manifolds is the parametrization method, a powerful mathematical technique that provides a unified perspective on different computational approaches~\cite{Haro2016}. Within this framework, earlier methods can be understood as distinct "styles" of parametrization. The approach originating with Shaw and Pierre~\cite{Shaw1993} is known as the graph style, where the goal is to find an explicit functional relationship expressing the slave modal coordinates in terms of the master coordinates. While intuitive, this style is inherently limited as it cannot parametrize a manifold beyond a folding point \cite{Touze2021}. In contrast, the normal form style, developed extensively for vibratory systems by Touzé and coworkers \cite{Touze2003a, Touze2004,Touze2014}, uses a nonlinear change of coordinates for the entire phase space to simplify the reduced dynamics as much as possible, retaining only the essential resonant terms. This results in a more complex nonlinear mapping but yields a simpler reduced model and is capable of accurately capturing complex manifold geometries, including folds. Recent advancements from both schools of thought have focused on developing direct methods—such as Direct Normal Form (DNF) \cite{Vizzaccaro2021,Opreni2021} and direct SSM computation \cite{Jain2021,Veraszto2020}—that apply these advanced concepts directly to large-scale Finite Element (FE) models, making them suitable tools for modern engineering analysis. Despite their different computational philosophies, these approaches all belong to a powerful geometric framework built upon two fundamental assumptions: first, the system dynamics are smooth (i.e., sufficiently differentiable), and second, a unique linearization exists at the equilibrium point, providing the LNM subspaces to which the NNMs are tangent.

A significant class of mathematical models for problems in engineering, however, are characterized by nonsmooth dynamics, and thus violates these smoothness conditions at their core. Of particular interest are PWL systems, which arise in mechanical systems with bilinear stiffness properties such as cracked beams and rotating shafts \cite{Zuo1994, Chati1997, Vestroni2008}, and oscillators with clearance and unilateral elastic contact \cite{Chen1996,Butcher1999}.

The development of analysis methods for NNMs in PWL systems has evolved along two distinct paths dictated by the nature of the system's equilibrium. For the class of inhomogeneous or piecewise affine (PWA) systems, the equilibrium resides within a region of purely linear behavior, with nonlinearities activated only when oscillations exceed a certain clearance or gap. This crucial feature causes the foundational assumption of the smooth theory—the existence of a local, unique linearization—to remain valid. Consequently, the invariant manifold framework of Shaw and Pierre can be adapted and extended to this nonsmooth setting. Early work by Chen and Shaw \cite{Chen1996} successfully applied this invariant manifold approach, using the switching hyperplane as a Poincaré section to construct the NNM via asymptotic expansions. This line of research was further advanced by Jiang, Pierre, and Shaw \cite{Jiang2004}, who developed a Galerkin-based method to compute accurate large-amplitude NNMs, and by Butcher \cite{Butcher1999}, who investigated the precise influence of the clearance-to-amplitude ratio on the NNM frequency. More recently, computational techniques such as the Harmonic Balance Method in the context of inhomogeneous PWL systems have been employed by Kim et al. \cite{Kim2005} and Moussi et al. \cite{Moussi2015}. A common thread in these contributions is their reliance on the locally linear nature of the equilibrium as an anchor point for their constructive methods. 

In stark contrast, \textit{homogeneous} PWL systems present a qualitatively different challenge. Here, the equilibrium is located on the switching manifold, which causes the Jacobian to jump, thereby rendering the system nonlinearizable, since a local linearization of the system in the form $\dot{x} = A x$ cannot be obtained. Classical approaches that rely on Taylor series expansions around an equilibrium are thus fundamentally inapplicable, as the vector field is non-differentiable precisely at the point of expansion. This challenge, however, is accompanied by a unique structural property: positive homogeneity. This property, which leads to the remarkable feature of amplitude-independent NNM frequencies, was first explored by Zuo and Curnier \cite{Zuo1994}, who extended Rosenberg's definition to unforced conservative systems with homogeneous PWL restoring forces. Along these lines, Chati et al. \cite{Chati1997} analyzed cracked beam models and introduced the notion of bilinear frequency, an analytical approximation of the energy-independent eigenfrequencies of homogeneous PWL systems, whose accuracy is guaranteed only for weak nonlinearities as shown by Attar et al. in \cite{Attar2017}. The investigations by Vestroni and Casini et al. \cite{Casini2011,Casini2012,Casini2013,Vestroni2008} further deepened the understanding of these systems by applying perturbation methods to uncover complex bifurcation phenomena and the existence of superabundant modes. While these works successfully characterized the periodic orbits, they did not aim to develop a general geometric framework for their full parametrization for MOR purposes.

The preceding analyses of PWL systems, across both inhomogeneous and homogeneous classes, share a crucial, underlying assumption: the continuity of the system's vector field. The introduction of discontinuities—for instance, from a unilateral damper in addition to a unilateral spring—poses an additional challenge, as the dynamics are no longer described by an ordinary differential equation with a continuous right-hand side. Significant progress in understanding discontinuous dynamics has been made in the context of vibro-impact systems, where discontinuities manifest as instantaneous, event-driven velocity jumps governed by impact laws \cite{Thorin2015, Thorin2017,LeThi2018}. However, a geometric manifold-based framework for MOR in systems with state-dependent, discontinuous restoring forces, characteristic of certain friction models or unilateral dampers, remains largely unexplored, particularly within the challenging nonlinearizable context of homogeneous PWL systems.

The conceptual key to a unified approach for homogeneous PWL systems can be found in stability theory. The pioneering work of Carmona et al. \cite{carmona2005bifurcation}, later extended to higher-dimensional continuous and discontinuous homogeneous PWL systems by Küpper and Hosham \cite{Kuepper2008, Kuepper2011, Hosham2016}, identified invariant cones as the fundamental organizing structures for the dynamics of these nonlinearizable systems. Their framework was designed to reduce the high-dimensional stability problem to a far simpler one. By identifying invariant directions on a Poincaré section (typically the switching manifold), they showed that orbits starting on these directions span a 2D surface—the cone—and ultimately return to the same direction. This allowed the stability of the equilibrium to be determined by a simple 1D map governing the radial dynamics along the cone.
It is critical, however, to distinguish this type of stability analysis from the goal of MOR. While their work has used the cone's invariance to track stability, its objective has not been to provide a geometric basis onto which the full system dynamics can be projected to yield a reduced-order model. A paradigm shift is thus required: from using invariant cones for stability analysis to developing a constructive method for their full geometric parametrization, a necessary step for deriving explicit ROMs.

The present work aims to bridge this gap: we shift the focus from stability analysis of trajectories on the invariant cone to geometric parametrization and model order reduction. Our main goal is to obtain the complete geometric shape of the entire invariant cone, which can then be used to derive lower-dimensional dynamical models that capture the full system's behavior. We propose two distinct approaches: an intuitive projection method, reminiscent of the graph-style parametrization method of smooth invariant manifold, and a more robust, general arc-length method.

The paper is structured as follows. We begin in Section \ref{section:background} by defining the class of homogeneous and generally discontinuous PWL systems under study. Section \ref{section:invcones} then reviews the foundational concept of their invariant cones. Our constructive approach unfolds in two stages. First, in Section \ref{section:graphstyle}, we develop an intuitive projection method adapted from the classical graph-style parametrization of smooth manifolds. We then perform a critical analysis in Section \ref{Section_example_GS_limitations}, revealing the geometric limitations of this approach and thereby motivating the need for a more robust technique. This leads to our primary methodological contribution in Section \ref{Section_ArcLength}: the arc-length parametrization capable of accurately capturing complex cone geometries, including folds where the projection method fails. The efficacy of this framework is demonstrated in Section \ref{Section_NumValidation} through numerical examples of both continuous and discontinuous systems. There, we derive explicit two-dimensional ROMs and validate their accuracy against the full system dynamics. Finally, Section \ref{Section_Conclusion} summarizes our contributions and outlines avenues for future research.

\section{Homogeneous Piecewise Linear Mechanical Systems}\label{section:background}
\subsection{Mechanical system}
\begin{figure}
	\centering	
	\includegraphics[scale=1]{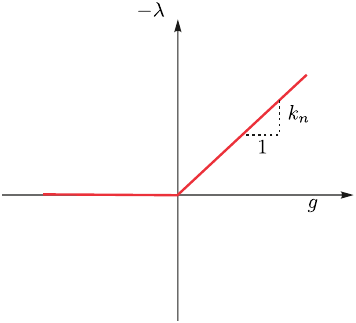}	
	\caption{Characteristic of the unilateral elastic contact law ($c_n=0$).}
	\label{fig_unilateral_elastic_force_law}
\end{figure}

We consider a class of autonomous mechanical systems, which essentially consist of a linear multi-degree-of-freedom oscillator with constant system matrices, which is additionally equipped with a single unilateral visco-elastic support placed exactly at the equilibrium. This unilateral support models contact within the system (e.g., a crack in a beam) or a one-sided support where damping can be considered in the contact mechanism. Such systems may be expressed by
\begin{align}\label{eq:eom_2ndorder}
	\vM \ddot{\vq} + \vC \dot{\vq} + \vK \vq = \vf_n,
\end{align}
where $\vq \in \mathbb{R}^N$ is the generalized position vector, $\vM=\vM\T \in \mathbb{R}^{N \times N}$ is the positive definite mass matrix, $\vC=\vC\T \in \mathbb{R}^{N \times N}$ is the damping matrix, $\vK=\vK\T \in \mathbb{R}^{N \times N}$ is the positive definite stiffness matrix. The right-hand side consists of the nonsmooth function $\vf_n$ describing the  contact force, which arises from the interaction of the linear structure with a massless discontinuous support consisting of a unilateral spring-damper force element. The unilateral force element is described within the framework of nonsmooth dynamics \cite{Glocker2001}. Herein, a generic force element is considered with kinematic scalar quantity $g(\vq)$ and its kinetic quantity $\lambda$ such that the contribution to the virtual work of the system is 
\begin{equation}
	\delta W = \lambda \delta g.
\end{equation}
Correspondingly, the contribution of the force element to the equation of motion is expressed by the generalized force 
\begin{equation}
		\vf_n = \vw \lambda, 
\end{equation}
where 
\begin{equation}
	\vw = (\pd{g}{\vq})\T
\end{equation}
is the generalized force direction ensuring 
\begin{equation}
	\delta W = \lambda \delta g = \lambda \pd{g}{\vq} \delta \vq = \vf_n\T \delta \vq.
\end{equation}
Within the context of this paper, we will consider a visco-elastic support modeled as a unilateral spring damper element which can only transmit a compressive force. As kinematic quantity, we choose $g$ to denote the compression of the element, i.e., the penetration depth of the visco-elastic support. In turn, the negative of kinetic quantity $-\lambda\geq 0$ has the physical meaning of contact force. We assume a linear visco-elastic element with contact force 
\begin{equation}
	-\lambda = k_n g + c_n \dot{g} \geq 0, \quad \mbox{(contact)}
\end{equation}
where $k_n>0$ and $c_n>0$ are the spring and damper constants and $\dot{g}$ denotes the rate of penetration. Furthermore, we will assume a linear kinematic dependence of the penetration on the generalized coordinates, i.e., 
\begin{equation}
	g = \vw\T \vq,
\end{equation}
with constant generalized force direction $\vw$ such that $\dot{g}=\vw\T \dot{\vq}$.
%This function is defined as
%\begin{equation}
%	\vf_n = \vw \lambda,
%\end{equation}
%where $\vw$ is the generalized force direction, and $\lambda$ is the scalar local force quantity acting on the unilateral support. If we consider that the support sustains only compressive forces, we can introduce the kinematic variables $g = \vw\T \vq$ describing the compression length and the compression rate $\dot{g}=\vw\T \dot{\vq}$ with the corresponding kinetic quantity $\lambda$ being defined as
%\begin{equation}
%	-\lambda = k_n g + c_n \dot{g} > 0.
%\end{equation}
%Herein, $k_n>0$ and $c_n>0$ are the stiffness and damping parameters for the visco-elastic support. 
In a state of contact, the conditions $g>0$ (penetration) and $-\lambda>0$ (compressive force) hold. The bilinear characteristic of the unilateral force element is shown in Figure \ref{fig_unilateral_elastic_force_law} for the simple case where $c_n=0$.
A mechanical system with such a support has two possible modes: the unilateral spring and damper are active during contact and inactive in the no-contact phases. The support itself being a first-order system relaxes to its equilibrium state if the contact is inactive. If we assume that the relaxation time of the support is much smaller than the time between two contact events, then we can neglect the free motion of the support and assume it is always at rest at the moment a contact is initiated \cite{Lein2004}. This assumption implies that $\vf_n=\vf_n(\vq,\dot{\vq})$ resulting in a system of second-order differential equations with a state-space partitioned by two scalar switching functions:
\begin{enumerate}
	\item \textbf{Position-based switching function, $h_{\alpha}$:} This function determines if the system is at or beyond the contact boundary. It is defined by the penetration depth:
	\begin{align}
		h_{\alpha}(\vq,\dot{\vq}) = g = \vw\T\vq.
	\end{align}
	\item \textbf{Force-based switching function, $h_{\beta}$:} This function represents the magnitude of the contact force, which must be compressive for the support to be active. It depends on both generalized position and velocity:
	\begin{align}
		h_{\beta}(\vq, \dot{\vq}) =  k_n g + c_n \dot{g} =  k_n(\vw\T\vq) + c_n(\vw\T\dot{\vq}).
	\end{align}
\end{enumerate}
These functions define the two distinct regions of the system:
\begin{itemize}
	\item \textbf{Contact Region $\cVplus$}: The system is in contact with the support, which exerts a repulsive force. This requires both a positive generalized displacement ($h_{\alpha} > 0$) and a positive contact force ($h_{\beta} > 0$), i.e., $\cVplus = \{ (\vq,\dot{\vq}) \in \mR^{2N} \, | \, h_{\alpha}(\vq) > 0 \land h_{\beta}(\vq,\dot{\vq}) >0 \}$. 
	\item \textbf{No-Contact Region $\cVminus$}: The system is not in active contact with the support. This region is the complement of the contact region, where the nonsmooth force vanishes, $\vf_n(\vq,\dot{\vq}) = 0$, i.e., $\cVminus = \{ (\vq,\dot{\vq}) \in \mR^{2N} \, | \, h_{\alpha}(\vq) < 0 \lor h_{\beta}(\vq,\dot{\vq}) < 0 \}$. 
\end{itemize}
Crucially, our modeling assumption implies that contact is initiated only when a mass hits the support at its rest position, i.e., through a zero crossing of the function $h_\alpha$. The no-contact mode corresponds to the conditions $h_\alpha<0$ or $h_\beta<0$.
Substituting the expression for the nonsmooth force into Eq. \eqref{eq:eom_2ndorder} yields the complete PWL equations of motion. By rearranging terms, we can express the dynamics in each region as a homogeneous linear system with different system matrices:
\begin{equation} \label{eq:eom_piecewise}
	\begin{cases}
		\vM\ddot{\vq} + (\vC + c_n\vw\vw\T)\dot{\vq} + (\vK + k_n\vw\vw\T)\vq = 0 & \text{for } (\vq,\dot{\vq}) \in \cVplus \\
		\vM\ddot{\vq} + \vC\dot{\vq} + \vK\vq = 0 & \text{for } (\vq,\dot{\vq}) \in \cVminus
	\end{cases}
\end{equation}
The switching boundary $\Sigma$, which divides the state-space $\mR^{2N}$ in the subspaces $\mathcal{V}_\pm$ consists of the union of two surfaces $\Sigma_\alpha$ and $\Sigma_\beta$. The switching boundary $\Sigma_\alpha$ is defined by 
\begin{equation}
	\Sigma_\alpha = \{ (\vq,\dot{\vq}) \in \mR^{2N} \, | \, h_\alpha(\vq) = 0 , h_\beta(\vq,\dot{\vq}) \geq 0  \},
\end{equation}
and has the normal $\vn_\alpha(\vq,\dot{\vq}) = \nabla h_\alpha(\vq,\dot{\vq}) = [ \vw\T , \bm{0}_{N\times 1}\T ]\T$. The switching boundary $\Sigma_\beta$ is defined by 
\begin{equation}
	\Sigma_\beta = \{ (\vq,\dot{\vq}) \in \mR^{2N} \, | \, h_\alpha(\vq) \geq 0 , h_\beta(\vq, \dot{\vq}) = 0 \},
\end{equation} 
and has the normal $\vn_\beta = \nabla h_\beta(\vq,\dot{\vq}) = [ k_n \vw\T , c_n \vw\T]\T$. The switching boundary $\Sigma$ is nonsmooth at the origin. 
The system is discontinuous upon entering the contact region at  the boundary $\Sigma_\alpha$, where it holds that the switching function $h_\alpha = \vw\T \vq = 0$. Thus, the contact force jumps from zero to a strictly positive value $- \lambda_n= c_n \vw\T \dot{\vq}$ due to the positive penetration velocity ($\dot{g}>0$). This jump on the right-hand side of the equation of motion, i.e. due to the term $c_n\vw\vw\T$, causes the discontinuity in the dynamics, thereby classifying the system as a discontinuous PWL system.

An important special case is that of a continuous piecewise linear (CPL) system, which is recovered by setting the unilateral damping coefficient to zero, $c_n = 0$. Physically, this corresponds to a purely elastic unilateral support (a unilateral spring). In this limit, the switching function $h_{\beta}$ becomes purely position-dependent: $h_{\beta}(\vq) = k_n(\vw\T\vq)$. The two conditions for contact, $h_{\alpha} \ge 0$ and $h_{\beta} \ge 0$, become identical, simplifying to the single condition $\vw\T\vq \ge 0$ (since $k_n > 0$). The term causing the discontinuity in the damping matrix vanishes, and the equations of motion simplify to:
\begin{equation}\label{eq:eom_damp_CPL}
	\vM\ddot{\vq} + \vC\dot{\vq} + \vK\vq = \begin{cases} - \vw (k_n \vw\T\vq) & \text{if } g = \vw\T\vq \geq 0 \\ 0 & \text{if } g = \vw\T\vq < 0 \end{cases}.
\end{equation}
This can be written more compactly using the $\max$ function as
\begin{align}\label{eq:eom_cons_CPL}
	\vM \ddot{\vq} + \vC \dot{\vq} + \vK \vq = - \vw \max( k_n \vw\T\vq, 0),
\end{align}
which is the formulation used for CPL systems in previous works such as \cite{Karoui2025}. In this continuous case, only the stiffness matrix changes between regions, while the mass and damping matrices remain constant, ensuring the continuity of the vector field across the switching boundary. Our general formulation in Eq. \eqref{eq:eom_piecewise} therefore naturally includes CPL systems modeling a purely elastic unilateral contact as a special case.
\begin{remark}[Underlying Conservative Dynamics]
	A crucial aspect of our analysis concerns the distinction between the full dissipative system and its underlying conservative dynamics. The discontinuity in the vector field of Eq.~\eqref{eq:eom_piecewise} arises exclusively from the velocity-dependent term $c_n\vw\vw\T\dot{\vq}$ associated with the damping in the unilateral support. The underlying undamped dynamics are obtained by setting all damping coefficients to zero, i.e., $\vC=\vO$ and $c_n=0$. In this limit, the system reduces to the CPL case described by
	\begin{align}
		\vM \ddot{\vq} + \vK \vq = - \vw \max( k_n \vw\T\vq, 0).
	\end{align}
	Therefore, while the damped system is discontinuous, its underlying conservative dynamics are inherently continuous. This distinction is fundamental for our work, as we will later define and compute the system's \textit{conservative NNMs} based on the periodic solutions of this underlying CPL system.
\end{remark}

\subsection{State-Space Representation and Filippov Framework}
The second-order equations of motion~\eqref{eq:eom_piecewise} can be rewritten in a first-order form. Defining the state vector as $\vx = [\vq\T, \dot{\vq}\T]\T \in \mR^n$ with $n=2N$, the dynamics are governed by the discontinuous, homogeneous PWL system
\begin{align}\label{eq:eom_firstorder}
	\dot{\vx} = \vF(\vx) = \begin{cases} \vA^- \vx & \text{for } \vx \in \cVminus \\ \vA^+ \vx & \text{for } \vx \in \cVplus \end{cases},
\end{align}
where the system matrices for the no-contact ($\vA^-$) and contact ($\vA^+$) regions are given by
\begin{align}\label{eq:sys_matrices}
	\vA^- = \begin{bmatrix} \vO & \vI \\ -\vM^{-1}\vK & -\vM^{-1}\vC \end{bmatrix}, \quad \vA^+ = \begin{bmatrix} \vO & \vI \\ -\vM^{-1}(\vK + k_n\vw\vw\T) & -\vM^{-1}(\vC + c_n\vw\vw\T) \end{bmatrix}.
\end{align}
The dynamics described by Eq.~\eqref{eq:eom_firstorder} belongs to the class of Filippov systems, which provide a rigorous mathematical framework for analyzing differential equations with discontinuous right-hand sides \cite{Filippov1988}. The solution $\vx(t)$ in the sense of Filippov is an absolutely continuous function of time that satisfies the differential inclusion
\begin{align}
	\dot{\vx} \in \left\{ 
	\begin{alignedat}{2} % {2} indicates two pairs of alignable columns
		& \vA^- \vx , && \qquad \vx \in \cVminus \\ 
		& \overline{\text{co}}\{\vA^-\vx, \vA^+\vx\}, && \qquad \vx \in \Sigma = \Sigma_\alpha \cup \Sigma_\beta \\ 
		& \vA^+ \vx , && \qquad \vx \in \cVplus
	\end{alignedat}
	\right.
\end{align}
almost everywhere (i.e., not for a Lebesgue negligible set of time-instants), where $\overline{\text{co}}$ denotes the closed convex hull. For any state $\vx$ not on the switching boundary $\Sigma$, the set on the right-hand side contains a single vector field. However, for states on the boundary, the behavior of trajectories is more complex and depends on the orientation of the vector fields relative to the boundary.

Let us consider a generic locally smooth switching surface defined by $h(\vx)=0$ with normal vector $\vn(\vx) = \nabla h(\vx)$. The behavior of trajectories at the surface is determined by the projection of the vector fields $\vA^{\pm}\vx$ onto the normal direction $\vn$. Three distinct scenarios are possible \cite{Lein2004}:
\begin{itemize}
	\item \textbf{Transversal Crossing:} If both vector fields point to the same side of the switching surface, a trajectory crosses it transversally. This occurs when the projections of the vector fields onto the normal have the same sign:
	\begin{align}
		(\vn\T \vA^-\vx)(\vn\T \vA^+\vx) > 0.
	\end{align}
	\item \textbf{Attractive Sliding Mode:} If the vector fields on both sides of the surface point towards it, trajectories become trapped and slide along the boundary. This occurs when:
	\begin{align}
		\vn\T \vA^-\vx > 0 \quad \text{and} \quad \vn\T \vA^+\vx < 0.
	\end{align}
	The dynamics during sliding are governed by a sliding vector field, $F_{\Sigma}$, which is a specific convex combination of the two vector fields that is tangent to the surface.
	\item \textbf{Repulsive Sliding Mode:} If both vector fields point away from the surface, trajectories are repelled from it. This condition is met if:
	\begin{align}
		\vn\T \vA^-\vx < 0 \quad \text{and} \quad \vn\T \vA^+\vx > 0.
	\end{align}
	In this case, solutions starting on the boundary are non-unique in forward time.
\end{itemize}
For the specific class of mechanical systems considered here, we can prove that trajectories can only exhibit transversal crossing behavior, and thus that sliding modes are excluded.

\begin{proposition}[Absence of Sliding Modes]\label{prop:no_sliding}
	The homogeneous PWL system defined by Eq. \eqref{eq:eom_firstorder} does not exhibit sliding modes on either of its switching boundaries, $\Sigma_{\alpha}$ or $\Sigma_{\beta}$. Trajectories can only cross these boundaries transversally or, in specific instances where the velocity component normal to the boundary is zero, graze them tangentially.
\end{proposition}
\begin{proof}
	We analyze the behavior at each boundary separately.
	\begin{enumerate}
		\item \textit{Behavior at the position-based switching boundary $\Sigma_{\alpha}$:} The switching function is $h_{\alpha}(\vx) = \vw\T\vq = [ \vw\T, \bm{0}\T ]\vx$, with normal $\vn_{\alpha} = [ \vw\T, \bm{0}\T]\T$. The projection of an arbitrary vector field $\vA^{\pm}\vx$ onto this normal direction is
		\begin{align}
			\vn_{\alpha}\T (\vA^{\pm}\vx) = [ \vw\T, \vO\T] \begin{pmatrix} \dot{\vq} \\ \ddot{\vq}^{\pm} \end{pmatrix} =  \vw\T\dot{\vq}.
		\end{align}
		Crucially, the result is independent of whether the contact or no-contact dynamics are considered. Since the projections are identical, $\vn_{\alpha}\T (\vA^+\vx) = \vn_{\alpha}\T (\vA^-\vx)$, their product can never be negative. The condition for a sliding mode is therefore never met. Trajectories must cross if $ \vw\T \dot{\vq} \neq 0$ or graze $\Sigma_{\alpha}$ if $ \vw\T \dot{\vq} = 0$.
		
		\item \textit{Behavior at the force-based switching boundary $\Sigma_{\beta}$:} The switching function is $h_{\beta}(\vx) = k_n(\vw\T\vq) + c_n(\vw\T\dot{\vq}) = [ k_n\vw\T,  c_n\vw\T]\vx$. For any state $\vx \in \Sigma_{\beta}$, the condition $h_{\beta}(\vx) = 0$ holds by definition. Examining the difference between the two vector fields from Eq.~\eqref{eq:sys_matrices}, we find
		\begin{align}
			(\vA^+ - \vA^-)\vx = \begin{bmatrix} \vO & \vO \\ -\vM^{-1}(k_n\vw\vw\T) & -\vM^{-1}(c_n\vw\vw\T) \end{bmatrix} \vx = \begin{bmatrix} \vO \\ - \vM^{-1}\vw h_{\beta}(\vx) \end{bmatrix}.
		\end{align}
		For any $\vx \in \Sigma_{\beta}$, this difference is identically zero. Therefore, the vector field $\vF(\vx)$ is \textit{continuous} across the boundary $\Sigma_{\beta}$. This continuity trivially excludes the possibility of sliding modes.
	\end{enumerate}
\end{proof}
The absence of sliding modes is a crucial feature of the systems under investigation. It implies that the system's evolution is always governed by one of the two linear dynamics, with instantaneous switching between them. This property, combined with the positive homogeneity of the vector field, allows for the existence of invariant cones and makes the system amenable to the parametrization techniques developed in the subsequent sections. The remainder of this paper will therefore proceed under this established no-sliding mode assumption.
\subsection{Positive Homogeneity}
A fundamental property of the class of systems defined by Eq.~\eqref{eq:eom_firstorder}, which is essential for the study of their NNMs, is positive homogeneity.
A function $\vf: \mR^n \mapsto \mR^m $ is said to be positively homogeneous of degree $k \in \mR_0^+$ if $\vf(\beta \vx) = \beta^k \vf(\vx), \, \forall \beta > 0$. Correspondingly, a vector field $\vF(\vx)$ is said to be positively homogeneous of degree one if it satisfies the scaling property
\begin{align}\label{eq:homogeneity_def}
	\vF(\beta \vx) = \beta \vF(\vx) \quad \forall \vx \in \mR^n, \forall \beta > 0.
\end{align}
A direct consequence of this property is that the system's trajectories scale in the same manner: if $\vx(t)$ is a solution, then $\tilde{\vx}(t) = \beta \vx(t)$ is also a solution for any $\beta > 0$. We now verify that our general discontinuous system preserves this property.

First, we examine the switching functions $h_{\alpha}(\vx)$ and $h_{\beta}(\vx)$. Due to their linear dependence on the state vector $\vx$, they are inherently homogeneous of degree one:
\begin{align}
	h_{\alpha}(\beta \vx) = [ \vw\T, \vO\T](\beta \vx) = \beta h_{\alpha}(\vx), \\
	h_{\beta}(\beta \vx) = [ k_n\vw\T, c_n\vw\T](\beta \vx) = \beta h_{\beta}(\vx).
\end{align}
This implies that the conditions defining the contact and no-contact regions, $\cVplus$ and $\cVminus$, are invariant to scaling by a positive constant. If a state $\vx$ belongs to a given region, so does the scaled state $\beta \vx$. Geometrically, this means that the regions $\mathcal{V}_{\pm}$ are cones in the state space $\mR^{2N}$.

Next, we verify the homogeneity of the vector field $\vF(\vx)$ itself. The matrices $\vA^{\pm}$ are constant, thus the vector fields in both regions are linear functions of the state. It immediately follows that:
\begin{equation}
	\text{If } \vx \in \mathcal{V}^{\pm}, \quad \vF(\beta \vx) = \vA^{\pm}(\beta \vx) = \beta(\vA^{\pm}\vx) = \beta \vF(\vx).
\end{equation}
Since the vector field is homogeneous in each piecewise domain and the domains themselves are conic, the entire dynamics governed by Eq.~\eqref{eq:eom_firstorder} is positively homogeneous of degree one.

The property of positive homogeneity implies that any trajectory of the system can be scaled to generate a one-parameter family of geometrically similar trajectories. Consequently, any invariant set containing oscillatory motion must be an \textit{invariant cone} in the state space, formed by the union of these scaled solutions \cite{Kuepper2008, Karoui2025}.

Furthermore, the absence of sliding modes, established in Proposition~\ref{prop:no_sliding}, ensures that trajectories transition between the linear subsystems without being trapped on the boundary, thus simplifying the geometric shape of the cones across the switching boundaries. These two properties, positive homogeneity and the absence of sliding modes, provide the theoretical foundation for our approach developed in subsequent sections, in which we aim at finding explicit geometric parametrizations of invariant cones to use for model order reduction purposes. Before we develop our proposed techniques, the next section is dedicated at briefly introducing the concept of invariant cones for homogeneous $n-$dimensional PWL systems.

\section{Invariant Cones in Homogeneous PWL Systems: From Elementary to Multi-Crossing Structures}\label{section:invcones}
Systems of the form \eqref{eq:eom_firstorder}, which are characterized by the homogeneity property, allow for the existence of special invariant sets known as cones. Analogous to the role of center manifolds in reducing the complexity of smooth systems, invariant cones capture the essential dynamics, allowing for a rigorous stability and bifurcation analysis on a lower-dimensional surface. The concept was formally developed for general PWL systems (not necessarily continuous) by Küpper et al.~\cite{Kuepper2008, Kuepper2011, Weiss2012} and has been extensively studied in the context of three-dimensional CPL systems by Carmona et al.~\cite{carmona2005bifurcation, Carmona2006} and Huan et al.~\cite{Huan2016,Huan2017,Huan2017a}. Here, we restrict our attention to systems characterized by trajectories without sliding behavior on the switching boundary; for a detailed treatment of invariant cones in the presence of sliding modes, the reader is referred to \cite{Weiss2015}. 

This section aims to provide a comprehensive overview of this topic while also formalizing a key generalization. First, we will provide a review of the established theory of \textit{elementary invariant cones}, which arise from the simplest oscillatory trajectories involving a single crossing of the switching hyperplane into each linear region per cycle. We will give their definition via the Poincaré map and discuss their use to reduce the stability analysis.
We will then build upon this foundation to provide a formal treatment of the \textit{$k$-crossing invariant cones}. This concept, a straightforward generalization of elementary invariant cones, has been used implicitly in the literature \cite{Kuepper2013} but has not been formally detailed. This framework is based on a generating cycle of rays on the switching manifold and allows for the classification of more complex invariant structures arising from orbits with multiple crossings per full cycle. Although not strictly practical from a computational point of view, the concept will enable us to treat more complex types of NNMs resulting from bifurcations and internal resonances as shown in Section \ref{Section_NumValidation}. 
\subsection{The Elementary Invariant Cone: A Review of the Established Concept}
For PWL systems with two zones, the state space can be decomposed in two subregions $\cV^\pm$ separated by the switching hyperplane $\Sigma$ as follows
\begin{equation}
	\mR^n = \cVminus \cup \Sigma \cup \cVplus, \quad \mbox{with } \quad  \cV^\pm := \{ \vx \in \mR^n \mid \ve_1\T \vx \gtrless 0 \}, \quad \mbox{and } \quad  \Sigma := \{ \vx \in \mR^n \mid \ve_1\T \vx = 0 \}. \nonumber 
\end{equation}
In the following, oscillatory mechanical systems are of interest in this work and therefore, we assume the presence of complex eigenvalues in both subsystems, which in turn introduces rotation in the state space and excludes trajectories confined to only one linear subsystem. Following this argumentation and assuming that no sliding behavior can occur, a solution that evolves in one subregion and reaches $\Sigma$ necessarily crosses it towards the other subregion. Let $\vph(t,t_0,x_0) := \vx(t)$ be the solution of \eqref{eq:eom_firstorder} with initial condition $\vx(t_0)= \vx_0$.
Depending on the direction of crossing, the switching hyperplane $\Sigma$ can be decomposed in the disjunct subsets 
\begin{equation}
	\Sigma^\pm := \{ \vx \in \Sigma \mid \exists \, \tau > 0 \, : \, \vph(t,0,\vx) \in \cV^\pm, \, \forall \,  t \in (0,\tau) \},
\end{equation}
and the origin $\{\bm{0}\}$ being the equilibrium of the system, with $\Sigma = \Sigma^- \cup \Sigma^+ \cup \{\bm{0}\}$. The subsets $\Sigma^\pm$ are defined such that if $\vx \in \Sigma^\pm$, then the solution will continue in $\cV^\pm$.  

The intersection point of the solution and the switching hyperplane can then be used as initial condition for the following part of the orbit. Since the solution of the homogeneous PWL system \eqref{eq:eom_firstorder} in each linear subregion of the state space can be described using matrix exponential functions, we can capture the sequential motion by Poincaré half-maps defined as follows
\begin{definition}[Poincaré Half-Maps]
	Let $\vxi \in \Sigma^-$. The \textbf{negative Poincaré half-map} $\mathcal{P}^-: \Sigma^- \to \Sigma^+$ is defined as
	$$\mathcal{P}^-(\vxi) = \exp(\vA^- t^-(\vxi)) \vxi,$$
	where the first-return time is $t^-(\vxi) := \min \{ t> 0 \mid \ve_1\T \exp(t \vA^-) \vxi = 0 \}$. Similarly, for $\veta \in \Sigma^+$, the \textbf{positive Poincaré half-map} $\mathcal{P}^+: \Sigma^+ \to \Sigma^-$ is
	$$\mathcal{P}^+(\veta) = \exp(\vA^+ t^+(\veta)) \veta,$$
	where the first-return time is $t^+(\veta) := \min \{ t> 0 \mid \ve_1\T \exp(t \vA^+) \veta = 0 \}$.
\end{definition}

Throughout this paper, we adopt the following convention for classifying invariant cones. Without loss of generality, as the system is autonomous, we assume a trajectory generating a cone initiates in the crossing region $\Sigma^-$. A 'crossing' is then counted for each time the trajectory subsequently enters the positive half-space $\cVplus$ through $\Sigma^+$. The elementary Poincaré map $\mathcal{P} = \mathcal{P}^+ \circ \mathcal{P}^-$, which describes a trajectory's evolution through $\cVminus$ and then through $\cVplus$, therefore corresponds to a single crossing event. Consequently, the elementary invariant cone, which we are about to define, will also be referred to as a \textit{1-crossing invariant cone}. This convention provides a natural basis for classifying the multi-crossing structures introduced later.

A fundamental consequence of the system's homogeneity is that the crossing times $t^\pm$ are positive homogeneous of degree zero (i.e., $t^\pm(\beta\vxi) = t^\pm(\vxi)$ for $\beta>0$), making them constant on rays emanating from the origin. Consequently, the Poincaré half-maps are positive homogeneous of degree one, mapping rays to other rays. The elementary Poincaré map $\cP = \cP^+ \circ \cP^-$ for one full cycle is the composition of these half-maps and is thus also positively homogeneous of degree one, i.e., 
\begin{equation}\label{eq_homog_prop_poincare}
	\mathcal{P}^- (\beta \bm{\xi}) = \beta \mathcal{P}^- ( \bm{\xi}), \quad \mathcal{P}^+ (\beta \bm{\eta}) = \beta \mathcal{P}^+ ( \bm{\eta}), \quad 
	\mathcal{P} (\beta \bm{\xi}) = \beta \mathcal{P} ( \bm{\xi}), \quad \forall \beta >0.
\end{equation}

\begin{definition}[Elementary Invariant Cone]
	An \textbf{elementary invariant cone} is an invariant set of the flow of \eqref{eq:eom_firstorder} generated by an invariant ray on the switching hyperplane. It is determined by a non-zero vector $\overline{\vxi} \in \Sigma^-$ that is a solution to the nonlinear eigenvalue problem
	\begin{equation}\label{eq_def_invco_poincare}
		\mathcal{P}(\overline{\vxi}) = (\mathcal{P}^+ \circ \mathcal{P}^-)(\overline{\vxi}) = \mu \overline{\vxi},
	\end{equation}
	for some positive real scalar $\mu \in \mathbb{R}^+$, called the characteristic multiplier.
\end{definition}

\begin{figure}
	\centering	
	\includegraphics[scale=0.5]{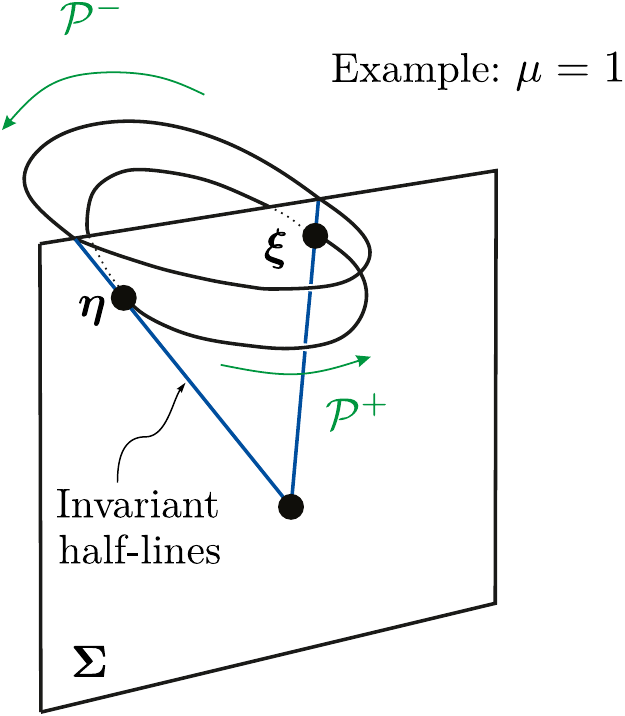}	
	\caption{A pictorial view of an invariant cone foliated by periodic orbits with $\mu=1$.}
	\label{fig_invcone_periodic}
\end{figure}

Such rays on $\Sigma$ which are invariant under $\cP$ are straight half-lines through the origin and will be referred to as invariant half-lines. The positive homogeneity property of the Poincaré map is essential for the following discussion and has an intuitive geometric interpretation. From \eqref{eq_homog_prop_poincare}, one can directly see that a point on an invariant half-line in $\Sigma$ that goes through the origin is mapped by $\mathcal{P}$ to a point on the same half-line. This is illustrated in Figure \ref{fig_invcone_periodic} for the case $\mu = 1$, where the invariant cone is foliated by periodic orbits. 

The multiplier $\mu$ governs the dynamics of trajectories \textit{on the cone}. If $\mu < 1$, trajectories spiral towards the origin; if $\mu > 1$, they spiral away; and if $\mu = 1$, the cone is foliated by a family of periodic orbits. The stability of the cone itself, i.e., its ability to attract nearby trajectories, is a separate issue.

Formally, since the Poincar\'e map $\mathcal{P}$ maps a subset of $\mathbb{R}^n$ (namely $\Sigma^- \subset \Sigma$) to itself, its Jacobian $D\mathcal{P} = \frac{\partial \mathcal{P}}{\partial \vxi}$ is an $n \times n$ matrix. Because the range of $\mathcal{P}$ is constrained to the $(n-1)$-dimensional subspace $\Sigma$, this Jacobian is rank-deficient and thus has at least one eigenvalue $\lambda_n = 0$. This eigenvalue, whose eigenvector $\ve_1$ is normal to $\Sigma$, is a trivial consequence of the map's geometry and has no relevance for the dynamics on the cone.

The dynamically relevant behavior is governed by the remaining $n-1$ eigenvalues, whose eigenvectors lie in the subspace $\Sigma$. Differentiating the homogeneity condition \eqref{eq_homog_prop_poincare} yields
\begin{equation}
	\frac{\mathrm{d}}{\mathrm{d}\beta} (\mathcal{P}(\beta \vxi)) \stackrel{!}{=} \frac{\mathrm{d}}{\mathrm{d}\beta} (\beta \mathcal{P}(\vxi)) = \mathcal{P}(\vxi), 
\end{equation}
and by use of the chain rule
\begin{equation}
	\frac{\mathrm{d}}{\mathrm{d}\beta} (\mathcal{P}(\beta \vxi)) = \frac{\partial \mathcal{P}}{\partial \vxi}\Bigr|_{\beta \vxi} \vxi \stackrel{!}{=} \mathcal{P}(\vxi),
\end{equation}
which shows that the generating vector $\overline{\vxi}$ is one of these $(n-1)$ dynamically relevant eigenvectors, with its corresponding eigenvalue being $\mu$. The cone is said to be attractive if the remaining $n-2$ eigenvalues $\lambda_j$ of this Jacobian satisfy the condition \cite{Weiss2012}
\begin{equation}\label{eq_attractivity}
	|\lambda_j| < \min(1,\mu), \quad j = 1, \cdots, n-2.
\end{equation}
Under this condition, the cone serves as a nonsmooth analogue of a stable center manifold.

To find such a cone numerically, one must solve a system of nonlinear equations for the unknowns $(\vxi, t^-, t^+, \mu)$. This system combines the eigenvalue problem with the crossing-time conditions and a normalization constraint \cite{Kuepper2008}:
\begin{equation}\label{eq_classic_inv_cone_problem}
	\vF(\vxi,t^-,t^+, \mu) = \begin{pmatrix}
		e^{t^+ \vA^+} e^{t^- \vA^-} \vxi - \mu \vxi \\
		\ve_1\T e^{t^- \vA^- } \vxi \\ 
		\ve_1\T \vxi \\ 
		\vxi\T \vxi - 1 
	\end{pmatrix}\stackrel{!}{=} \bm{0}.
\end{equation}
The formal analysis in the literature has largely focused on these elementary cones. This focus is a natural consequence of the significant analytical challenges involved in proving the existence, number, and stability of even these simplest invariant structures, which often require a cumbersome analysis of transcendental equations \cite{carmona2005bifurcation, Huan2016}. While fundamental, this framework does not explicitly address more complex oscillatory behaviors where trajectories may cross the switching hyperplane multiple times per cycle. Although such motions are known to occur in PWL systems, they have lacked a systematic classification within the context of invariant cone theory.

\noindent%

\subsection{Generalization to $k$-Crossing Invariant Cones}

The framework of elementary invariant cones effectively describes the simplest oscillatory behaviors. To address the limitation noted at the end of the previous subsection and to build a theory for more complex dynamics, we now generalize this concept to trajectories that cross the switching hyperplane multiple times before returning to their initial ray. The fundamental object for this generalization is not a single invariant ray, but a finite set of rays that are cyclically permuted by the elementary Poincaré map.

\begin{definition}[Generating Cycle of Rays]
	A set of $k$ non-zero vectors $\{\overline{\vxi}_1, \overline{\vxi}_2, \dots, \overline{\vxi}_k\}$ in $\Sigma$ that generate $k$ distinct rays\footnote{Two vectors $\vxi_i$ and $\vxi_j$ generate distinct rays if $\vxi_i \neq \beta \vxi_j$ for any $\beta > 0$.} constitutes a \textbf{generating cycle of rays of period k} under the map $\mathcal{P}$ if there exists a single positive scalar $\mu \in \mathbb{R}^+$, the characteristic multiplier, such that:
	\begin{align}
			\mathcal{P}(\overline{\vxi}_i) &= \overline{\vxi}_{i+1} \quad \text{for } i = 1, \dots, k-1 \label{eq:cycle_map} \\ 
			\mathcal{P}(\overline{\vxi}_k) &= \mu \overline{\vxi}_1 \label{eq:cycle_closure} 
	\end{align}
	Following the convention established previously, the cycle starts in the region leading to the negative half-space, $\overline{\vxi}_1 \in \Sigma^-$.
\end{definition}
A crucial property of such a cycle is that the characteristic multiplier $\mu$ is an invariant of the cycle itself, shared by all of its generating vectors. We state this formally as follows.

\begin{proposition}
	Let $\{\overline{\vxi}_1, \dots, \overline{\vxi}_k\}$ be the set of vectors defining a generating cycle of rays of period $k$ with characteristic multiplier $\mu$. Then, each vector $\overline{\vxi}_j$ in this set is an eigenvector of the k-th iterate map $\mathcal{P}_k := \mathcal{P}^k = P \circ \cdots \circ \cP$ ($k$-times) with the same eigenvalue $\mu$.
\end{proposition}

\textit{Proof.} The eigenvalue relation for the first vector, $\overline{\vxi}_1$, follows directly from iterating the relations in the definition:
\begin{equation*}
	\mathcal{P}_k(\overline{\vxi}_1) = \mathcal{P}^{k-1}(\mathcal{P}(\overline{\vxi}_1)) = \mathcal{P}^{k-1}(\overline{\vxi}_2) = \dots = \mathcal{P}(\overline{\vxi}_k) = \mu \overline{\vxi}_1.
\end{equation*}
To show that any other vector in the cycle, say $\overline{\vxi}_j$ for $j>1$, shares this same eigenvalue, we use its relation to $\overline{\vxi}_1$, which is $\overline{\vxi}_j = \mathcal{P}^{j-1}(\overline{\vxi}_1)$. We can then write:
\begin{align*}
	\mathcal{P}_k(\overline{\vxi}_j) &= \mathcal{P}^k(\mathcal{P}^{j-1}(\overline{\vxi}_1)) &&\text{by definition of } \overline{\vxi}_j \\
	&= \mathcal{P}^{j-1}(\mathcal{P}^k(\overline{\vxi}_1)) &&\text{since } \mathcal{P}^{j-1} \text{ and } \mathcal{P}^k \text{ commute} \\
	&= \mathcal{P}^{j-1}(\mu \overline{\vxi}_1) &&\text{using the eigenvalue relation for } \overline{\vxi}_1 \\
	&= \mu \mathcal{P}^{j-1}(\overline{\vxi}_1) &&\text{by homogeneity of } \mathcal{P} \\
	&= \mu \overline{\vxi}_j &&\text{by definition of } \overline{\vxi}_j
\end{align*}
This holds for all $j=1, \dots, k$, proving that all vectors in the generating cycle are eigenvectors of $\mathcal{P}_k$ with the same characteristic multiplier $\mu$. \hfill $\square$

This allows us to define the $k$-crossing invariant cone as the complete geometric object generated by this entire cycle of rays.

\begin{definition}[$k$-crossing Invariant Cone]
	For an integer $k \ge 1$, a \textbf{$k$-crossing invariant cone} is the full invariant set of the flow generated by a generating cycle of $k$ rays on the switching manifold. It is the union of $k$ cone segments, where each segment is the surface traced by the flow originating from one of the rays in the generating cycle. The entire composite structure is characterized by a single multiplier $\mu$.
\end{definition}

For $k=1$, this definition recovers the elementary (or 1-crossing) invariant cone, which is generated by a single invariant ray that is a fixed point of $\mathcal{P}$. For $k>1$, it describes a composite object. A trajectory starting on the ray of $\overline{\vxi}_1$ will land on the ray of $\overline{\vxi}_2$ after one elementary mapping cycle, then on the ray of $\overline{\vxi}_3$ after the next, and so on, returning to the initial ray of $\overline{\vxi}_1$ (scaled by $\mu$) only after $k$ elementary cycles have been completed.

A $k$-crossing invariant cone is generated by a cycle of $k$ distinct rays on the switching manifold. The transition from one ray in this cycle to the next is governed by the elementary Poincaré map $\mathcal{P}$, which itself consists of two distinct temporal segments: an excursion through the negative half-space $\cVminus$ and a subsequent excursion through the positive half-space $\cVplus$. Because the first-return times are state-dependent, each of these $2k$ segments will in general have a different duration.

To formalize this, let the generating cycle of rays be represented by the vectors $\{\overline{\vxi}_1, \overline{\vxi}_2, \dots, \overline{\vxi}_k\}$. The full sequence of events on the switching manifold is as follows:
\begin{enumerate}
	\item A trajectory starts on the first ray at $\overline{\vxi}_1 \in \Sigma^-$. It evolves in the negative half-space $\cVminus$ for a duration of $t_1^- = t^-(\overline{\vxi}_1)$, striking $\Sigma$ at an intermediate point $\overline{\veta}_1 = \mathcal{P}^-(\overline{\vxi}_1) = \exp(\vA^- t_1^-)\overline{\vxi}_1$. By definition, $\overline{\veta}_1 \in \Sigma^+$. From there, it evolves in $\cVplus$ for a duration of $t_1^+ = t^+(\overline{\veta}_1)$, returning to $\Sigma$ at the second generating ray, $\overline{\vxi}_2 = \mathcal{P}^+(\overline{\veta}_1)$.
	\item This process repeats, starting from $\overline{\vxi}_2$. The trajectory evolves for a new duration $t_2^- = t^-(\overline{\vxi}_2)$, reaches an intermediate point $\overline{\veta}_2$, evolves for a time $t_2^+ = t^+(\overline{\veta}_2)$, and lands on the third ray, $\overline{\vxi}_3$.
	\item The sequence continues for $k$ full transitions. The final transition, from $\overline{\vxi}_k$, completes the cycle by mapping back to the first ray, scaled by the characteristic multiplier $\mu$.
\end{enumerate}
This process involves a total of $2k$ possibly distinct, state-dependent first-return times, $\{t_1^-, t_1^+, \dots, t_k^-, t_k^+\}$. In the special case of a cone foliated by periodic orbits ($\mu=1$), the total period $T$ is the sum of all these times: $T = \sum_{j=1}^k (t_j^- + t_j^+)$. This correspondence is summarized in Table \ref{tab:crossing_times_k}.

\begin{table}[h!]
	\centering
	\caption{Temporal Structure of a $k$-crossing Invariant Cone}
	\label{tab:crossing_times_k}
	\begin{tabular}{c|c|c|c|c}
		\hline
		\textbf{Transition} & \textbf{Starting Ray} & \textbf{Negative Time ($t_j^-$)} & \textbf{Positive Time ($t_j^+$)} & \textbf{Ending Ray} \\
		\hline
		1 & $\overline{\vxi}_1$ & $t^-(\overline{\vxi}_1)$ & $t^+(\mathcal{P}^-(\overline{\vxi}_1))$ & $\overline{\vxi}_2$ \\
		2 & $\overline{\vxi}_2$ & $t^-(\overline{\vxi}_2)$ & $t^+(\mathcal{P}^-(\overline{\vxi}_2))$ & $\overline{\vxi}_3$ \\
		$\vdots$ & $\vdots$ & $\vdots$ & $\vdots$ & $\vdots$ \\
		k & $\overline{\vxi}_k$ & $t^-(\overline{\vxi}_k)$ & $t^+(\mathcal{P}^-(\overline{\vxi}_k))$ & $\mu \overline{\vxi}_1$ \\
		\hline
	\end{tabular}
\end{table}

To find such a structure, it is not necessary to treat all $k$ generating vectors as independent unknowns. The entire generating cycle is determined by the choice of a single starting ray, $\overline{\vxi}_1$, and the sequence of $2k$ crossing times. This leads to a compact formulation of the root-finding problem with a minimal set of primary unknowns: $(\overline{\vxi}_1, t_1^-, t_1^+, \dots, t_k^-, t_k^+, \mu)$, for a total of $n + 2k + 1$ scalar variables.

The problem is to find a vector of unknowns $\mathbf{X} = (\overline{\vxi}_1, t_1^-, \dots, t_k^+, \mu)$ that solves the system $\vF(\mathbf{X}) = \bm{0}$, where $\vF$ is constructed from $n + 2k + 1$ independent scalar equations:
\begin{itemize}
	\item \textbf{The Closure Condition ($n$ equations):} The vector equation ensuring the cycle closes on the initial ray.
	\begin{equation*}
		\mathcal{P}_k(\overline{\vxi}_1; t_1^-, \dots, t_k^+) - \mu\overline{\vxi}_1 = \bm{0},
	\end{equation*}
	where $\mathcal{P}_k(\cdot)$ is the full composition $\exp(t_k^+\vA^+)\exp(t_k^-\vA^-) \cdots \exp(t_1^+\vA^+)\exp(t_1^-\vA^-)$.
	
	\item \textbf{The Independent Crossing Conditions ($2k-1$ equations):} The scalar equations that implicitly define the first $2k-1$ crossing times.
	\begin{align*}
		\ve_1\T \exp(t_j^- \vA^-) \overline{\vxi}_j &= 0, \quad \text{for } j=1, \dots, k \\
		\ve_1\T \exp(t_j^+ \vA^+) \overline{\veta}_j &= 0, \quad \text{for } j=1, \dots, k-1,
	\end{align*}
	where the vectors $\overline{\vxi}_j$ and $\overline{\veta}_j$ are recursively defined from $\overline{\vxi}_1$ and the preceding times. The final crossing condition for $t_k^+$ is redundant and thus omitted.
	
	\item \textbf{Initial Vector Constraints (2 equations):}
	\begin{equation*}
		\ve_1\T \overline{\vxi}_1 = 0 \quad \text{and} \quad \overline{\vxi}_1\T \overline{\vxi}_1 - 1 = 0.
	\end{equation*}
\end{itemize}
The total number of scalar equations is $n + (2k-1) + 2 = n + 2k + 1$, which perfectly matches the number of unknowns. The system is therefore a complete set of nonlinear equations, and its solution characterizes the $k$-crossing invariant cone.

\section{Graph-Style Parametrization of Invariant Sets: From Smooth Manifolds to Cones}\label{section:graphstyle}
This section presents an extension of the classical graph-style parametrization of NNMs for computing the invariant cones of PWL systems. The development is organized into two main parts. First, Section \ref{Subsection_invManifApproach} reviews the foundational concepts of the classical method, wherein NNMs are interpreted as invariant manifolds in the phase space—a concept pioneered by Shaw and Pierre \cite{Shaw1991, Shaw1993}. The summary of the established theory provided therein is based on the comprehensive reviews by Shaw \cite{Shaw2014} and Touzé and Vizzaccaro \cite{Touze2024}. Building on this foundation, Section 4.2 introduces our novel extension of the framework to homogeneous PWL systems and discusses its limitation and failure modes.
\subsection{The Invariant Manifold Approach for NNMs in Smooth Systems}\label{Subsection_invManifApproach}
Linear vibration analysis relies on the eigenspaces spanned by LNMs, which are orthogonal hyperplanes in phase space sharing the invariance property, i.e., trajectories starting on the LNM remain on it indefinitely. The orthogonality is a direct consequence of the uncoupling in the modal space. Each LNM is a two-dimensional invariant plane filled with periodic orbits of the conservative linear dynamics or decaying orbits spiraling to the stable equilibrium in the case of a dissipative linear system with modal damping. Extending this geometric picture beyond infinitesimal amplitudes is essential for nonlinear systems, where energy exchange, bifurcations and frequency-amplitude dependence are all characteristic features of nonlinear behavior and are related to the curvature of the invariant manifolds that replace the planar LNMs. Based on this argument, a NNM of an $N$-degree-of-freedom (DOF) mechanical system can therefore be defined as a smooth two-dimensional invariant manifold $\mathcal{W} \subset \mR^{2N}$ that is tangent at the origin to a selected modal eigenspace and persists for finite oscillation amplitudes \cite{Shaw1991}. For undamped mechanical structures, the definition coincides with the Lyapunov sub-centre manifold (LSM), while extending the earlier NNM definition of Rosenberg \cite{Rosenberg1962} to dissipative settings. Even though several key properties of LNMs are lost in the extension process, the crucial invariance property is inherited by the NNMs. This leads to NNMs playing a central role in model order reduction, since they provide a rigorous basis for a mathematically justifiable reduction of the dynamics as opposed to linear MOR approaches, where linear subspaces which are not invariant to the full nonlinear dynamics are used as a reduction basis. Thanks to the invariance property, an orbit initiated on $\mathcal{W}$ never leaves it; the reduced dynamics on $\mathcal{W}$ therefore reproduces the backbone curve and frequency–energy dependence of the full system with a minimal two-degree-of-freedom ROM.
In the following, we will describe how the invariant manifold - more specifically its parametrization $\mathcal{W}$ - can be obtained and used to derive a ROM. To keep notation simple, we will restrict the next developments to the autonomous setting, but a similar approach can be applied in the presence of explicitly time-dependent forcing \cite{Jiang2005,Shaw2014}. Let's consider an $N$-DOF system with its autonomous dynamics written in first-order form as
\begin{equation}
\begin{aligned}
\dot{x}_i &= y_i, \\
\dot{y}_i &= f_i(\bar{\vx}, \bar{\vy}), \quad i = 1,\ldots, N
\end{aligned}\label{eq:eom}
\end{equation}
where displacements and velocities are grouped in $\bar{\vx}=(x_1, x_2, \ldots, x_N)\T$, and $\bar{\vy} = (y_1,y_2,\ldots,y_N)\T$, respectively, and $f_i$ is a smooth nonlinear function. The core idea is to express the system's response on the NNM as being governed by a pair of master coordinates, typically the displacement and velocity of a single chosen mode, say $(x_m, y_m)$, where $\dot{x}_m=y_m$. All other slave coordinates $(x_j,y_j)$ for $j\neq m$, are assumed to be functions of these master coordinates, which explains the choice of the name as \textit{graph style} parametrization, i.e., the NNM is sought as a graph: 
\begin{equation}
x_{j}=X_{j}(x_{m},y_{m}),\qquad
y_{j}=Y_{j}(x_{m},y_{m}),\qquad \forall j = 1, \ldots n, \quad \mbox{with } j \neq m,
\label{eq:graph}
\end{equation}
These relations define the NNM as a two-dimensional, generally curved surface (invariant manifold) within the full $2N-$dimensional phase space. This graph formulation paves the way for a systematic computation, as detailed next. For this surface to be invariant, any trajectory starting on it must remain on it for all time. This imposes a mathematical constraint. By applying the chain rule to take the time derivative of the parametrization equations \eqref{eq:graph} and substituting the system's equations of motion \eqref{eq:eom}, we eliminate the explicit time dependence and arrive at the so-called \textit{invariance equations}:
\begin{equation}\label{eq:inv}
	\begin{aligned}
		\frac{\partial X_{j}} {\partial x_{m}} \, y_{m}
		+\frac{\partial X_{j}}{\partial y_{m}} f_{m}(\bar{\vx},\bar{\vy})
		&=Y_{j}(x_{m},y_{m}), \\  
		\frac{\partial Y_{j}}{\partial x_{m}} y_{m}
		+\frac{\partial Y_{j}}{\partial y_{m}} f_{m}(\bar{\vx},\bar{\vy})
		&= f_{j}(\bar{\vx},\bar{\vy}),
	\end{aligned}
\end{equation}
for every slave index $j \neq m$. Equations \eqref{eq:inv} constitute $2N-2$
coupled first-order PDEs in the free variables $(x_m,y_m)$, whose solution describes the geometry of the NNM. Since solving these nonlinear PDEs exactly is generally intractable in closed form, local solutions around the equilibrium point are sought using power-series expansions in the master coordinates $(x_m, y_m)$. For instance, a third-order truncation is taken in the original work of Shaw and Pierre \cite{Shaw1993} as
\begin{equation}
	\begin{aligned}
	x_{j}&= a_{j1} x_m + a_{j2} y_m + a_{j3}x_{m}^{2}+a_{j4}x_{m}y_{m}+a_{j5}y_{m}^{2}+a_{j6}x_{m}^{3}+\cdots, \\ 
	y_{j}&=b_{j1} x_m + b_{j2} y_m + b_{j3}x_{m}^{2}+b_{j4}x_{m}y_{m}+b_{j5}y_{m}^{2}+b_{j6}x_{m}^{3}+\cdots \, .
	\end{aligned}\label{eq:polyexpansion}
\end{equation}
Substituting these expansions into the invariance equations \eqref{eq:inv} and collecting monomials of same order provides a sequence of algebraic equations to solve for the unknown coefficients at each order (linear, quadratic, cubic, etc.). This procedure systematically constructs a polynomial approximation of the invariant manifold. 
Once the manifold functions $(X_j,Y_j)$ are approximated, they are substituted back into the equation of motion for the master coordinate. This yields a planar ordinary differential equation that governs the dynamics on the manifold, i.e. the reduced-order model:
\begin{equation}
\begin{aligned}
\dot{x}_m &= y_m \\ 
\dot{y}_m &= f_m(X_1(x_m,y_m), \ldots, x_m, \ldots,X_N(x_m,y_m), Y_1(x_m,y_m), \ldots, y_m, \ldots,Y_N(x_m,y_m) ).
\end{aligned}
\end{equation}
In subsequent developments, numerical methods were proposed to solve the invariance equations PDEs \eqref{eq:inv}, as in Pesheck et al. \cite{Pesheck2001a, Pesheck2001, Pesheck2002}, Blanc et al. \cite{Blanc2013}, and Renson et al. \cite{Renson2014b}. Furthermore, Jiang et al. also extended this approach to harmonically forced vibrating systems, incorporating time-dependent invariant manifolds \cite{Jiang2005}. In this case, a nonautonomous term $f_t \cos(\Omega t)$ is added to the dynamics \eqref{eq:eom}. The phase $\varphi= \Omega t$ is then appended as an additional state, $\dot{\varphi}=\Omega$, thereby recovering an autonomous system with an extra variable. Equations \eqref{eq:inv} then extend to a system parametrized by $(x_m,y_m, \varphi)$ in which $\varphi$ plays the role of a cyclic coordinate; the resulting manifold yields an invariant forced response surface that collapses to the unforced NNM when no harmonic forcing is applied. Furthermore, Chen and Shaw \cite{Chen1996} and Jiang et al. \cite{Jiang2004} extended this graph style invariant manifold approach to systems with piecewise-linear restoring forces with a nonvanishing clearance gap, i.e., PWA systems. This restriction ensures that the system is \textit{linearizable} in the sense that the linear subspaces required for the NNM continuation are well defined. In the present study, we will focus on the case of \textit{nonlinearizable} PWL systems, i.e., with a vanishing clearance gap. In this setting, the Jacobian jumps exactly at the equilibrium, where no properly defined eigenvalues exist. This leads to the absence of any underlying linear subspace to anchor a nonlinear invariant manifold and therefore calls for the conic generalization introduced next. 

\subsection{Graph style invariant cone parametrization}
In this section, we will adapt the classical graph-style invariant-manifold method to homogeneous PWL systems in the form \eqref{eq:eom_firstorder}. We will show how the positive homogeneity property collapses the PDE-type invariance equations to an ODE on a generating line. Similar to the classical procedure, we start by choosing a pair of master coordinates $(x_m,y_m)$, which is then transformed to a polar form:
\begin{equation}\label{eq:master_polartrafo}
\begin{pmatrix}
x_m \\ y_m 
\end{pmatrix} = \begin{pmatrix}
u \\ v 
\end{pmatrix} = r \begin{pmatrix}
\cos(\theta) \\ \sin(\theta)
\end{pmatrix}, \quad \mbox{with } \quad r = \sqrt{u^2+v^2} > 0.
\end{equation}

The dynamics of the master pair can be given in terms of $(u,v)$ as $\dot{u}=v$ and $\dot{v}=f_m$, or equivalently using the new polar coordinates $(r,\theta)$ as follows: 
\begin{equation}\label{eq_ROM_in_polarCoordinates}
\begin{aligned} \dot{r} &= r \sin(\theta) \cos(\theta) + \sin(\theta) f_m(\vp), \\ \dot{\theta} &= - \sin(\theta)^2 + \frac{\cos(\theta)}{r} f_m(\vp), \end{aligned}
\end{equation}
where $\vp=(u, v, \ldots, x_j, y_j, \ldots )\T, \quad \forall j \neq m  $.

Similar to the smooth case, we assume that the rest of the state variables are functionally dependent on the master pair, i.e., the NNM is sought as a graph
\begin{equation}
	\begin{aligned}
	x_j &= X_j(r,\theta) \\ 
	y_j &= Y_j(r,\theta), \quad \forall j \neq m.
	\end{aligned}
\end{equation} 
Applying the chain rule and eliminating the explicit time dependence in the slave dynamics yields the invariance equations expressed in polar coordinates as
\begin{equation}
\begin{aligned}\label{eq:PDEs1_augmentedstate}
\frac{\partial X_{j}} {\partial r} \, \dot{r}
+\frac{\partial X_{j}}{\partial \theta} \dot{\theta}
&=Y_{j}(r,\theta), \\  
\frac{\partial Y_{j}}{\partial r} \dot{r}
+\frac{\partial Y_{j}}{\partial \theta } \dot{\theta}
&= f_{j}(\vp(r,\theta)), \quad \forall j \neq m,
\end{aligned}
\end{equation}
where the state vector is now expressed in terms of the polar master coordinates as 
\begin{equation}
\vp(r,\theta)=(r \cos(\theta), r\sin(\theta), \ldots, X_j(r,\theta), Y_j(r,\theta), \ldots )\T, \quad \forall j \neq m. 
\end{equation}
Due to the lack of any underlying linear subspace tangent to the NNM at the equilibrium, we cannot take a polynomial series expansion as in \eqref{eq:polyexpansion} to solve the PDEs \eqref{eq:PDEs1_augmentedstate}. However, we can leverage the positive homogeneity property of the dynamics \eqref{eq:eom_firstorder} to compensate for this lack of smoothness. As discussed in the previous section, homogeneous PWL systems admit invariant cones, which we denote here by $\mathcal{C}$. The fundamental geometric property of an invariant cone $\mathcal{C}$ is that it consists of the origin and a collection of half-lines, or \textit{rays}, passing through the origin. This means that if a point $\vp$ belongs to the cone, then the ray $R(\vp) := \{ \alpha \vp \, \mid \, \alpha > 0\}$ through that point also belongs to the cone. Consequently, the cone's geometry is not defined by the magnitude of its points, but purely by the set of \textit{directions} it encompasses.

The concept of directionality of the cone is formalized by considering the \textit{Ray Set} of the cone, denoted $\mathcal{R}(\mathcal{C})$ and defined as
\begin{equation}
	\mathcal{R}(\mathcal{C}) := \{  R(\tilde{\vp}), \, \forall \tilde{\vp} \in \mathcal{C}, \, \tilde{\vp} = \frac{\vp}{||\vp||}, \, \vp \neq \bm{0} \},
\end{equation}
where we used the normalized point $\tilde{\vp}$ to highlight the independence of the rays on the magnitude of their points. 
Each element of $\mathcal{R}(\mathcal{C})$ is a single ray, representing one of the cone's directions, with the origin $\bm{0}$ excluded from being a direction per se. Formally, $\mathcal{R}(\mathcal{C})$ is the set of equivalence classes on $\mathcal{C} \setminus \{0\}$, where two points are equivalent if they are positive scalar multiples of each other.

To parameterize the cone, we must select a unique representative point for each ray. This is achieved by intersecting the cone with a suitable hypersurface, whose choice is not uniquely defined. Since our analysis is centered on the master pair $(u,v)$, we strategically choose a hypersurface defined in terms of these master coordinates; the cylinder $\mathcal{H}_{\textrm{cyl}} = \{ \vp \in \mathbb{R}^{2N} \, | \, u^2+v^2 = 1 \}$, which has a unit-circle projection in the master space. The intersection defines the \textit{generating set} $\mathcal{G} = \mathcal{C} \cap \mathcal{H}_{\text{cyl}}$. This procedure establishes a formal \textit{Ray-to-Generator Map}, $\Phi$, which assigns to each ray its intersection point on the cylinder:
\begin{equation}
	\begin{aligned}
		\Phi: \mathcal{R}(\mathcal{C}) &\to \mathcal{G} \\
		R(\vp) &\mapsto R(\vp) \cap \mathcal{H}_{\text{cyl}}.
	\end{aligned}
\end{equation}
For our method to be well-posed, this map $\Phi$ must be a bijection. The cases where it fails will be detailed Section 4.2.2.

\begin{figure}
	\centering	
	\includegraphics[scale=0.6]{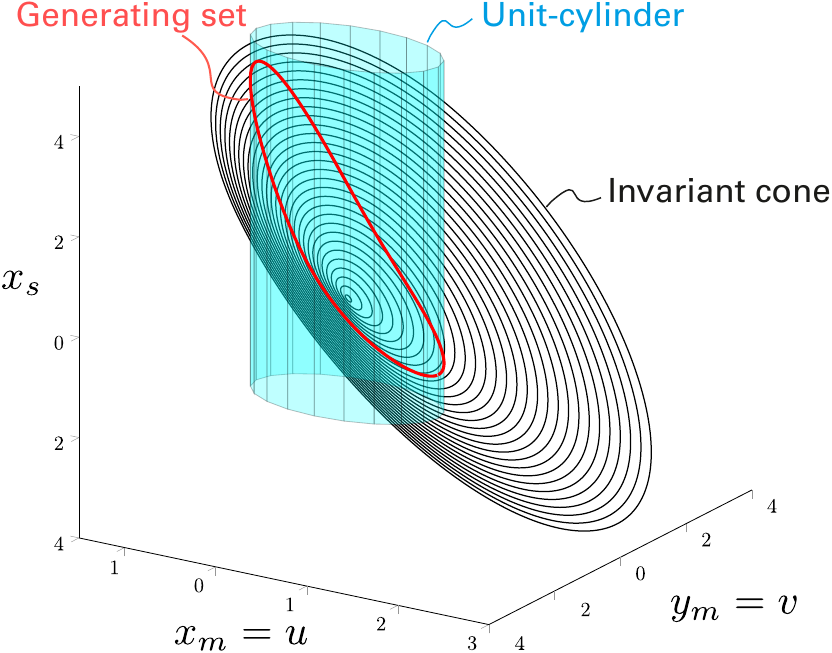}	
	\caption{A $3-$D representative sketch of an invariant cone and its generating set obtained as its intersection with the unit-cylinder.}
	\label{fig_GeneratingSet}
\end{figure}

Assuming for now that $\Phi$ is a well-defined bijection, the entire cone can be reconstructed by scaling every point in the generating set by all possible non-negative lengths:
\begin{equation}
	\mathcal{C} = \{ \gamma \vg \, | \, \vg \in \mathcal{G}, \, \gamma \ge 0 \}, \quad \text{with } \quad \mathcal{G} = \mathcal{C} \cap \mathcal{H}_{\text{cyl}}.
\end{equation}
The specific choice of the cylinder $\mathcal{H}_{\textrm{cyl}}$ has two powerful consequences. First, it gives a concrete physical meaning to the general scaling factor $\gamma$ needed to reconstruct the cone. For any point $\vp=(u,v,\vw) \in \mathcal{C}$, its corresponding point on the generating set is $\vg= (1/\gamma) \vp$. For $\vg$ to lie on the cylinder $\mathcal{H}_{\textrm{cyl}}$, its master coordinates must satisfy $(u/\gamma)^2 + (v/\gamma)^2=1$, which uniquely determines the scaling factor as $\gamma = \sqrt{u^2+v^2}$. This is precisely the radial coordinate $r$ of the master space. Thus, our choice of the cylinder forces the abstract scaling factor $\gamma$ to become the physically meaningful master coordinate $r$. Second, since the generating set $\mathcal{G}$ lies on a surface whose projection onto the master plane is the unit circle $S^1$, it is natural to parametrize $\mathcal{G}$ using the angular coordinate $\theta$. This transforms the abstract set $\mathcal{G}$ into a function $\mathcal{G}(\theta)$ that describes the generating curve obtained as the intersection of the cone $\mathcal{C}$ with the cylinder $\mathcal{H}_{\textrm{cyl}}$. 
This process effectively separates the cone's geometry into its directional shape, captured by $\mathcal{G}(\theta)$, and its radial extent, captured by the scaling factor $\gamma=r$, obtained as the magnitude of the projection onto the master plane $(u,v)$ through the special choice of $\mathcal{H}_{\textrm{cyl}}$ as intersecting hypersurface.
Assuming for now that the parametrization is valid, we can describe any point $\vp \in \mathcal{C}$ by scaling a point on the generating set. We further assume that $\mathcal{G}(\theta)$ is a single-valued function. The possible failure of this assumption and its implications are discussed Section 4.2.2. This leads to the conic \textit{Ansatz} for the slave variables:
\begin{equation}
	\begin{aligned} \label{eq:conicAnsatz}
		x_j &= r P_j(\theta) \\ 
		y_j &= r Q_j(\theta), \quad \forall j \neq m,
	\end{aligned}
\end{equation}
where $P_j(\theta)$ and $Q_j(\theta)$ are the unknown functions describing the shape of the generating curve. The full parametrization of a point on the cone is thus:
\begin{equation}\label{eq:cone_param}
	\begin{aligned} \mathcal{W} : [0, \infty) \times S^1 &\mapsto \mathbb{R}^{n} , \\ (r,\theta) &\mapsto \vp := \mathcal{W}(r,\theta) = r \mathcal{G}(\theta) = r \begin{pmatrix} \cos(\theta) \\ \sin(\theta) \\  \vdots \\ P_j(\theta) \\ Q_j(\theta) \\ \vdots \end{pmatrix} , \quad \mbox{with } \quad \mathcal{G}(\theta) = \mathcal{C} \cap \mathcal{H}_{\textrm{cyl}} . \end{aligned}
\end{equation}
This parametrization allows us to rewrite the polar master dynamics \eqref{eq_ROM_in_polarCoordinates} in a simpler form by leveraging the positive homogeneity property of the dynamics $(f_m(r\mathcal{G}(\theta))=r f_m(\mathcal{G}(\theta)))$, yielding:
\begin{equation}\label{eq_ROM_in_polarCoordinates_without_r}
	\begin{aligned} \tilde{R}(\theta) := \frac{\dot{r}}{r} &= \sin(\theta) \cos(\theta) + \sin(\theta) f_m(\mathcal{G}(\theta)), \\ \tilde{\Theta}(\theta) := \dot{\theta} &= - \sin(\theta)^2 + \cos(\theta) f_m(\mathcal{G}(\theta)), \end{aligned}
\end{equation}
To sum up, the parametrization of the generating set $\mathcal{G}(\theta)$ defines the set of all points of the cone having a unit-circle projection in the master space, which can then be scaled by the positive radius~$r$ in any angular direction $\theta \in [0, 2\pi]$. A three-dimensional sketch representing this idea is shown in Figure \ref{fig_GeneratingSet}.
Substituting the conic Ansatz \eqref{eq:conicAnsatz} into the slave dynamics and leveraging the positive homogeneity of the vector field ($f_j(r\vp) = r f_j(\vp)$) allows us to divide out by $r>0$, yielding the invariance equations as
\begin{equation}\label{eq:geometryODE}
\begin{aligned} P_j'(\theta) &= (Q_j(\theta) - \frac{\dot{r}}{r} P_j(\theta) ) \frac{1}{\dot{\theta}} &&=  (Q_j(\theta) - \tilde{R}(\theta) P_j(\theta) ) \tilde{\Theta}^{-1}(\theta)  \\ Q_j'(\theta) &= (\frac{1}{r}f_j(r \mathcal{G}(\theta)) - \frac{\dot{r}}{r} Q_j(\theta) ) \frac{1}{\dot{\theta}} &&= (f_j(\mathcal{G}(\theta)) - \tilde{R}(\theta) Q_j(\theta) )\tilde{\Theta}^{-1}(\theta) ,
\end{aligned} 
\end{equation}
which we will refer to as the geometry ODE of the invariant cone.
This represents a dramatic simplification of the invariance equations compared with the smooth nonlinear case. Instead of a set of PDEs, we end up with a set of ODEs describing the geometry of the generating set $\mathcal{G}(\theta) = \mathcal{C} \cap \mathcal{H}_{\textrm{cyl}}$ belonging to the invariant cone. 
The ODE \eqref{eq:geometryODE} is accompanied by periodic boundary conditions $P_j(0)=P_j(2\pi)$ and $Q_j(0)=Q_j(2\pi)$, as the generating set (if $\Phi$ is bijective) forms a closed orbit on the unit cylinder. The periodic boundary conditions motivate a truncated Fourier series expansion as approximation for the unknown functions $P_j$ and $Q_j$:
\begin{equation}
\begin{aligned}\label{eq:FourierAnsatz}
P_j(\theta) &= a_{j,0} + \sum_{k=1}^{N_h} \left[ a_{j,k} \cos(k\theta) + b_{j,k} \sin(k\theta) \right], \\ 
Q_j(\theta) &= c_{j,0} + \sum_{k=1}^{N_h} \left[ c_{j,k} \cos(k\theta) + d_{j,k} \sin(k\theta) \right],
\end{aligned}
\end{equation}
where $N_h$ is the harmonic order of truncation. 
Substituting the Fourier expansion \eqref{eq:FourierAnsatz} into the geometry ODE \eqref{eq:geometryODE} and applying a harmonic balance method yields a system of nonlinear algebraic equations for the unknown Fourier coefficients in \eqref{eq:FourierAnsatz}, which can be solved using a Newton-type root finding scheme. Once the Fourier coefficients are obtained, the harmonic terms $\cos(k\theta)$ and $\sin(k\theta)$ in \eqref{eq:FourierAnsatz} can be reformulated as functions of the physical master coordinates $(u,v)$ using the Chebyshev polynomials of the first and second kind, which are defined by the following relationships:
\begin{equation}
	\begin{aligned}
	\cos(k\theta) &= T_k\bigl(\cos\theta\bigr) &&= T_k\bigl(\frac{u}{r}\bigr) \\ 
	\sin(k\theta) &= \sin\theta \; U_{k-1}\bigl(\cos\theta\bigr) &&= \frac{v}{r} U_{k-1}\bigl(\frac{u}{r}\bigr) \\ \mbox{with } r &= \sqrt{u^2 + v^2 }
	\end{aligned}
\end{equation}
Substituting back into the Fourier expansion yields: 
\begin{equation}\label{eq:PjQj_in_uv} 
\begin{aligned}
	P_j(u,v) = a_{c,0} + \sum_{k=1}^{N_h} \left[ a_{j,k} T_k\bigl(\frac{u}{r}\bigr) + b_{j,k} \frac{v}{r} U_{k-1}\bigl(\frac{u}{r}\bigr)  \right], \\ 	
	Q_j(u,v) = c_{c,0} + \sum_{k=1}^{N_h} \left[ c_{j,k} T_k\bigl(\frac{u}{r}\bigr) + d_{j,k} \frac{v}{r} U_{k-1}\bigl(\frac{u}{r}\bigr)  \right] .
\end{aligned}
\end{equation}
Therefore, we end up with an explicit parametrization of the cone in terms of the master coordinates $(u,v)$ as follows: 
\begin{equation}\label{eq:Slaves_in_uv}
	\begin{aligned} x_j &= \sqrt{u^2+v^2} P_j(u,v) \\ 
	y_j &= \sqrt{u^2+v^2} Q_j(u,v)  . \end{aligned}
\end{equation}
These functions describe the dependence of the slave variables on the master pair $(u,v)$ explicitly, and when plugged into the master dynamics lead to the following $1-$DOF ROM: 
\begin{equation}
	\begin{aligned} \dot{u} &= v \\ \dot{v} &= f_m(\tilde{\mathcal{W}}(u,v)),   \end{aligned}\label{eq:ROM_uv}
\end{equation}
where $\tilde{\mathcal{W}}(u,v)$ denotes the full state vector parameterized by $(u,v)$ obtained using the Chebyshev polynomials as demonstrated in \eqref{eq:Slaves_in_uv} and \eqref{eq:PjQj_in_uv}. Another form of the ROM is given by \eqref{eq_ROM_in_polarCoordinates} which expresses the reduced dynamics in terms of the polar coordinates $(r,\theta)$, wherein the state vector $\vp = \mathcal{W}(r,\theta) = r \mathcal{G}(\theta)$ is known after solving for the Fourier coefficients of \eqref{eq:FourierAnsatz}. 

To summarize, we have derived a parametrization method for invariant cones of PWL systems, very similar in spirit to the classical graph-style parametrization method for smooth invariant manifolds. By leveraging the positive homogeneity property of the cone and choosing a cylinder with a unit-circle projection on the master plane as an intersecting hypersurface to characterize the generating set of the cone, we were able to obtain a closed-form parametrization of the invariant cone, and thus derive a mathematically justifiable ROM explicitly. This ROM represents the system dynamics constrained to the invariant cone, i.e., the NNM of the PWL system. As is well-known for smooth nonlinear systems, the graph-style parametrization method is based on a graph assumption and breaks down in several scenarios where the invariant manifold possesses a complex geometry. In Section 4.2.2, we will focus on the potential problems that our cylinder-based graph-style parametrization of the invariant cone could encounter. 

\subsubsection{Choice of the Master Coordinates}
The core idea behind the extension of LNMs to NNMs, is that during an NNM motion, there is a functional dependence between all degrees of freedom. This motion can thus be parametrized by a single displacement-velocity pair, where the choice of the master coordinates $(u,v) = (x_m, v_m)$, i.e., the nonlinear modal coordinates, is a priori arbitrary \cite{Renson2013}. This holds, obviously, as long as the master pair is not a node of the sought NNM, in which case the graph-assumption would collapse from the get-go. Unfortunately, general good practices to identify the most suitable set of master coordinates \cite{Shaw2014}. Due to the special nature of our PWL systems and the special structure of the geometry ODEs obtained from the invariance equations in \eqref{eq:geometryODE}, we should highlight a remarkable point regarding the choice of the master pair. The class of PWL systems with two zones governed by dynamics of the form \eqref{eq:eom_firstorder} and \eqref{eq:sys_matrices} is characterized by a localized nonsmooth behavior, i.e., the nonsmoothness induced by the switch affects a restricted set of DOFs upon which the discontinuous support may act, while the dynamics of the remaining state variables is purely linear. For instance, if the elastic contact is placed at a free end of an oscillator, only one DOF is affected by the PWL discontinuity, and thus, the nonsmoothness affects the time evolution of one single velocity belonging to that DOF governing the switch. Taking the displacement and velocity of that specific switching DOF as master coordinates leaves all the slave dynamics $f_j,\, \forall j \neq m$ to be purely linear. In this case, the geometry ODE \eqref{eq:geometryODE} is linear in the unknowns $P_j,Q_j$, where the rest of the terms are functions of $\theta$. This choice is not absolutely necessary but hugely simplifies solving the ODE for a periodic (closed) generating curve using HBM. If this geometry ODE is linear, the HBM is guaranteed to converge to a unique solution, and thus the numerical computation of the unknown Fourier coefficients of the slaves $P_j,Q_j$ becomes simpler compared with the case of a PWL ODE obtained if the switching DOF is handled as a slave variable. 

\subsubsection{Conditions for Validity and Potential Failure Modes}\label{sec:cone_limitations}
The validity of the proposed graph-style parametrization hinges on two critical assumptions: first, that the ray-to-generator map $\Phi$ is a total bijection, and second, that the generating set $\cG$ can be represented as a single-valued function of the angular coordinate $\theta$. A violation of these assumptions leads to a breakdown of the method, which can manifest in three distinct ways. 

First, the graph-style parametrization fails if a ray on the cone is parallel to the cylinder's axis, preventing it from intersecting the cylinder. In this case, the projection $\pi$ of a point $\vp_0$ on this ray onto the master space is zero, i.e., $\exists \vp_0 \in \mathcal{C} \setminus \{0\} = \mathcal{R}(\mathcal{C})$ such that $\pi(\vp_0) = (u_0, v_0) = (0,0)$. The map $\Phi$ thus becomes a partial function (i.e., a function that is not defined for all elements in its domain of definition), leading to what we will refer to as a nodal failure. Physically, this corresponds to a scenario where the chosen master coordinates are at a nodal point of the NNM, resulting in a ROM that cannot capture dynamics occurring purely in the slave coordinate space. 

Second, the method becomes ill-posed if the generating curve~$\cG$ is not a single-valued function of~$\theta$. This occurs when two distinct points on the generating set, $\vg_1, \vg_2 \in \mathcal{G}$, share the same angular coordinate $\theta_0$, meaning the projection $\pi|_{\mathcal{G}}: \mathcal{G} \to S^1$ is not injective. We term this a multi-valued failure. Such a multi-valued geometry signals that the angle $\theta$ is not a suitable parameter to describe the generating set and makes the geometry ODE \eqref{eq:geometryODE} ill-posed as it does not have unique solutions. This failure may simply indicate that the master coordinates are not suitably chosen, or it may characterize a true cone geometry with multiple "sheets" over the same master-space direction. We will see in Section 7 that this case may correspond to a NNM at internal resonance, whose motion takes place on a multiple-crossing invariant cone as defined in Section 3. 

Third, a local breakdown of the parametrization occurs if the angular velocity vanishes, i.e., $\dot{\theta}=0$, which renders the geometry ODE \eqref{eq:geometryODE} singular. This conditions signifies that at a point $\vp$ on the cone, the projection of the velocity vector $\dot{\vp}$ onto the master plane is purely radial. Geometrically, this corresponds to a turning point in the generating curve $\mathcal{G}(\theta)$, a situation we will call a fold-point failure. This local failure of $\theta$ as a parameter for $\cG$ indicates a complex cone geometry, due to a loop in the generating curve causing a fold in the cone surface. In Section 5, a numerical example is shown where the invariant cone folds due to a strong nonsmooth effect, leading to the breakdown of the graph-style parametrization.

These potential failure modes -nodal, multi-valued, and fold-point- motivate the need for a more general and robust approach. The following sections will use a numerical example to illustrate these challenges and set the stage for introducing a more robust arc-length parametrization technique. 

\section{Application and Limitations of the Graph-Style Parametrization}\label{Section_example_GS_limitations}
\begin{figure}[h!]
	\centering
	\includegraphics[scale=0.7]{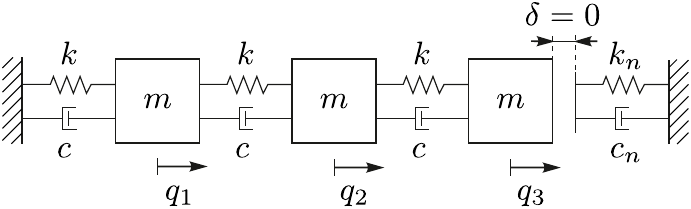}
	\caption{Schematic of the 3-DOF chain of oscillators with a discontinuous unilateral support. This system generalizes the benchmark from Attar et al. \cite{Attar2017} by including internal linear damping (via coefficient $d$) and unilateral contact damping (via coefficient $c_n$).}
	\label{fig_SimpleMechSys3DOF}
\end{figure}

To demonstrate the efficacy and expose the limitations of the graph-style parametrization, we introduce a benchmark problem that will serve as our primary example throughout the remainder of the paper. We consider a 3-DOF chain of oscillators with a unilateral contact, as depicted in Figure~\ref{fig_SimpleMechSys3DOF}. This system is an extension of the conservative model analyzed by Attar et al. in \cite{Attar2017}, which exhibits rich dynamics, including bifurcations and internally resonant NNMs alongside the fundamental NNM branches. Our version generalizes this benchmark by incorporating both internal linear damping and a discontinuous unilateral support, thus fitting precisely into the broader class of systems defined in Section 2. This simple yet insightful system serves as an ideal testbed for assessing the domain of viability for the graph-style parametrization and clearly demonstrates how its theoretical failure modes, discussed in Section 4, may manifest in numerical applications.

The system consists of three masses, each of mass $m$ connected by linear springs of stiffness $k$. The third mass is constrained by a unilateral support with spring stiffness $k_n$ and damping coefficient $c_n$. Following the general formulation of Eq.\eqref{eq:eom_2ndorder}, the system matrices are given by:
\begin{align}
	\vM &= m \vI_{3\times 3}, \quad
	\vC = \begin{pmatrix}
		2c & -c & 0 \\ -c & 2c & -c \\ 0 & -c & c
	\end{pmatrix}, \quad
	\vK = \begin{pmatrix}
		2k & -k & 0 \\ -k & 2k & -k \\ 0 & -k & k
	\end{pmatrix}.
\end{align}
The unilateral support acts only on the third mass, meaning the generalized force direction is ${\vw\T = \begin{pmatrix}
	0 & 0 & 1
\end{pmatrix}}$. The resulting PWL dynamics are fully described by Eq.\eqref{eq:eom_piecewise}, with the system being discontinuous whenever the unilateral damping $c_n>0$.
The original conservative benchmark from Attar et al. \cite{Attar2017} is recovered by setting all damping coefficients to zero ($c=0$ and $c_n=0$). As established in Remark 2.1, this underlying conservative system is of fundamental importance, as its periodic solutions define the conservative NNMs that form the backbone of the dynamics. Therefore, to construct and analyze these structures, we will first focus on this undamped, continuous case. In subsequent sections, non-zero damping will be considered, with the full discontinuous case ($c_n>0$) being analyzed in Section 7.

In smooth, nonlinear, undamped systems, NNM backbone curves typically present the nonlinear frequencies as a function of an energy-like quantity (e.g. oscillation amplitude). However, for the class of CPL systems of the form \eqref{eq:eom_cons_CPL}, the nonlinear frequencies are energy-independent and depend instead on the stiffness of the unilateral contact. For this reason, Attar et al. \cite{Attar2017} proposed a modified representation of NNM backbones as frequency-stiffness-plots (FSPs), which track the evolution of the nonlinear frequencies against the unilateral contact stiffness. The authors studied the three branches of fundamental NNMs alongside other bifurcating NNM branches and internal resonance tongues. In this work, we will take the first NNM as a benchmark to show where the graph-style parametrization of invariant cones is viable and where it breaks down. 
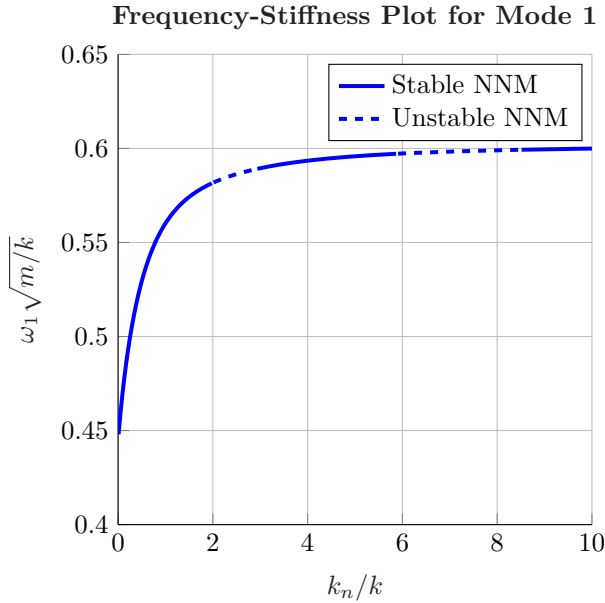
\begin{figure}[h]
	\centering	
	\tikzsetnextfilename{NNMbranch}
	\input{FSP_1.tikz}%
	\caption{Frequency-stiffness plot for the first NNM of the undamped $(c=0)$ 3-DOF chain of oscillators with no initial gap $(\delta=0)$. The vertical axis shows the dimensionless frequency, while the horizontal axis shows the dimensionless contact stiffness. Reference data from Attar et al. \cite{Attar2017}.}
	\label{fig_NNMbranch}
\end{figure}%
For the undamped case $(c=c_n=0)$, the first NNM branch is shown in Figure~\ref{fig_NNMbranch}. In the following subsections, we will select three distinct values of the contact stiffness from this branch to depict different scenarios for the cylinder-based graph-style parametrization and to assess its domain of validity. 
\subsection{Successful Application of the Graph-Style Parametrization}
\begin{figure}
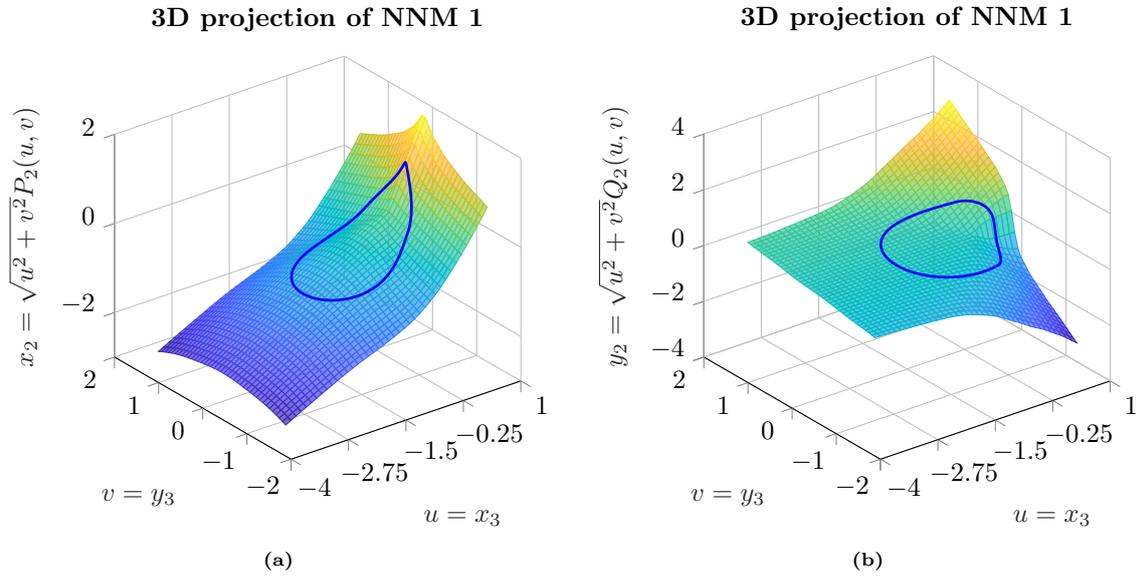

	\begin{subfigure}{0.49\linewidth}
		\centering
		\tikzsetnextfilename{consNNM1-sldisp}
		\input{NNM1_3D_sldisp.tikz}%
		\caption{}
		\label{fig:consNNM1_sldisp}
	\end{subfigure}
	\begin{subfigure}{0.49\linewidth}
		\centering
		\tikzsetnextfilename{consNNM1-slvelo}
		\input{NNM1_3D_slvelo.tikz}
		\caption{}
		\label{fig:consNNM1_slvelo}
	\end{subfigure}
	\caption{Invariant cone of the first NNM for the conservative system $(c=0)$, computed via the graph-style parametrization. The manifold geometry is shown for slave coordinates (a) $x_2$ and (b) $y_2$ as a function of the master coordinates $(x_3,y_3)$. A periodic orbit is shown in blue on the manifold surface.}
	\label{fig:consNNM1}
\end{figure}

\begin{figure}
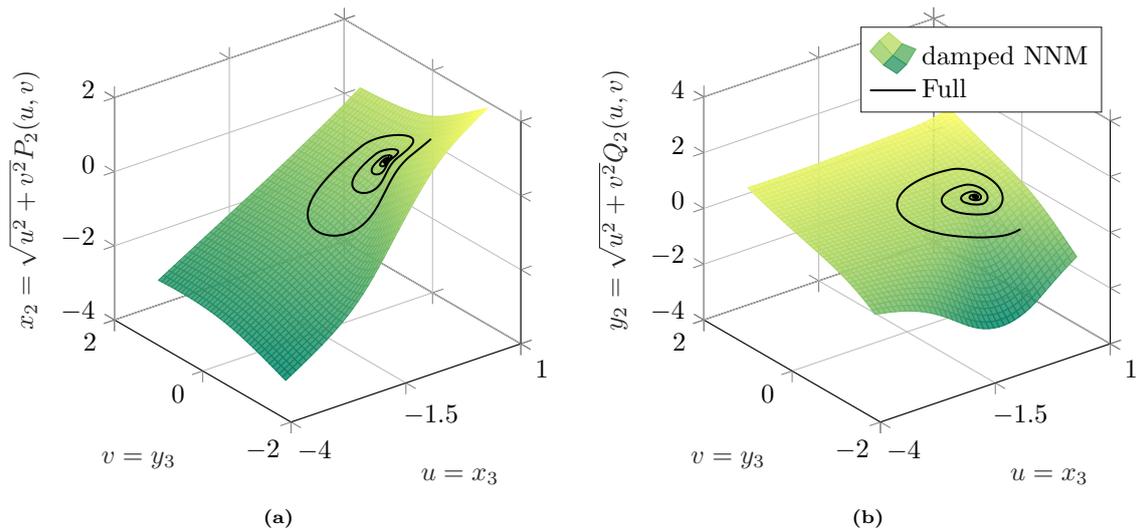

	\begin{subfigure}{0.49\linewidth}
		\centering
		\tikzsetnextfilename{NNM1-sldisp-damped}
		\input{NNM1_damped_traj_disp.tikz}%
		\caption{}
		\label{fig:NNM1_sldisp_damped_vs_conservative}
	\end{subfigure}
	\begin{subfigure}{0.49\linewidth}
		\centering
		\tikzsetnextfilename{NNM1-slvelo-damped}
		\input{NNM1_damped_traj_velo.tikz}
		\caption{}
		\label{fig:NNM1_slvelo_damped_vs_conservative}
	\end{subfigure}
	\caption{Invariant cone of the first NNM for the dissipative system ($c=1$). The geometry is shown for slave coordinates(a) $x_2$ and (b) $y_2$. A trajectory of the full system (black line), initiated on the manifold, remains on the surface as it decays towards the origin, numerically confirming the invariance property of the cone.}
	\label{fig:NNM1_damped}
\end{figure}

We first consider a case with a moderate contact stiffness ratio of $k_n/k = 1.5$. For this parameter set, the system's dynamics do not trigger the theoretical limitations outlined in Section 4.2.2. Consequently, the graph-style assumption holds, defining the slave coordinates as a single-valued function of the master coordinates. For this analysis, the displacement $x_3(t)$ and velocity $y_3(t) = \dot{x}_3(t)$ of the degree-of-freedom subject to the unilateral contact are chosen as the master pair. This choice simplifies the geometry ODEs governing the shape of the generating set. The parametrization of the invariant cone is then computed, yielding a functional relationship that expresses all slave variables in terms of the physical master coordinates $(x_3,y_3)$. The computation was performed using a harmonic order of $N_h=20$ for the slave variables. Figure \ref{fig:consNNM1} shows the resulting invariant cone for the conservative case $(c=c_n=0)$. The cone's surface is visualized, with a sample periodic orbit shown lying upon it, illustrating the geometry of the first NNM. The damped NNM for the continuous PWL obtained by setting $c=0.5$ and $c_n=0$ is depicted in Figure \ref{fig:NNM1_damped}. To numerically validate the invariance property of the computed cone, a trajectory of the full damped system was initiated on the cone and is shown in black. As expected, the trajectory remains on the parametrized surface for all forward time.

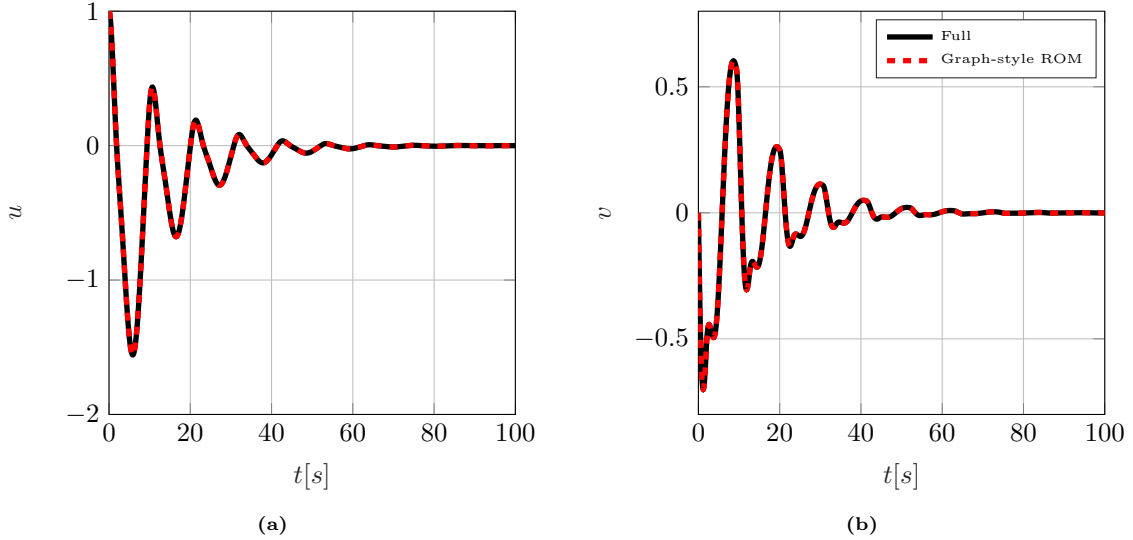
\begin{figure}
	\begin{subfigure}{0.49\linewidth}
		\centering
		\tikzsetnextfilename{NNM1-damped-timehistory-u}
		\input{NNM1_damped_timehist_u.tikz}%
		\caption{}
		\label{fig:NNM1_damped_timehistory_u}
	\end{subfigure}
	\begin{subfigure}{0.49\linewidth}
		\centering
		\tikzsetnextfilename{NNM1-damped-timehistory-v}
		\input{NNM1_damped_timehist_v.tikz}
		\caption{}
		\label{fig:NNM1_damped_timehistory_v}
	\end{subfigure}
	\caption{Time histories of the master coordinates on the damped NNM 1. Comparison between the full-order 6-DOF system (black) and the 2-DOF ROM (red dashed). (a) Master displacement $x_3(t)$. (b) Master velocity $y_3(t)$.}
	\label{fig:NNM1_damped_timehistories}
\end{figure}

Once the parametrization is computed, the reduced-order model (ROM) is obtained as an explicit 2D system of ordinary differential equations, following Eq. \eqref{eq:ROM_uv}. To demonstrate the accuracy of this ROM, Figure~\ref{fig:NNM1_damped_timehistories} depicts the time histories of the master variables obtained from both the full damped system (black) and the ROM (red). The excellent agreement between the two simulations confirms the high fidelity of the ROM in this case. 

\subsection{Failure Modes and the Need for a More Robust Approach}
\begin{figure}[ht!]
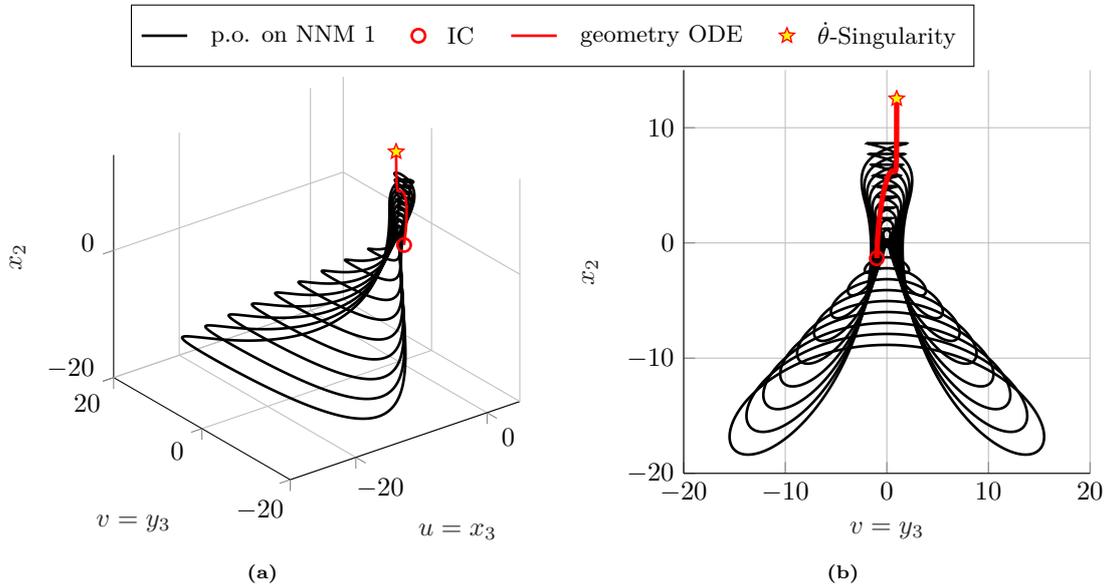
 % Added ! to encourage placement
	\centering
	
	% --- FIX 1: Isolate the legend in its own centered paragraph ---
	\parbox{\linewidth}{\centering
		\begin{tikzpicture}[font=\small]
			% --- FIX 3: Enable '&' as a column separator ---
			\matrix[
			matrix of nodes,
			ampersand replacement=\&, % <-- Added
			draw, thin,
			column sep=5pt,
			row sep=3pt,
			nodes={anchor=west}
			] (legend) {
				% Row 1: Line
				\draw[black, solid, line width=1.0pt] (0,0) -- (0.6,0); \& p.o. on NNM 1 \&\&
				% --- FIX 2: Use '\path plot' for markers ---
				% Row 2: Marker
				\path plot[color=red, line width=1.0pt, only marks, mark size=2.5pt, mark=o, mark options={solid, red, line width=1.0pt}] coordinates {(0,0)}; \& IC \&\&
				% Row 3: Line
				\draw[red, solid, line width=1.0pt] (0,0) -- (0.6,0); \& geometry ODE \&\&
				% Row 4: Marker
				\path plot[mark=mystar, only marks, mark size=3pt] coordinates {(0,0)}; \& $\dot{\theta}$-Singularity \\
			};
		\end{tikzpicture}
	}
	% --- Your Subfigures ---
	\begin{subfigure}{0.49\linewidth}
		\centering
		\tikzsetnextfilename{HighNonlin-thetadot-singularity}
		\input{folding_thetadot_singularity_3D.tikz}%
		\caption{}
		\label{fig:HighNonlin_thetadot_singularity}
	\end{subfigure}% <--- IMPORTANT: No blank line after this
	\begin{subfigure}{0.49\linewidth}
		\centering
		\tikzsetnextfilename{HighNonlin-thetadot-VX2Portrait}
		% Make sure you have REMOVED the legend from inside this .tikz file
		\input{folding_thetadot_singularity_2D.tikz}
		\caption{}
		\label{fig:HighNonlin_thetadot_singularity_VX2Portrait}
	\end{subfigure}
	
	\caption{Failure of the graph-style parametrization for a large contact stiffness ratio $(k_n/k)=10$. The invariant cone features a fold, leading to a singularity $(\dot{\theta})$ in the geometry ODE. (a) Three-dimensional view showing periodic orbits foliating the cone. (b) Projection onto the $(v,x_2)$ plane.} % Shortened caption for brevity
	\label{fig:HighNonlin_Folding}
\end{figure}

We now consider the case of a high contact stiffness ratio of $k_n/k=10$. For this highly nonlinear scenario, a periodic orbit on the first conservative NNM ($c=c_n=0$) is computed using the shooting method. Due to the positive homogeneity property of the system, the invariant cone is foliated by scaled versions of this periodic orbit. 
The resulting cone, shown in Figure~\ref{fig:HighNonlin_Folding}, exhibits a complex folded geometry when projected into the three-dimensional space spanned by the master coordinates $(x_3,y_3)$ and a slave coordinate (e.g. $x_2$). This folding is a direct consequence of the looping trajectories of the periodic solutions in the configuration space. As a result, the generating set - defined as the intersection of the cone with the cylinder $\mathcal{H}_{\textrm{cyl}}$ - has to fold on itself in order to reproduce the folding structure of the cone. This failure is demonstrated numerically. We integrate the geometry ODEs \eqref{eq:geometryODE} from an initial point on the generating set near the fold. The solution trajectory, shown in red in Fig.\ref{fig:HighNonlin_Folding}, successfully traces the generating set until it reaches the fold point, marked by a yellow star. At this point, the angular velocity in the cylindrical parametrization vanishes $(\dot{\theta} \rightarrow 0)$, creating a singularity in the geometry ODEs that halts the integration. This numerical breakdown is a direct manifestation of the fold-point failure, as theoretically described in Subsection 4.2.2. The inability of the graph-style method to handle such folded geometries necessitates a more robust approach.

\section{Arc-length parametrization for folding cones}\label{Section_ArcLength}
The failure modes of the cylinder-based graph-style method, established theoretically in Section 4 and demonstrated numerically in Section 5, can be traced to two fundamental choices in its formulation: the intersecting hypersurface and the curve parameter $\theta$. The arc-length parametrization, derived herein, is specifically designed to remedy both deficiencies. 
First, the nodal failure, i.e., the case where the ray-to-generator mapping becomes a partial function, arises because a ray on the cone can be parallel to the cylinder's axis, thus failing to produce an intersection. This is a geometric flaw in using the cylinder $\mathcal{H}_{\textrm{cyl}}$ as the intersecting hypersurface. We resolve this problem by replacing the cylinder with the unit hypersphere $\mathcal{S}^{(n-1)}$. A sphere centered at the origin is guaranteed to intersect every ray from the origin exactly once, ensuring that the generating set $\mathcal{G}$ is always well-defined. 
Second, the remaining two failure modes, namely the multi-valued nature of the generating curve and the fold-point singularity of the geometry ODE \eqref{eq:geometryODE}, are consequences of using the angle $\theta$ of the master pair $(u,v)$ as the curve parameter. This parametrization fails whenever the generating curve folds or loops with respect to $\theta$. We can overcome this by using the curve's intrinsic arc-length, $s$. The arc-length is, by definition, a single-valued and strictly increasing quantity that traces the curve's path independently of its geometry, which makes the parametrization robust even in the presence of folds. 
This new arc-length method is therefore built upon these two key ideas: intersection with a unit sphere to define the generating set, and an arc-length parametrization to trace it. The following derivation details this approach. 
\subsection{Derivation}
We begin by defining the invariant cone as the set
\begin{equation}\label{eq:arclength_Ansatz_conedefinition}
	\mathcal{C} = \{  \vC(r,s)= r \vG(s) \in \mR^{n} \, | \, \vG \in \mathcal{G} \subset \mathcal{C}, \, r\in [0, \infty) \},
\end{equation}
where $\vG$ is the generating curve of the cone and $s$ is its arc-length parameter. The scalar $r$ represents the magnitude of a point on the cone, i.e., its Euclidean distance from the origin. This 
definition implies two fundamental geometric constraints: 
\begin{equation}
	r = ||\vC(r,s)||, \quad \mbox{and } \quad ||\vG(s) || = 1 .
\end{equation}
This first constraint signifies that the generating curve $\mathcal{G}$ lies on the unit hypersphere $\mathcal{S}^{(n-1)}$ in $\mR^n$. This is a crucial departure from the graph-style method, which used a cylinder $\mathcal{H}_{\textrm{cyl}}$ defined over a $2$D master subspace. 
The second constraint arises from the definition of $s$ as the arc-length parameter. A curve parametrized by its own arc-length moves with unit speed along its path. Therefore the tangent vector to the generating curve must have a unit norm: 
\begin{equation}
	||\vG'(s)|| = 1,
\end{equation}
where $(\cdot)'$ denotes differentiation with respect to $s$. These two geometric conditions, $||\vG(s)|| \stackrel{!}{=} 1$ and $||\vG'(s)|| \stackrel{!}{=} 1 $ form the foundation of this parametrization. 
Our goal is to derive a system of ODEs that defines the generating curve $\vG(s)$. This is achieved by imposing the invariance condition. Any trajectory $\vC(t)$ that starts on the cone must remain on it. Its time evolution is governed by the dynamics $\dot{\vC} = \vf(\vC)$, where $\vf$ is the positively homogeneous, piecewise-linear vector field from Eq. \eqref{eq:eom_firstorder}. Applying the chain rule to $\vC(t) = r(t)\vG(s(t))$ gives:
\begin{equation}
	\dot{\vC} = \dot{r} \vG(s) + r \vG'(s) \dot{s}.
\end{equation}
Substituting this into the system dynamics and using the homogeneity of $\vf$ (i.e., $\vf(r\vG) = r\vf(\vG)$) yields:
\begin{equation}
	\dot{r} \vG(s) + r \vG'(s) \dot{s} = r \vf(\vG(s)).
\end{equation}
Dividing by the scalar magnitude $r>0$ (for points not at the origin), we obtain the central invariance equation:
\begin{equation}\label{eq:geomODE_arclength}
	\vG'(s) = \frac{\vf(\vG) - (\dot{r}/r) \vG(s)}{\dot{s}}.
\end{equation}
This equation separates the problem into two parts: finding the geometry of the manifold, $\vG(s)$, and finding the reduced dynamics on the manifold, described by the evolution of the new coordinates $(r,s)$.

First, we derive the reduced dynamics. The evolution of the radial coordinate, $\dot{r}$, is found by differentiating $r^2 = \vC\T\vC$:
\begin{equation}
	2r\dot{r} = \frac{\mathrm{d}}{\mathrm{d}t}(\vC\T\vC) = 2\vC\T\dot{\vC} = 2(r\vG)\T \vf(r\vG) = 2r^2 \vG\T\vf(\vG).
\end{equation}
This gives the rate of radial expansion as the projection of the normalized vector field onto the radial direction:
\begin{equation}\label{eq:radial_dynamics}
	\frac{\dot{r}}{r} = \vG(s)\T\vf(\vG(s)).
\end{equation}

Next, we substitute this radial dynamics back into the invariance equation \eqref{eq:geomODE_arclength}:
\begin{equation}\label{eq:definitionT}
	\vG'(s)\dot{s} = \vf(\vG(s)) - \left[ \vG(s)\T \vf(\vG(s)) \right] \vG(s) =: \vT(\vG(s)).
\end{equation}
The vector field $\vT(\vG(s))$ has a clear geometric meaning: it is the component of the normalized dynamics $\vf(\vG)$ that is tangent to the unit sphere $\mathcal{S}^{(n-1)}$ at the point $\vG(s)$. Taking the norm of this equation and applying the unit-speed constraint $||\vG'(s)||=1$, we find the rate of evolution along the generating curve:
\begin{equation}\label{eq:arclength_dynamics}
	\dot{s} = ||\vT(\vG(s))||.
\end{equation}

Finally, substituting this expression for $\dot{s}$ back into the invariance equation provides the explicit, singularity-free geometry ODE:
\begin{equation}\label{eq:arclength_geomODE}
	\vG' = \frac{\vT(\vG)}{||\vT(\vG)||}.
\end{equation}
This is the main result of our derivation. This first-order ODE system governs the geometry of the generating curve $\vG(s)$ defining the invariant cone. Unlike the graph-style formulation, its right-hand side is always well-defined and has a unit norm, except at fixed points of the normalized dynamics where $\vT(s)=0$. It is crucial to analyze this singular condition. By its definition in Eq. \eqref{eq:definitionT}, and recalling that $||\vG(s)||=1$, the term $\left[ \vG(s)\T \vf(\vG(s)) \right] \vG(s)$ is the orthogonal projection of the vector field $\vf(\vG(s))$ onto the radial direction $\vG(s)$. The vector $\vT(\vG(s))$ is thus the component of $\vf(\vG(s))$ that is tangent to the unit hypersphere. The condition $\vT(\vG(s)) = \vzero$ therefore implies that the vector field is purely radial, i.e., $\vf(\vG(s)) = \gamma \vG(s), \quad \text{where } \gamma = \vG(s)\T \vf(\vG(s)) \in \mathbb{R}$. For the PWL systems under consideration, this means $\vA^\pm \vG(s) = \gamma \vG(s)$. This is the definition of $\vG(s)$ being a real eigenvector of a subsystem matrix ($\vA^+$ or $\vA^-$) and $\gamma$ being its corresponding real eigenvalue.The invariant cones we seek to parametrize, however, represent the oscillatory dynamics (the NNMs) of the system, which are associated with the complex-conjugate eigenvalues of the subsystem matrices $\vA^\pm$. By assuming these subsystems are purely oscillatory (i.e., the matrices possess no real eigenvalues), the condition $\vT(s) = \vzero$ can never occur for any $\vG(s) \neq \vzero$. Therefore, for the oscillatory invariant cones of interest, the denominator $||\vT(\vG(s))||$ is strictly positive, and the geometry ODE is singularity-free.
This formulation is therefore robust against folds and provides a more robust approach for parametrizing complex invariant cone geometries.

The invariant cone $\mathcal{C}$ is, by definition, a union of rays that is invariant under the flow. This means that any trajectory $\vC(t)$ starting on the cone remains on it for all time. For the oscillatory PWL system under consideration, trajectories on the cone are typically periodic orbits or trajectories that decay or grow while spiraling on the cone's surface. In all these cases, a trajectory starting at $\vC(0)$ will return to the same ray after some time $T>0$. This implies that the direction vector is periodic, even if the magnitude is not. This can be expressed as $\vC(T) = \mu \vC(0)$, for some scaling factor $\mu \in \mR^+$. Applying our arc-length parametrization, $\vC(t) = r(t) \vG(s(t))$, to this return map condition gives: 
\begin{equation}\label{eq:vector_return_map}
	r(T) \vG(s(T)) = \mu \, r(0) \vG(0).
\end{equation}
Since $\vG$ is a unit vector, taking the Euclidean norm of both side yields $r(T) = \mu r(0)$. Substituting back into \eqref{eq:vector_return_map} and canceling the non-zero scalar magnitudes yields the fundamental condition on the geometry:
\begin{equation}
	\vG(s(T)) = \vG(s(0)).
\end{equation}
This proves that the generating set $\vG(s)$ must be a closed loop on the unit hypersphere $\mathcal{S}^{(n-1)}$. We therefore seek a periodic solution to the geometry ODE \eqref{eq:arclength_geomODE} with period $s(T)=L$, which corresponds to the total length of the generating curve. The corresponding arc-length frequency is $\Omega_s = \frac{2\pi}{L}$. While this periodic solution can be found numerically using a boundary value problem (BVPs) solver, such as shooting, this work exploits the periodicity of $\vG(s)$ and employs a spectral approach  to construct a computationally efficient ROM. The generating curve $\vG(s)$ is represented by a truncated Fourier series: 
\begin{equation}\label{eq:FourierGs}
	\vG(s) = \vX_{0} + \sum_{k=1}^{N_h} \left[ \vX_{c,k} \cos(k \Omega_s s) + \vX_{s,k} \sin(k \Omega_s s) \right].
\end{equation}
This procedure begins with a one-time, offline computation where the Fourier coefficient vectors $\vX_{0}, \vX_{c,k}, \vX_{s,k}$  are determined up to a desired harmonic order, $N_h$. These coefficients provide a compact analytical representation of the manifold's geometry, replacing the need for pointwise interpolation during simulation. In the subsequent online stage, the two-dimensional reduced dynamics $(\dot{r}, \dot{s})$ is numerically integrated, wherein the Fourier coefficients are used to evaluate the generating set $\vG(s)$ and also the dynamics $\vf(\vG(s))$ needed for the right-hand side in equations \eqref{eq:radial_dynamics} and \eqref{eq:arclength_dynamics}. The full-order system state can then be efficiently reconstructed at any desired time instant by evaluating the truncated Fourier series \eqref{eq:FourierGs} at $s(t)$ and scaling the result by $r(t)$. 

\subsection{Utility of the Geometric Approach and the Reduced-Order Model}
Having derived a singularity-free ODE for the cone's generating set, Eq. \eqref{eq:arclength_geomODE}, we now address its solution and the subsequent advantages of this geometric framework. Since the generating curve $\vG(s)$ must be a closed loop on the unit hypersphere for any invariant cone, we seek its periodic solutions. These can be found by formulating the problem as a two-point BVP, enforcing the periodicity condition $\vG(L)=\vG(0)$, and solving it with a numerical method such as the shooting method or the HBM.
At first glance, solving for a single periodic orbit on the cone might seem comparable in effort to solving for one periodic orbit of the full $N$-dimensional system. Due to the positive homogeneity of the autonomous CPL system, all other trajectories are simply scaled versions of this one solution, from which one may have the impression that a subsequent 2D ROM is redundant. This, however, overlooks the advantages of this approach, particularly for dissipative and forced systems.
The arc-length parametrization mainly decouples the cone's geometry from whether the system's dynamics is damped or conservative. The BVP solver applied to the arc-length geometry ODE solves for the shape of the invariant cone - the generating set $\vG(s)$. The usefulness of this fact becomes clear in the damped case. As opposed to a conservative system, a damped autonomous system does not possesses periodic orbits. Standard shooting methods for periodic orbits on the full system are therefore not applicable. However, the generating curve $\vG(s)$ remains periodic even for damped systems. Thus, using a periodic solution solver such as shooting or HBM to compute the periodic generating set is still viable and allows for the parametrization of the invariant cone containing the decaying trajectories. Furthermore, from this parametrization, we can directly compute the scalar $\mu= r(T)/r(0)$ characterizing the stability of trajectories on the cone, via a simple line integral over the known curve $\vG(s)$: 
\begin{equation}\label{eq:mu_integral}
	\ln(\mu) = \int_0^L \frac{\vG(s)\T \vf(\vG(s))}{||\vf(\vG(s)) - \left( \vG(s)\T \vf(\vG(s)) \right)\vG(s)||} \diff{s}
\end{equation}
Most importantly, this framework lays the essential groundwork for analyzing forced systems. The positive homogeneity property, which allows simple scaling of trajectories, is broken by the addition of an external forcing term. In the smooth framework, model order reduction methods based on the mathematically rigorous direct parametrization method of invariant manifolds often regard the time-dependent forced manifolds as perturbations of the autonomous invariant manifolds computed for the underlying unforced dynamics. For PWL systems, similar ideas can be followed to determine forced perturbations of invariant cones, using the autonomous parametrization investigated here as a starting point. While a detailed analysis of the forced problem is beyond the scope of the present work, this methodology provides the foundation for such investigations and is currently the subject of ongoing research by the authors.  

\section{Numerical Validation}\label{Section_NumValidation}
The advantages of the arc-length parametrization over the traditional graph-style approach are now demonstrated through a series of targeted numerical tests. Each test is designed to address a known limitation of the graph-style method. We first investigate the parametrization of an internally resonant NNM as a multiple-crossing invariant cone to address the multi-valued failure. Second, we analyze the method's performance to parametrize folding invariant cones for high-stiffness conservative and dissipative continuous systems. Finally, the robustness of the framework is tested on a discontinuous system, which combines a geometrically folded manifold with a discontinuous PWL vector field. In all scenarios, the arc-length method provides a consistent parametrization of the invariant manifold, leading to accurate two-dimensional ROMs.

\subsection{Application to an Internally Resonant NNM}
\begin{figure}[htbp]
	\centering
	\begin{subfigure}{0.43\linewidth}
		\includegraphics[width=\linewidth]{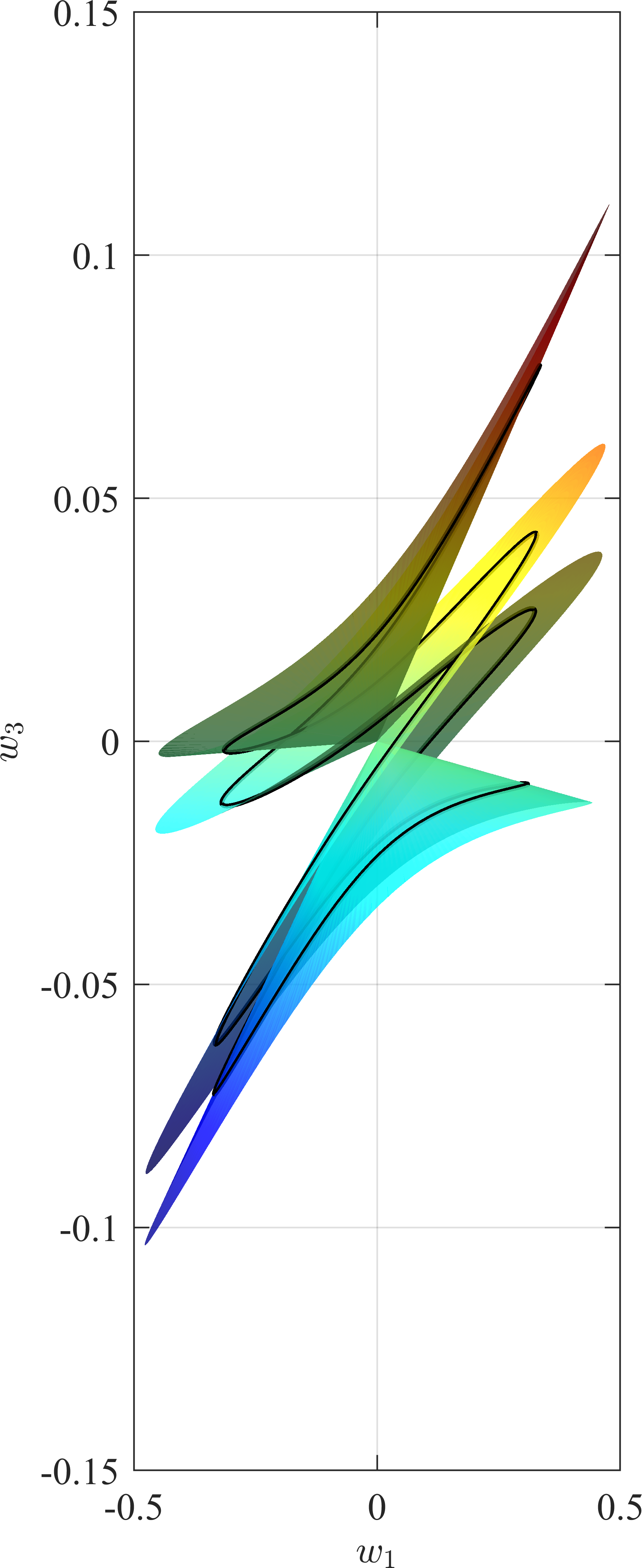}%
		\caption{Projection onto the $(w_1, w_3)$ plane.}
		\label{fig:stretchedCone_w1w3}
	\end{subfigure}
	\hfill % Adds a horizontal space that pushes the subfigures to the edges
	\begin{subfigure}{0.55\linewidth}
		\includegraphics[width=\linewidth]{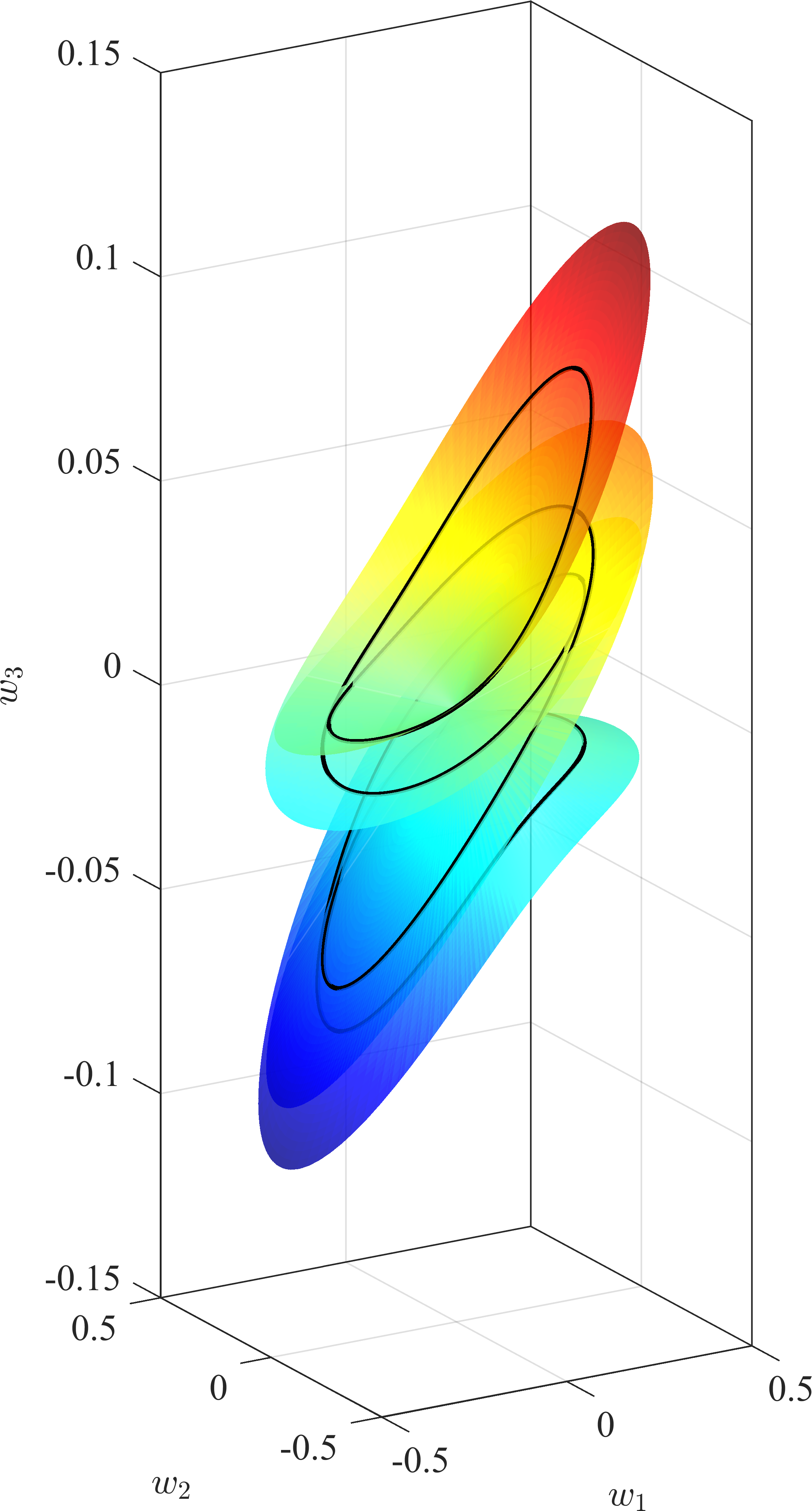}
		\caption{Three-dimensional view.}
		\label{fig:stretchedCone}
	\end{subfigure}
	\caption{Figure shows the invariant cone projected onto a rotated coordinate system $(w_1,w_2,w_3)$. The transformation is defined by $w_1 = (x_2-x_1)/\sqrt{2}, w_2 = x_5$, and $w_3 = (x_2+x_1)/\sqrt{2}$, which aligns the $w_3$-axis perpendicular to the cone's central plane. To reveal the geometric details, the data aspect ratio has been set to $(1 1 0.1)$, effectively stretching the $w_3$-axis by a factor of $10$ for visualization purposes.}
	\label{fig:resonantCone}
\end{figure}

\begin{figure}[htbp]
	\centering	
	\newlength\figurewidth
	\newlength\figureheight
	\setlength\figurewidth{13cm}
	\setlength\figureheight{10.4cm}
	\tikzsetnextfilename{InternalRes-timehistories}
	\input{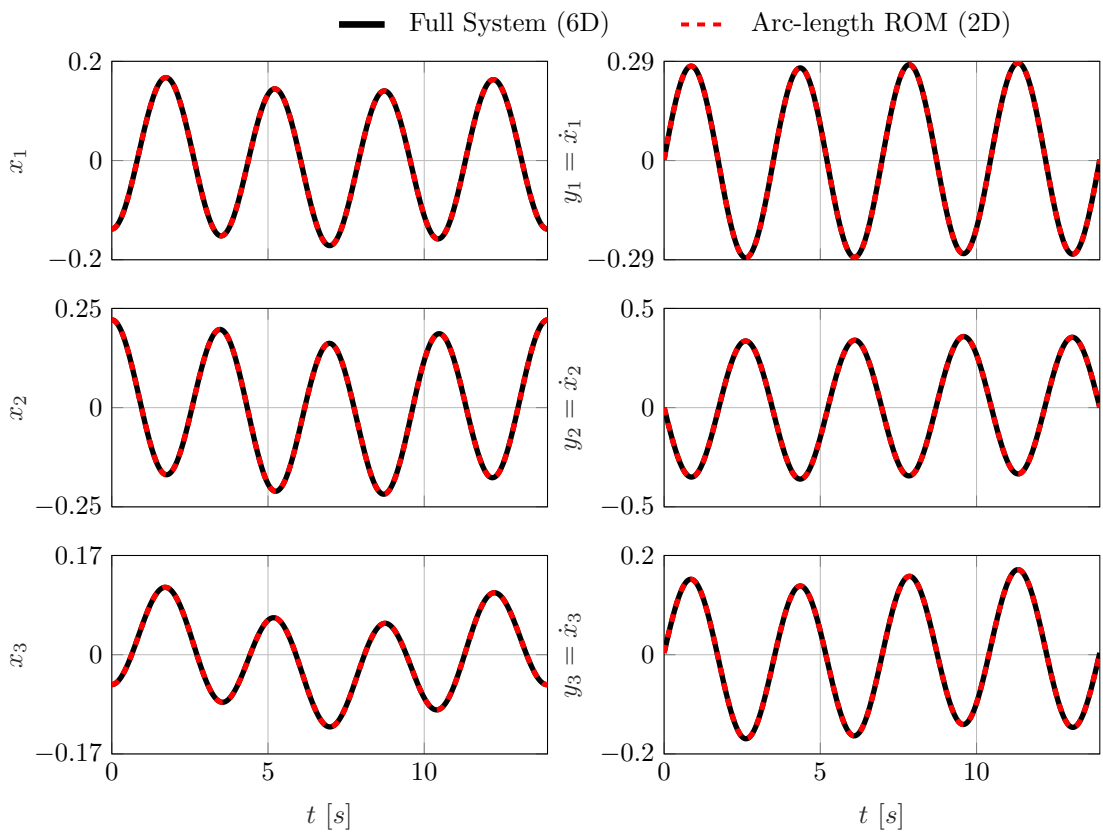}	
	\caption{Time histories of the full system and ROM on the folding cone at the internal resonance.}
	\label{fig_InternalRes_timehistories}
\end{figure}

We first investigate a case known to exhibit complex modal interactions. Attar et al. identified a 4:1 internal resonance between the first and third NNMs of the undamped 3-DOF system for a contact stiffness ratio near $k_n/k=0.018$ and with $c=c_n=0$. At this resonance, the geometry of the first NNM undergoes a bifurcation, leading to trajectories that exhibit multiple crossings per period. For a graph-style approach, this manifests as a multi-valued generating curve, rendering the parametrization invalid. Here, we demonstrate that the arc-length parametrization successfully captures this complex, multiple-crossing invariant cone. Using HBM as a BVP solver for the arc-length geometry ODE \eqref{eq:geomODE_arclength}, we compute the generating set of the first NNM at the internal resonance point. The resulting invariant cone is visualized in Figure~\ref{fig:resonantCone}. To better reveal its intricate structure, the cone is projected onto a rotated coordinate system $(w_1,w_2,w_3)$ where the $w_3$-axis is orthogonal to the cone's central plane, and the aspect ratio is adjusted for clarity. Figure \ref{fig:stretchedCone_w1w3} shows a 2D projection that clearly reveals the self-intersecting nature of the cone, with a periodic orbit (black line) exhibiting multiple loops. The 3D view in Figure~\ref{fig:stretchedCone} confirms the geometry is a non-trivial, multi-sheeted surface. The ability to parametrize this cone demonstrates a key advantage of the arc-length approach: by tracing the intrinsic path of the generating curve, it remains robust even when the cone's projection onto a lower-dimensional subspace becomes multi-valued.
Once the cone's geometry is parametrized, the 2D ROM is constructed using equations \eqref{eq:radial_dynamics} and \eqref{eq:arclength_dynamics}. Figure \ref{fig_InternalRes_timehistories} shows the time histories of all six state variables, comparing the trajectory from the full 6D system with that predicted by the 2D arc-length ROM. The agreement is excellent, confirming that the 2D model accurately captures the dynamics on this complex, internally resonant manifold. This result is significant, as it shows that even when modal interactions are strong, the essential dynamics can still be captured by a two-dimensional model, provided the underlying geometry is correctly parametrized.
\subsection{Parametrization of a Folding Invariant Cone for Conservative and Dissipative Systems}
\begin{figure}[htbp]
	\begin{subfigure}{0.49\linewidth}
		\centering	
		\includegraphics{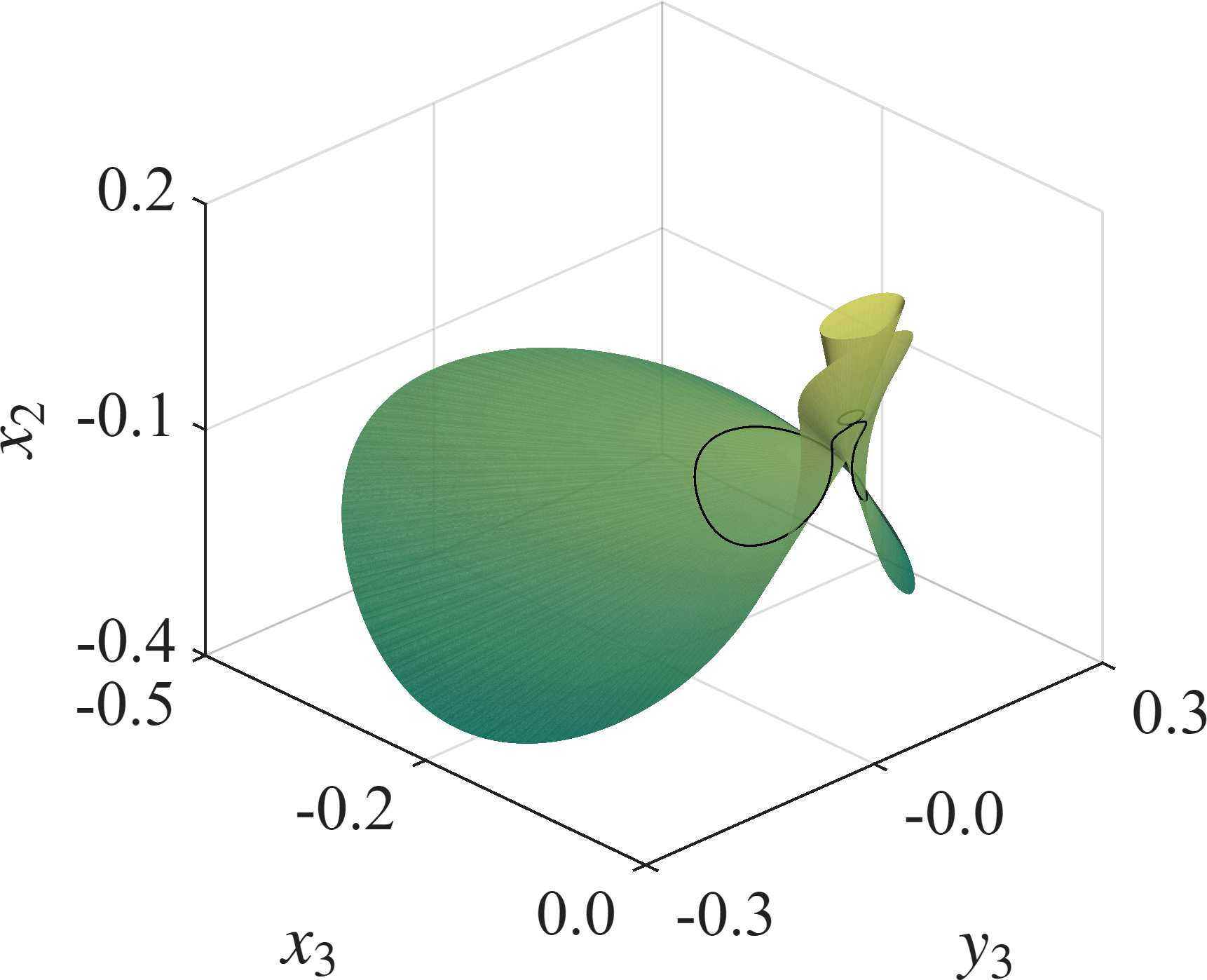}%
		\caption{}
		\label{fig:HighNonlin_arclength}
	\end{subfigure}%
	\begin{subfigure}{0.49\linewidth}
		\centering	
		\tikzsetnextfilename{HighNonlin-timehist-x3}
		\input{Highnonlin_timehist_x3_new.tikz}%
		\caption{}
		\label{fig:HighNonlin_timehist_x3}
	\end{subfigure}%
\caption{Validation of the arc-length ROM for the high-stiffness conservative case ($k_n/k=10, \, c=0$). (a) The computed folding invariant cone for the first NNM, with the generating curve $\vG(s)$ shown in black. The geometry is parametrized using $N_h=12$ harmonics. (b) Time history of the switching variable $x_3(t)$, comparing the full 6D system with the 2D arc-length ROM, showing excellent agreement.}
\label{fig:Highnonlin}
\end{figure}

\begin{figure}[htbp]
	\begin{subfigure}{0.49\linewidth}
		\centering	
		\includegraphics{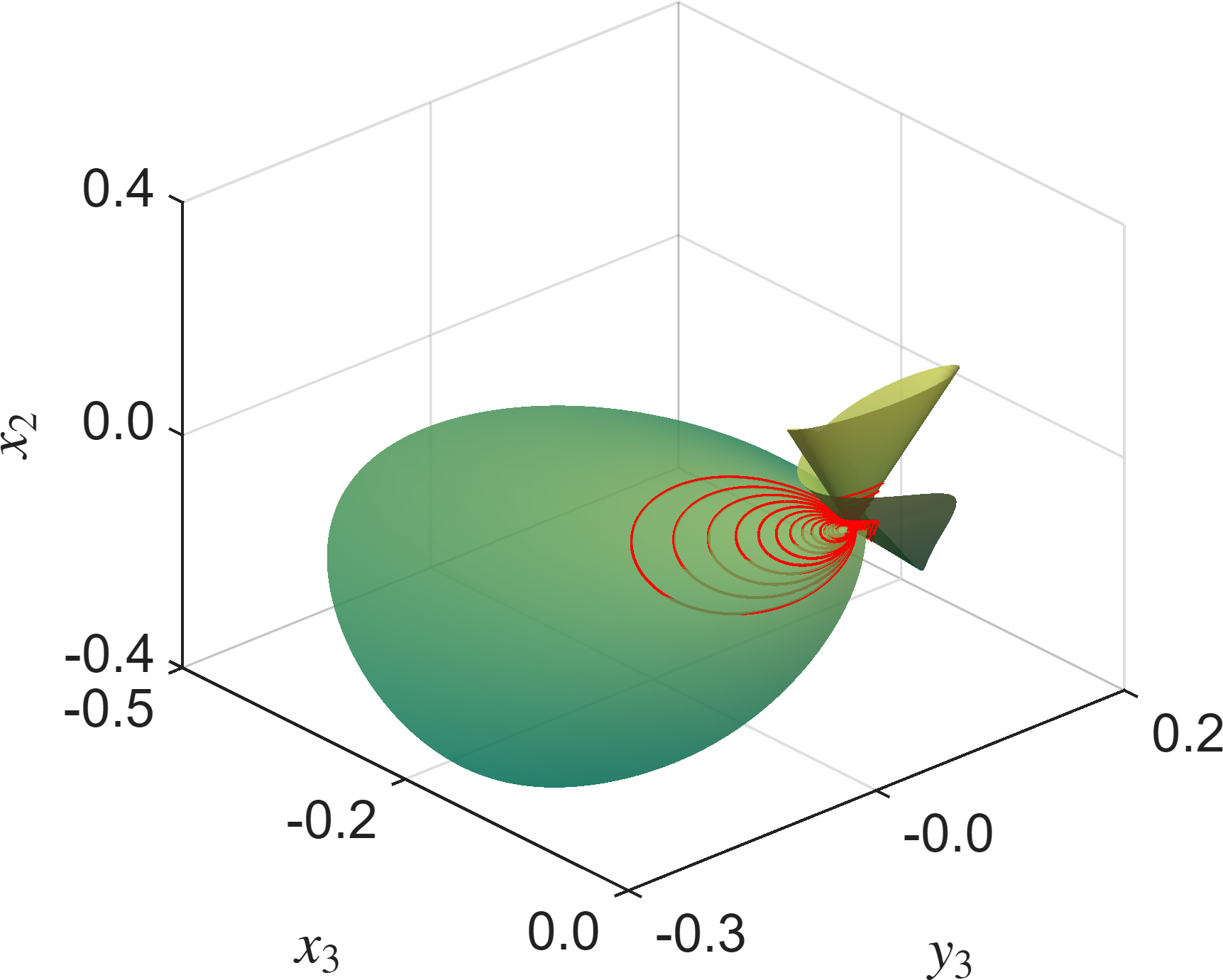}%
		\caption{}
		\label{fig:HighNonlin_arclength_Damped}
	\end{subfigure}%
	\begin{subfigure}{0.49\linewidth}
		\centering	
		\tikzsetnextfilename{HighNonlin-timehist-x3-damped}
		\input{Highnonlin_0x0295damping_x3timehist.tikz}%
		\caption{}
		\label{fig:HighNonlin_timehist_x3_damped}
	\end{subfigure}%
	\caption{Validation of the arc-length ROM for the high-stiffness dissipative case ($k_n/k=10, \, c=0.0295$). (a) The damped invariant cone, showing a trajectory (red spiral) of the full system initiated on the manifold decaying towards the origin. (b) Time history of the switching variable $x_3(t)$, comparing the full system to the 2D ROM.}
	\label{fig:Highnonlin_damped}
\end{figure}

We now return to the high-stiffness case $k_n=10 k$ that led to the breakdown of the graph-style method due to a fold in the invariant cone. By applying the arc-length parametrization, we can overcome the $\dot{\theta}=0$ singularity and accurately characterize the manifold's geometry for both conservative and dissipative scenarios. First, for the undamped system $c=c_n=0$, the geometry ODE, Eq. \eqref{eq:arclength_geomODE}, is solved using HBM with a harmonic order of $N_h=12$ to find the periodic generating curve $\vG(s)$ on the unit hypersphere. The resulting invariant cone is successfully parametrized over its entire folded geometry without encountering any singularities, as shown in Figure~\ref{fig:HighNonlin_arclength}. This three-dimensional projection reveals a highly distorted "fin"-like structure, which is a direct consequence of the stiff unilateral contact significantly restricting motion in the positive $x_3$ direction. With the cone's geometry now fully characterized, we validate the corresponding 2D ROM. Figure~\ref{fig:HighNonlin_timehist_x3} compares the time history of the switching variable, $x_3(t)$, as computed by the full 6D system and the 2D arc-length ROM. The two trajectories are visually indistinguishable, demonstrating the accuracy of the ROM even in this highly nonlinear regime. Next, we demonstrate the method's capability to handle dissipative systems, a case where traditional periodic-orbit-based methods fail. We introduce linear damping with a coefficient of $c=0.0295$ while keeping the high contact stiffness and a unilateral support with only a unilateral spring and no damping $c_n=0$. As argued in Section 6.2, the generating curve remains periodic and is computed again from the arc-length geometry ODE \eqref{eq:geomODE_arclength}. The resulting damped invariant cone is shown in Figure \ref{fig:HighNonlin_arclength_Damped}. While retaining the overall folded topology, the geometry is noticeably altered by the presence of damping. To confirm its invariance, a trajectory of the damped continuous system is initiated on the manifold (red spiral). The trajectory remains on the cone's surface as it decays towards the stable equilibrium at the origin. Furthermore, the computed parametrization allows for a direct, semi-analytical calculation of the decaying rate of trajectories on the cone. By evaluating the integral in Eq. \eqref{eq:mu_integral} over the computed generating curve $\vG(s)$, we obtain the cone's characteristic multiplier $\mu=0.8066$, which governs the rate of decay on the cone. The accuracy of the resulting 2D ROM for the damped continuous system is confirmed in Figure \ref{fig:HighNonlin_timehist_x3_damped}, where the time history of the switching variable $x_3(t)$ shows excellent agreement between the full system and the ROM. This result confirms that the arc-length parametrization provides a robust and reliable framework for model order reduction, successfully capturing the dynamics on complex invariant cones. 

\subsection{Validation on a Nonsmooth System with a Folding Invariant Manifold}

Having computed the explicit parametrization of the invariant cones for the continuous version of the 3-DOF system, both with and without linear damping, the final validation of our proposed arc-length parametrization is presented here, addressing a system that combines a geometrically folding invariant cone with a discontinuous vector field. 

For this, the 3-DOF system is configured as follows. A high unilateral contact stiffness, $k_n/k=2.5377$, is chosen to induce the geometric fold in the manifold. To introduce the discontinuity, the unilateral damper is activated with $c_n = 1$, while all internal linear damping is canceled ($c=0$). The damping discontinuity at the contact interface is therefore the sole source of dissipation, and the system's dynamics are governed by the Filippov framework.

\begin{figure}[htbp]
	\begin{subfigure}{0.49\linewidth}
		\centering	
		\includegraphics{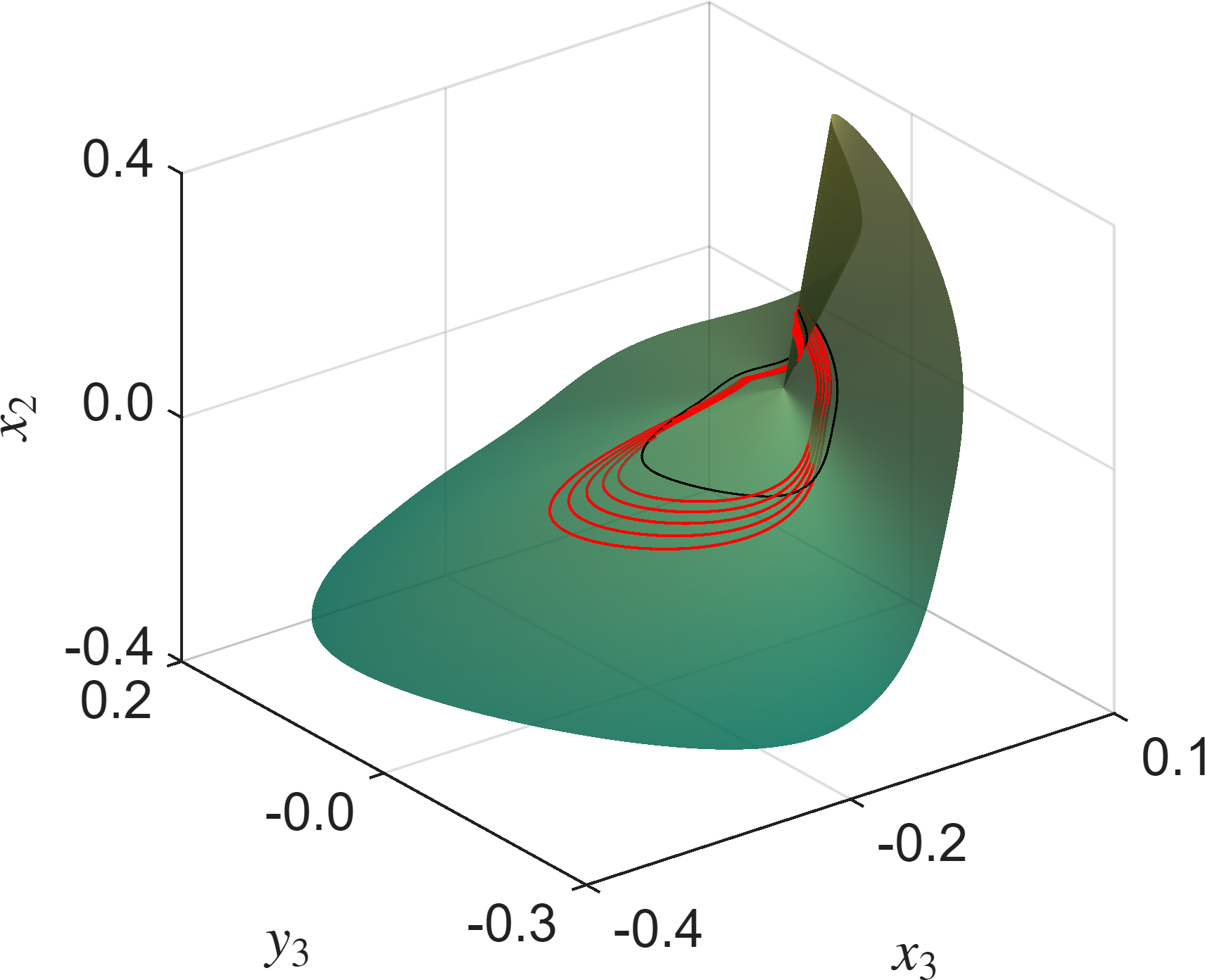}%
		\caption{}
		\label{fig:Discont_arclength_Damped}
	\end{subfigure}%
	\begin{subfigure}{0.49\linewidth}
		\centering	
		\tikzsetnextfilename{Discont-timehist-x6-damped}
		\input{discont_x6_timehist.tikz}%
		\caption{}
		\label{fig:Discont_timehist_y3_damped}
	\end{subfigure}%
	\caption{Validation of the arc-length ROM for the discontinuous case with ($k_n/k=2.5377, \, c=0, \, c_n = 1$). (a) The damped invariant cone, showing a trajectory (red spiral) of the full system initiated on the manifold decaying towards the origin. The red trajectory has an apparent kink at $x_3=0$ characterizing the impact with the discontinuous support  (b) Time history of the velocity of the impacting degree of freedom $y_3(t)$, comparing the full system to the 2D ROM.}
	\label{fig:Discont_damped}
\end{figure}

First, the arc-length technique described in Section \ref{Section_ArcLength} is applied to parametrize the invariant cone of the discontinuous PWL system using a moderate harmonic order $N_h=12$. The resulting manifold, depicted in Figure~\ref{fig:Discont_arclength_Damped}, retains the folded topology observed in the continuous case. However, its geometry is altered by the discontinuity; the decaying trajectory plotted on its surface (red spiral) exhibits a kink where it crosses the switching boundary $\Sigma_{\alpha}$ at $x_3=0$ with positive $y_3>0$, which signals a discontinuity in the contact force between the third mass and the support. 

The utility of this cone parametrization is demonstrated by its corresponding 2D ROM. Figure~\ref{fig:Discont_timehist_y3_damped} compares the time history of the impacting velocity, $\dot{x}_3(t)=y_3(t)$, between the full discontinuous PWL system and the arc-length 2D ROM. A close correspondence between the two is observed, not only with respect to the overall decay but also at the kinks. This result confirms the applicability of the arc-length parametrization of invariant cones even with a folding geometry and in the presence of a discontinuity. The ability to derive a two-dimensional ROM even for a strongly nonlinear setting and a folding geometry of the cone shows the robustness of the proposed arc-length parametrization technique even for discontinuous PWL systems.

\section{Conclusions}\label{Section_Conclusion}
In this work, we presented a systematic framework for computing invariant manifolds -- in this context, invariant cones -- for a class of homogeneous PWL systems. The primary outcome is the derivation of two-dimensional reduced-order models (ROMs) that accurately capture the dynamics of the full-order system.

A central aspect of this study is the demonstration that the positive homogeneity inherent to PWL systems simplifies the governing invariance equations. This property reduces the partial differential equations (PDEs) governing the invariant manifold's geometry to a more tractable system of ordinary differential equations (ODEs), which can then be solved numerically using periodic solvers. This approach separates the cone's geometry into its radial and directional components, allowing one to compute the NNMs as invariant cones for both conservative and damped systems and was implemented via two distinct parametrization techniques. 

First, a projection-based method analogous to the classical graph-style parametrization for smooth systems was developed. This technique uses a cylinder as an intersecting surface to fully characterize the generating set of the cone as a function of a pair of master coordinates. This intuitive graph-style approach provides a direct route to obtaining closed-form ROMs, but its limitations in handling complex geometries, particularly folding and multiple-crossing cones, were identified. 

To address these limitations, we derived a more robust arc-length parametrization. By using a unit hypersphere as the intersecting surface and the arc-length of the generating curve as its parameter, this method enables the computation of invariant cones with complex folding geometries. We demonstrated the effectiveness of this arc-length technique in a highly nonlinear setting with stiff unilateral contact forces, with both continuous and discontinuous vector fields. 

Ultimately, the ability to explicitly parametrize these invariant cones provides a rigorous basis for the nonlinear modal analysis of PWL systems. The resulting two-dimensional ROMs are shown to capture the system's behavior with high fidelity, even in challenging scenarios such as an internal resonance. In this context, our approach successfully parametrizes the resonant NNMs—which can be geometrically classified as multi-sheeted invariant cones—within the original two-dimensional reduced space, a task that would typically require increasing the number of master modes in smooth nonlinear MOR frameworks.

\section*{Conflict of interest}
The authors declare that they have no conflicts of interest.

\section*{Author Contributions}
Y.K. and R.L. contributed to the study conceptualization, methodology, and investigation. Y.K. was responsible for the formal analysis, validation, visualization, software, and writing the original draft. R.L. provided supervision and funding acquisition. Both authors reviewed and edited the manuscript. All authors read and approved the final manuscript.

\section*{Funding}
No funding has been received for this article.

%\begin{appendices}
%	\section{Appendix if needed}
%\end{appendices}

\bibliographystyle{abbrv}
\bibliography{bibliography_yassine}
\end{document}

%% file: FSP_1.tikz
% This file was created by matlab2tikz.
%
%The latest updates can be retrieved from
%  http://www.mathworks.com/matlabcentral/fileexchange/22022-matlab2tikz-matlab2tikz
%where you can also make suggestions and rate matlab2tikz.
%
\definecolor{mycolor1}{rgb}{0.12941,0.12941,0.12941}%
\begin{tikzpicture}

\begin{axis}[%
width=0.4\textwidth,
height=0.4\textwidth,
at={(0in,0in)},
scale only axis,
unbounded coords=jump,
xmin=0,
xmax=10,
xlabel style={font=\color{mycolor1}},
xlabel={$k_n/k$},
ymin=0.4,
ymax=0.65,
ylabel style={font=\color{mycolor1}},
ylabel={$\omega_1\sqrt{m/k}$},
axis background/.style={fill=white},
title style={font=\bfseries\color{mycolor1}},
title={Frequency-Stiffness Plot for Mode 1},
axis x line*=bottom,
axis y line*=left,
xmajorgrids,
ymajorgrids,
legend style={legend cell align=left, align=left}
]
\addplot [color=blue, line width=1.5pt]
  table[row sep=crcr]{%
0.00999999999999979	0.448042653700814\\
0.0501204819277117	0.459166789115519\\
0.0902409638554218	0.469023620600757\\
0.130361445783132	0.47783201539348\\
0.170481927710844	0.485757326300652\\
0.210602409638554	0.492928532717031\\
0.250722891566266	0.499448195385911\\
0.290843373493976	0.5053992656942\\
0.330963855421686	0.510849849804973\\
0.371084337349398	0.515856620459548\\
0.411204819277108	0.520467309318892\\
0.45132530120482	0.524722561343061\\
0.49144578313253	0.528657333929752\\
0.53156626506024	0.532301967762617\\
0.571686746987952	0.535683015353866\\
0.611807228915662	0.538823888372191\\
0.651927710843374	0.541745368020452\\
0.692048192771084	0.544466010685444\\
0.732168674698796	0.547002472745998\\
0.772289156626506	0.549369772842788\\
0.812409638554216	0.551581505082355\\
0.852530120481928	0.553650014064354\\
0.892650602409638	0.555586539768317\\
0.93277108433735	0.557401338721093\\
0.97289156626506	0.559103781795567\\
1.01301204819277	0.560702464780183\\
1.09325301204819	0.563619317671773\\
1.17349397590361	0.566207326541853\\
1.25373493975904	0.568513215491221\\
1.33397590361446	0.570576417508699\\
1.41421686746988	0.57243024808084\\
1.4944578313253	0.574102890999489\\
1.57469879518072	0.575618214081725\\
1.65493975903614	0.576996460366018\\
1.77530120481928	0.578843726413771\\
1.89566265060241	0.580468380936409\\
1.97590361445783	0.581446763511627\\
nan	nan\\
2.93879518072289	0.589169234173328\\
3.17951807228916	0.590396084137886\\
3.46036144578313	0.59161750690267\\
3.74120481927711	0.592657538036541\\
4.0621686746988	0.593667881593953\\
4.42325301204819	0.594621983290857\\
4.8244578313253	0.595501882681805\\
5.30590361445783	0.596362325191427\\
5.86759036144578	0.597161321267935\\
nan	nan\\
8.5155421686747	0.599296065907573\\
9.75927710843374	0.599838014243046\\
10	0.599925676570859\\
};
\addlegendentry{Stable NNM}

\addplot [color=blue, dashed, line width=1.5pt]
  table[row sep=crcr]{%
2.01602409638554	0.58190811648789\\
2.13638554216868	0.583192900470159\\
2.29686746987952	0.584705681226213\\
2.45734939759036	0.586029586569362\\
2.65795180722892	0.587469583574013\\
2.85855421686747	0.58871622171919\\
2.89867469879518	0.588945731489799\\
nan	nan\\
5.90771084337349	0.597211509177288\\
6.54963855421687	0.597914759831417\\
7.31192771084337	0.59855825478639\\
8.23469879518072	0.599147046356713\\
8.47542168674699	0.5992754955255\\
};
\addlegendentry{Unstable NNM}

\end{axis}

\begin{axis}[%
width=0.4\textwidth,
height=0.4\textwidth,
at={(0in,0in)},
scale only axis,
xmin=0,
xmax=1,
ymin=0,
ymax=1,
axis line style={draw=none},
ticks=none,
axis x line*=bottom,
axis y line*=left
]
\end{axis}
\end{tikzpicture}%

%% file: NNM1_damped_timehist_u.tikz
% This file was created by matlab2tikz.
%
%The latest updates can be retrieved from
%  http://www.mathworks.com/matlabcentral/fileexchange/22022-matlab2tikz-matlab2tikz
%where you can also make suggestions and rate matlab2tikz.
%
\begin{tikzpicture}

\begin{axis}[%
width=0.7\textwidth,
height=0.7\textwidth,
at={(0in,0in)},
scale only axis,
xmin=0,
xmax=100,
xlabel style={font=\color{white!15!black}},
xlabel={$t[s]$},
ymin=-2,
ymax=1,
ylabel style={font=\color{white!15!black}},
ylabel={$u$},
axis background/.style={fill=white},
xmajorgrids,
ymajorgrids
]
\addplot [color=black, line width=2.0pt]
  table[row sep=crcr]{%
0	1\\
0.0429179781207267	0.998943800381483\\
0.0876531023877902	0.995635081413496\\
0.132613125270979	0.990106316478702\\
0.178887001041204	0.982186690568383\\
0.224082537151645	0.972349498991235\\
0.276739425615418	0.958380540532175\\
0.327640238835883	0.942436218110657\\
0.384577817758171	0.921919797056816\\
0.438591092458182	0.900005629577748\\
0.495700968803661	0.874425706421249\\
0.556298772236531	0.844793991562014\\
0.618121467238936	0.812167182299248\\
0.687970632017908	0.772713693156149\\
0.766850051550094	0.725302222967315\\
0.854962703736561	0.669380412268183\\
0.950641283479527	0.605923995475592\\
1.07717191785697	0.519091995438075\\
1.48002664367016	0.240179772770986\\
1.59494598597641	0.164649052612077\\
1.69157273236482	0.103735709854291\\
1.78575855741971	0.0470278305779743\\
1.87913084493682	-0.00631333704932047\\
1.96443136900199	-0.0524810599927292\\
2.05729388909172	-0.100327940988905\\
2.1647928590161	-0.153157829222479\\
2.29509093814674	-0.214419641003076\\
2.45907263249016	-0.288742199230484\\
3.07363733057969	-0.564709327559299\\
3.26632639934091	-0.655244627832602\\
3.46872384400864	-0.752927048504716\\
3.75142632129123	-0.89219858776039\\
4.05648930926135	-1.04192258287657\\
4.22185115164976	-1.12057660367859\\
4.35439401082832	-1.18132603907006\\
4.47268851512418	-1.23327671096743\\
4.58406153467203	-1.27982946133901\\
4.68990018215618	-1.32162401137784\\
4.79189169370937	-1.35937217559562\\
4.88510792275925	-1.39149424478687\\
4.97404130307171	-1.41985922084399\\
5.05958849080928	-1.44491039511986\\
5.148491556935	-1.46849992703403\\
5.23079225216442	-1.48801807216452\\
5.31507585230534	-1.50560794050008\\
5.39786161428046	-1.52044918187973\\
5.47243476127581	-1.53169963355047\\
5.55150719455588	-1.54139195963819\\
5.63193114141329	-1.54885166293582\\
5.71377370419698	-1.55393066600999\\
5.79255783626265	-1.55640726887179\\
5.87347498466966	-1.55647760322911\\
5.95010145048916	-1.55423315104974\\
6.03210154991282	-1.5493494416954\\
6.11098216540007	-1.54224642962276\\
6.19264031855245	-1.53243440672158\\
6.27694194362064	-1.51971694799565\\
6.35164108424243	-1.50628844751114\\
6.42977733771312	-1.49011391178696\\
6.51154325268568	-1.47091133779401\\
6.59993132087146	-1.4476069364366\\
6.67960205135752	-1.42440456865178\\
6.76100495969072	-1.39862160110657\\
6.85815010907082	-1.36521539366487\\
6.95508332417278	-1.32915540841373\\
7.04623420013854	-1.29289372755311\\
7.14671595620902	-1.25043334514642\\
7.25693290225985	-1.2010627072153\\
7.36445796478178	-1.15029830883132\\
7.48371887491543	-1.09126265046694\\
7.60227085182052	-1.03002403032555\\
7.73414582360815	-0.959276927278452\\
7.88039824187942	-0.878057742446572\\
8.04262288617838	-0.785208308166915\\
8.23720793345372	-0.670973427670248\\
8.49924430198257	-0.514129116825856\\
9.05214554117103	-0.182346327385275\\
9.26305215127614	-0.0587849794321187\\
9.44809154592301	0.046994175659961\\
9.56707259245404	0.112477985991106\\
9.65939711247918	0.160872372519123\\
9.74187161596134	0.201781779263555\\
9.82432936742369	0.240103890386038\\
9.89713613879708	0.271536839560611\\
9.96950651094174	0.300348917663868\\
10.0414499032163	0.326431445319429\\
10.1129587918965	0.349704021024934\\
10.1840060060367	0.37011493981062\\
10.2457576121026	0.385609531505636\\
10.3070726304038	0.398901064337096\\
10.367820509774	0.410002101492765\\
10.4364398559499	0.420070964686843\\
10.5038711317932	0.427435233604868\\
10.5695468935503	0.432236052892222\\
10.6323101986165	0.434686326647324\\
10.7016506471287	0.435036032831647\\
10.7656175837525	0.433242760081953\\
10.8384546656183	0.428838916984816\\
10.9108722567634	0.42209919416652\\
10.9887434206276	0.412395827047277\\
11.0703148210941	0.399713591195123\\
11.1527986489022	0.384512098460007\\
11.2429475930533	0.365467550444933\\
11.3349773191706	0.343755543312156\\
11.4434889973284	0.3157077860795\\
11.5649483309714	0.281873083694961\\
11.7189244078566	0.23649744507334\\
12.2482678472872	0.0782844632574182\\
12.3877988955676	0.0399806163564591\\
12.5229986396988	0.00531413420095816\\
12.6520635032856	-0.0252995443454438\\
12.7998004044376	-0.0578192176263741\\
12.9770216919665	-0.0943236359666457\\
13.2359329365367	-0.145029140133715\\
13.7220627005529	-0.240034356910812\\
14.0052338901347	-0.298075779644577\\
14.3317499329943	-0.36749867455174\\
14.8399289220187	-0.47563900585132\\
15.0512814485483	-0.518009527626873\\
15.2252013288836	-0.550600301299824\\
15.3857062487773	-0.578306207676974\\
15.5284945673709	-0.600676478152309\\
15.6640479653234	-0.619661178448183\\
15.7945703904402	-0.635662634185948\\
15.9222954371799	-0.648992648808743\\
16.0449844932917	-0.65950285768298\\
16.1662655660578	-0.66758699776129\\
16.285813050736	-0.673242134368124\\
16.4082482764793	-0.676602681908051\\
16.5275593990663	-0.67748250147315\\
16.6420096459759	-0.676097578976169\\
16.7598963847036	-0.672399645612259\\
16.8760249953672	-0.666529888928025\\
16.9977401610059	-0.658051498865348\\
17.1139794739473	-0.647790662654003\\
17.2398105460131	-0.634382318791623\\
17.3652254609559	-0.618738912101477\\
17.4951686064583	-0.600254701839503\\
17.6244655928686	-0.579702779297847\\
17.7604936269762	-0.555921454289233\\
17.9140523806317	-0.52664325859655\\
18.0752327416092	-0.493428779244127\\
18.2468765351724	-0.455623163349912\\
18.4315162455647	-0.41259512378231\\
18.6454236014575	-0.360311682810618\\
18.9076799634629	-0.293660437398827\\
19.3071357490826	-0.189374786659656\\
19.7957178111621	-0.0622281177238051\\
20.0834684130136	0.0102537731627166\\
20.2584779831064	0.0523341642113309\\
20.3895791286438	0.081638808215061\\
20.5087975575114	0.105966849816653\\
20.6177742343014	0.125885568665552\\
20.7257826471239	0.143173295658499\\
20.8239375907692	0.156603422805532\\
20.9211531957172	0.167666057716019\\
21.017259695615	0.176361612589929\\
21.111880567765	0.182734652883653\\
21.2124682180026	0.18715466401791\\
21.3154004073935	0.189232464184556\\
21.4176000450966	0.18895361055219\\
21.5231135583835	0.186363534548406\\
21.6301403278429	0.181530985977872\\
21.7459480364346	0.174047515331907\\
21.8753099560895	0.163273462594461\\
22.0120922865486	0.149577008289469\\
22.1720372948588	0.131249812941093\\
22.3810209820357	0.104855203801847\\
23.0125876802413	0.0235021745668575\\
23.2155548915146	0.000638309934444692\\
23.4081313296167	-0.0187215539237826\\
23.6555326856127	-0.0411919367595885\\
24.137374226974	-0.0820079134058744\\
24.5937074244352	-0.121725619936285\\
25.0785067675977	-0.166481145679612\\
25.6242098233373	-0.216507967954641\\
25.9039024712856	-0.239764418556263\\
26.1379270585168	-0.256941150131908\\
26.3426491093326	-0.269775277864312\\
26.536489072769	-0.279724014734114\\
26.7237731793823	-0.287079198872362\\
26.9041560567497	-0.291916660292245\\
27.0870233551979	-0.29447195242011\\
27.2666512880687	-0.294630779451495\\
27.440846468243	-0.292555241231256\\
27.6172843621823	-0.288249892083286\\
27.7955557666079	-0.281719380759469\\
27.9828253555104	-0.272613645223217\\
28.1780940387897	-0.260833194468177\\
28.3838697831427	-0.246118009085251\\
28.6094171988724	-0.227591111798986\\
28.8494312594798	-0.205535267955142\\
29.1263624693986	-0.177690131917089\\
29.4638467042152	-0.141272804399222\\
29.9739936499624	-0.0835606793547612\\
30.6419556825287	-0.00837975470248864\\
30.9564034229524	0.0248902207350739\\
31.1451237494427	0.0426368008258891\\
31.3087587601427	0.055830482218397\\
31.4613078049308	0.0659182582664926\\
31.6117638753741	0.0735687334576625\\
31.75937971764	0.0787658494740242\\
31.9102930556971	0.0817197757785948\\
32.062940881106	0.0823751305500053\\
32.2227704737889	0.080735879761022\\
32.3948163163389	0.0766202824891735\\
32.5839266206675	0.0697574498824451\\
32.8090129849189	0.0591930941766918\\
33.1199874661793	0.0420879797635649\\
33.7379293661362	0.00775509091506876\\
34.0275661452184	-0.00572699393350717\\
34.3895251779268	-0.0200791460341208\\
35.7097106564491	-0.0705902609448685\\
36.4731838827147	-0.100582926898895\\
36.8413547343797	-0.112599238276275\\
37.1507247039394	-0.120445649252972\\
37.439805987737	-0.125477283667621\\
37.7084641418587	-0.127927611996256\\
37.9773883445611	-0.128122439642809\\
38.2428277691139	-0.126089658673777\\
38.5234547270388	-0.121621511382571\\
38.8185116815438	-0.11454526192162\\
39.1299005042242	-0.104747631886823\\
39.4833720474739	-0.0912345622528363\\
39.8930344196982	-0.0732105656711042\\
40.4799424984277	-0.0448547868944047\\
41.5677388787001	0.00793531276192994\\
41.8653180568733	0.0201711743566335\\
42.1091088967912	0.0279275837336428\\
42.3388570376256	0.0328989427470248\\
42.5605368695183	0.0354117300182963\\
42.7899630394696	0.0357324879726093\\
43.0378835947976	0.0337732023477173\\
43.3301113698048	0.0290454718461888\\
43.7165782408576	0.0203083546687708\\
44.631673946359	-0.00111073248879734\\
45.1753735594495	-0.0104982659302379\\
46.5730710104651	-0.0340280600121474\\
47.3814003612927	-0.0471886614374313\\
47.8644311642797	-0.0527537854618174\\
48.2949559535362	-0.055423501074813\\
48.7078858834572	-0.0556572206898096\\
49.1244260341441	-0.0535544138292181\\
49.5664311486241	-0.048963291153612\\
50.0591214345091	-0.0414857895064529\\
50.6826215517841	-0.0295642912758325\\
52.7220975892422	0.0113714032188028\\
53.0667410624933	0.0146557869031625\\
53.4128482920398	0.0156038971033183\\
53.7868688627275	0.0142828236485002\\
54.2764492299005	0.0100702989697368\\
55.6965703861725	-0.00346996749240702\\
58.6169054138633	-0.0232363504023851\\
59.2509167613692	-0.0243005415144069\\
59.8743642357545	-0.0230458029089959\\
60.5650027447877	-0.0192888763210846\\
61.4604210007209	-0.0119551543089216\\
63.4963397271365	0.00545458023211154\\
64.0088429221257	0.00678309427037505\\
64.5811617651813	0.00586245328763368\\
65.619705505132	0.00127944708857797\\
66.7960726102114	-0.00279396071790927\\
69.7295999790906	-0.0105406533397883\\
70.668566209443	-0.00982128934813886\\
71.7908126616624	-0.00653453455748831\\
74.6437245677502	0.00294124894833203\\
75.5384254616995	0.00208278729476774\\
79.2036120412269	-0.00363741094612635\\
80.8346562499446	-0.00456452752636949\\
82.401215787724	-0.00295010185972444\\
85.7370404225949	0.00121511633798832\\
93.4805212808885	-0.000993316831937818\\
96.8597003226148	0.000405846122660591\\
100.009601964004	-0.000545601840798327\\
};

\addplot [color=red, dashed, line width=2.0pt]
  table[row sep=crcr]{%
0	1\\
0.0446302474218356	0.998858247804293\\
0.0839415695318735	0.995993836291134\\
0.130886516219533	0.990358886444724\\
0.179019350880324	0.982161657466136\\
0.226529419477501	0.971760605994731\\
0.281962119217482	0.956854788426384\\
0.336023422943356	0.939591713643011\\
0.390906460347679	0.919475470288191\\
0.447424961204348	0.896210075788503\\
0.504535587606568	0.870264986986783\\
0.561314576476065	0.842240045002654\\
0.628693252469247	0.806371723302036\\
0.695921078862142	0.768073101794826\\
0.775316031203857	0.720064365238258\\
0.856168081806629	0.668608438527599\\
0.961618177130688	0.598514394819361\\
1.09155755399202	0.509102041515547\\
1.45659215726027	0.255951131455063\\
1.57708834402574	0.17620175996565\\
1.67818996687036	0.112034806415807\\
1.77244785540861	0.0548911727370012\\
1.86830909828848	-0.000265438033608234\\
1.95051905806814	-0.0450876304079344\\
2.04957427254013	-0.0964173486374023\\
2.1549327731757	-0.148392649487377\\
2.27457052602898	-0.204916699772895\\
2.42492712443499	-0.273423683645404\\
2.71340855306838	-0.401717136671579\\
2.96169275878323	-0.513275596011695\\
3.16177415514119	-0.605805036619159\\
3.34875986008758	-0.694744234738224\\
3.56871625033851	-0.801941362197084\\
4.13939644195257	-1.08172379934953\\
4.28697460347971	-1.15077304605047\\
4.4174628258384	-1.20937107501359\\
4.53117225791136	-1.25807799212605\\
4.6468092621765	-1.30495368225272\\
4.74529306051694	-1.34248138812306\\
4.84536913671406	-1.37811756719155\\
4.94667628385383	-1.41140449541605\\
5.02822616991892	-1.43601876849618\\
5.10980704454425	-1.45858538462755\\
5.19111141181499	-1.47893541335775\\
5.27195867777402	-1.49696927150947\\
5.35231781765863	-1.51264573526491\\
5.43229941694835	-1.52596031874273\\
5.51213869653539	-1.53692234170445\\
5.5921904966963	-1.54553505579659\\
5.67296115157149	-1.55177818121608\\
5.75524838139997	-1.55558693126183\\
5.84080932255318	-1.55680118389058\\
5.91671754028664	-1.55553135300711\\
5.99433088995094	-1.55195619903421\\
6.0643268626655	-1.54676772421286\\
6.13545954566196	-1.53960134716048\\
6.20956217283715	-1.53012653915494\\
6.28800974005037	-1.51789198235213\\
6.37217724316089	-1.5022877029507\\
6.45411550514569	-1.48468453978536\\
6.53439452255689	-1.4651862881736\\
6.61304740914879	-1.44398063628319\\
6.69010727867565	-1.42124694020275\\
6.79002167078731	-1.3889907168431\\
6.88761399970197	-1.35458471298692\\
6.98346498554535	-1.31814644660231\\
7.07816829325733	-1.27971239578257\\
7.17230165713346	-1.23925844411876\\
7.26641083890594	-1.19671783083247\\
7.38484139494106	-1.14043244398256\\
7.50461380986746	-1.08067472550263\\
7.62567586056521	-1.01768956651081\\
7.77097834341936	-0.939106095846441\\
7.91325629856742	-0.859491089538238\\
8.07445197988353	-0.766731030767772\\
8.26261294855296	-0.655904604433104\\
8.55291368245732	-0.481805208761827\\
9.02568292392131	-0.198035087776191\\
9.25276412965977	-0.064752786864716\\
9.43517372173845	0.0397224074810509\\
9.55412017183833	0.105498571632253\\
9.65513292214463	0.158693572573625\\
9.73850413067657	0.20016051440075\\
9.81941499659828	0.237900391723002\\
9.89596866953895	0.271054858178843\\
9.96563177700162	0.298874573522937\\
10.0393633128856	0.325717122558061\\
10.1092348968962	0.348564818710827\\
10.1775478335229	0.368379806322821\\
10.2377567696795	0.383728184764152\\
10.3080595319825	0.39910582333917\\
10.3732372109423	0.410900619004451\\
10.4432583612457	0.420938234177484\\
10.5058365588367	0.427622240459257\\
10.5720119234454	0.432381384418861\\
10.6362688690734	0.434782604311906\\
10.7048281918395	0.435004937050977\\
10.7742329860252	0.432861687183731\\
10.8463300981575	0.428228695422476\\
10.9172416541362	0.421410065157787\\
10.9933583570301	0.411756417366007\\
11.0700941175111	0.399762541123479\\
11.1493648432598	0.385200613266278\\
11.2403999021487	0.366048911441453\\
11.3413343904832	0.34219133790684\\
11.442975732507	0.315856213752241\\
11.5703711760048	0.280324523826479\\
11.7351104890523	0.231632762108731\\
12.2292652932608	0.0836795999082369\\
12.3735865573434	0.0437758543642843\\
12.5053016928087	0.00970489868998925\\
12.6373763433154	-0.0219288629440371\\
12.777194379806	-0.0529867755285807\\
12.9447949066236	-0.0878378914762408\\
13.1788468473984	-0.133987537044973\\
13.7491130864221	-0.24548090665499\\
14.0325817968669	-0.303811112407587\\
14.3830739137789	-0.378544139459748\\
14.8223736665815	-0.472031458552877\\
15.0257273554722	-0.513053400321823\\
15.1952747828336	-0.545189528057833\\
15.3498286551331	-0.572351321257955\\
15.4885868622828	-0.594670832470484\\
15.6303122674297	-0.615165902026604\\
15.7526738926686	-0.630795400911367\\
15.8747179684442	-0.644324553438437\\
15.9957445720628	-0.655585448092069\\
16.115804723953	-0.664534217316202\\
16.2356478886435	-0.671180350518384\\
16.3565963493751	-0.675513464395863\\
16.4815553245735	-0.677446514479513\\
16.5933434475112	-0.676971460782795\\
16.7096207266164	-0.674272828706762\\
16.8273857892714	-0.669269455669465\\
16.9552641053073	-0.661292600692377\\
17.0660499397465	-0.652291084213815\\
17.201319387069	-0.638752819675986\\
17.3317147640681	-0.623157669867183\\
17.457632995305	-0.605846250191149\\
17.6034133275066	-0.583209381782396\\
17.7461937497387	-0.558536415945071\\
17.887401490442	-0.531911485175499\\
18.0527920539664	-0.498205404378538\\
18.2206140593556	-0.461567192407145\\
18.4148229206235	-0.416585915043996\\
18.62792101461	-0.364681558101381\\
18.8812578609788	-0.300477129344145\\
19.267043807505	-0.199903381598702\\
19.8023630090781	-0.0605287491935229\\
20.0890384662928	0.01162866958515\\
20.2629234953583	0.0533662318431709\\
20.3958591452639	0.0829823460522903\\
20.5172132570857	0.107590995592403\\
20.6299738813728	0.127968173488412\\
20.728277488826	0.143545887859602\\
20.8309813267392	0.157484856837172\\
20.9262888848586	0.168191757041043\\
21.0281752329639	0.177212097343741\\
21.1201034480062	0.183191263940031\\
21.2166989744633	0.187293090374226\\
21.3122162772645	0.189209260682958\\
21.4157450191984	0.188983624292746\\
21.5192805917357	0.186501459394762\\
21.634913452027	0.18127178921155\\
21.7549669942126	0.17337904302893\\
21.8801956530088	0.162826358245198\\
22.0263254092472	0.148040318291663\\
22.1973814625017	0.128171788338662\\
22.4206310560787	0.0996706353354995\\
22.9571942243275	0.0301529040220885\\
23.1668700780264	0.00589091076724912\\
23.3614058001489	-0.0142098827808468\\
23.5899769422345	-0.0354264552460251\\
23.9563643610557	-0.0667770041491451\\
24.5402897366526	-0.116945728673883\\
25.0024032907311	-0.159354279124088\\
25.6264930201617	-0.216717760019122\\
25.9085051351144	-0.240134463920555\\
26.1435495324185	-0.25732921449827\\
26.3458843556241	-0.26996563814788\\
26.5495493568191	-0.280319772099773\\
26.7307358972983	-0.287319180839816\\
26.9104712329555	-0.292057202079789\\
27.0930850789245	-0.294526239810949\\
27.2726207518625	-0.294603961375984\\
27.4497831432406	-0.292396366019759\\
27.6241087075828	-0.288046181056799\\
27.8154235820642	-0.280868227171624\\
27.9998468316275	-0.271686254284816\\
28.1995403768953	-0.259414316714924\\
28.4149148711514	-0.243715524550822\\
28.650421473541	-0.223986267449774\\
28.8893243931806	-0.201671813711371\\
29.1793183772104	-0.172133290800105\\
29.5270384071032	-0.134244729700328\\
30.1490577289955	-0.0635849210431587\\
30.6859564411464	-0.00358114578651225\\
30.9796801951424	0.0272025234586124\\
31.169732032385	0.0447662328309093\\
31.3278355563027	0.0572166703999812\\
31.4818489632884	0.0671027985810326\\
31.6354213607148	0.0745588294514192\\
31.7802228711888	0.0793169416416788\\
31.9327360802344	0.0819621513570468\\
32.0904603409508	0.0822580645723292\\
32.2527897619263	0.0801852160650611\\
32.428723807272	0.0755563450433669\\
32.6216315902528	0.0681468878306788\\
32.859313529069	0.0565729852253298\\
33.2141219615172	0.0366666012103138\\
33.7027454946065	0.00954469956022308\\
34.0026922276793	-0.00465668698775801\\
34.3590476583624	-0.0189348297581802\\
35.8015749545912	-0.0743374544116193\\
36.4640720666682	-0.100261298018069\\
36.8327688033107	-0.112352227834947\\
37.1373686038531	-0.12016020149143\\
37.4200067773432	-0.125215029413368\\
37.7005992659205	-0.127892792295725\\
37.9589994014911	-0.128185243779001\\
38.2187806922134	-0.126366079664592\\
38.4885693585027	-0.122304940012768\\
38.7742183737506	-0.11575661265428\\
39.0884551872181	-0.106178552078376\\
39.442692121515	-0.0929015650864784\\
39.852891995803	-0.0750639646753797\\
40.4010076163496	-0.0487605001127207\\
41.603564839299	0.00953121710203675\\
41.8860871303508	0.020923530578429\\
42.1210616151029	0.02824575015579\\
42.3427274683718	0.032963281935551\\
42.5637000605813	0.035432572598026\\
42.7893106034854	0.0357355470788292\\
43.0430385040919	0.0337108983551389\\
43.3356606314785	0.028936407818108\\
43.7274428531406	0.0200407345999736\\
44.6197545726124	-0.000879265999316203\\
45.1496110214289	-0.0100852714536614\\
46.5976604418026	-0.0344634600813833\\
47.3812014939086	-0.047187694318211\\
47.8656825269712	-0.052766636032743\\
48.2871827634449	-0.0553992178319476\\
48.6865207332589	-0.0557038428586907\\
49.0921951224694	-0.0537991867928866\\
49.5334103566588	-0.0493836524046998\\
50.0320923480936	-0.0419497539569562\\
50.6314492144682	-0.0306134011301822\\
52.0511121519647	-0.000443342862553209\\
52.5916806477005	0.00954668957781735\\
52.94482290119	0.0137616560624849\\
53.2893917253336	0.0155288031446617\\
53.652452125664	0.0150004418950687\\
54.0884665637067	0.0119155000746929\\
55.0432833147847	0.00193640958862318\\
55.7571833030587	-0.00390173670598415\\
58.6690012972537	-0.0234111038724052\\
59.2963668361664	-0.0242865558687129\\
59.9198934984813	-0.022869280033035\\
60.6281382588088	-0.0188451638497042\\
61.5749396094973	-0.0109060573461335\\
63.4350622455091	0.00514563901256793\\
63.9557002540862	0.00674960558197313\\
64.5080226548083	0.00609082852692211\\
65.3866311461572	0.00234926927581114\\
66.5817566821611	-0.00214789156783013\\
69.7338803588589	-0.0105429304224174\\
70.6759491701379	-0.00980741945895147\\
71.8094386746304	-0.00646539268788615\\
74.6312129570287	0.00293692588380168\\
75.5119511487771	0.00213177092568628\\
79.3790552834905	-0.00385173460719557\\
80.9397191610302	-0.00452760812153485\\
82.5572835645984	-0.00269779202078269\\
85.6327548961344	0.00125315001044157\\
88.6244181038461	-0.000808462018554224\\
92.135376303254	-0.00181670660752786\\
98.763756614515	-0.000214657607259028\\
100.002997367862	-0.000543812658492016\\
};

\end{axis}

\begin{axis}[%
width=0.7\textwidth,
height=0.7\textwidth,
at={(0in,0in)},
scale only axis,
xmin=0,
xmax=1,
ymin=0,
ymax=1,
axis line style={draw=none},
ticks=none,
axis x line*=bottom,
axis y line*=left
]
\end{axis}
\end{tikzpicture}%

%% file: NNM1_damped_timehist_v.tikz
% This file was created by matlab2tikz.
%
%The latest updates can be retrieved from
%  http://www.mathworks.com/matlabcentral/fileexchange/22022-matlab2tikz-matlab2tikz
%where you can also make suggestions and rate matlab2tikz.
%
\begin{tikzpicture}

\begin{axis}[%
width=0.7\textwidth,
height=0.7\textwidth,
at={(0in,0in)},
scale only axis,
xmin=0,
xmax=100,
xlabel style={font=\color{white!15!black}},
xlabel={$t[s]$},
ymin=-0.8,
ymax=0.8,
ylabel style={font=\color{white!15!black}},
ylabel={$v$},
axis background/.style={fill=white},
xmajorgrids,
ymajorgrids,
legend style={legend cell align=left, align=left, font=\tiny, draw=white!15!black}
]
\addplot [color=black, line width=2.0pt]
  table[row sep=crcr]{%
0	0\\
0.0496254766668471	-0.056547613477349\\
0.0965620328627068	-0.108373863108014\\
0.141937426254472	-0.156817698757138\\
0.188037487451069	-0.204252452546712\\
0.233047391581877	-0.248734509624086\\
0.276739425615418	-0.29010293274203\\
0.319257394974059	-0.328581316905968\\
0.360531077188838	-0.364203741908966\\
0.400455226830474	-0.396999579853826\\
0.438591092458182	-0.426770525779261\\
0.47491265122936	-0.453690459496869\\
0.515680687647219	-0.482219400699563\\
0.550940725775817	-0.505443602619039\\
0.587716572585691	-0.52822545784899\\
0.624554433908415	-0.549567506283125\\
0.658644214969897	-0.567999624429433\\
0.695678088942344	-0.586591686094522\\
0.734911282811197	-0.604669438059005\\
0.766850051550094	-0.61816541587956\\
0.804658374526369	-0.63273898613609\\
0.838352922374213	-0.644459545556671\\
0.871773621795356	-0.654920509729735\\
0.908483435964769	-0.665094679261699\\
0.943338207750074	-0.673501210485128\\
0.980537289920136	-0.681152849259675\\
1.01963561794548	-0.687758602831067\\
1.05214430429024	-0.692156924068527\\
1.08559030319641	-0.695671375702204\\
1.11988492873149	-0.698238184629517\\
1.15496074757122	-0.699808413716795\\
1.1907663421014	-0.700345081774842\\
1.22726276464884	-0.69982107408012\\
1.26442107273169	-0.698217572545644\\
1.30222058903436	-0.695522839104939\\
1.34064766876818	-0.691731244766046\\
1.3796948387464	-0.686842471637064\\
1.41936022133244	-0.680860836742724\\
1.46975785378137	-0.67187144395831\\
1.52126365963714	-0.661188710931768\\
1.57375079225771	-0.648872811137053\\
1.62699343949443	-0.635036680326891\\
1.68072285077311	-0.619837585292387\\
1.7452856258251	-0.600120084573334\\
1.81379608683682	-0.577701135443803\\
1.94382046527875	-0.534184711998989\\
1.99360599695908	-0.51971749894787\\
2.04054344964032	-0.507388983090507\\
2.09157820493159	-0.495388559406834\\
2.13656723420259	-0.485985896679082\\
2.1838667662021	-0.477247538845205\\
2.23387511044673	-0.46923828230905\\
2.28497571577329	-0.462301362921025\\
2.33498459735092	-0.456672070764569\\
2.39334469630471	-0.4514683011447\\
2.44978000630346	-0.447741987539658\\
2.51077685025507	-0.445045165371781\\
2.57018202528387	-0.443631158129818\\
2.63849277803156	-0.443330288680841\\
2.70801107551924	-0.444301414948583\\
2.78647241602893	-0.446708615590182\\
2.8686534426278	-0.450436824481784\\
2.97187387161912	-0.456406431812553\\
3.11922129500046	-0.466328383890826\\
3.33385848315307	-0.480774758785827\\
3.44624411032483	-0.487055548890169\\
3.53620481509786	-0.490959968719437\\
3.62644001708023	-0.493610173887788\\
3.70582306410112	-0.49472392751143\\
3.78578807088589	-0.494559579111126\\
3.85492989188302	-0.493291621105158\\
3.92497971847229	-0.490876798524951\\
3.99601000360221	-0.487213288305\\
4.06875818187656	-0.482147879236095\\
4.14351445547499	-0.475521593106947\\
4.20858544139516	-0.468561247208214\\
4.27489036566659	-0.460316883340752\\
4.34114504791479	-0.450914753589018\\
4.40738967591572	-0.440355225223414\\
4.47268851512418	-0.428823270211808\\
4.54268879933966	-0.415241799848999\\
4.61285312205716	-0.400389898843912\\
4.68196364506581	-0.384582839086093\\
4.75670297905198	-0.366218532768997\\
4.83283838649476	-0.346212814874704\\
4.90876323528327	-0.325025157304083\\
4.99083512572079	-0.300818800193781\\
5.07712821335373	-0.274016484596984\\
5.16661054047553	-0.244894453175078\\
5.25870787849742	-0.213663648292595\\
5.36265575882371	-0.177092210626896\\
5.48223592099379	-0.133601841932688\\
5.62179150937025	-0.0814258130729115\\
5.82812559600742	-0.00274321337748518\\
6.09957849785388	0.100620096829147\\
6.24042761062999	0.152946898779419\\
6.36452822437403	0.197762636726949\\
6.47015958616451	0.234708227491467\\
6.56929347849534	0.268210825220336\\
6.66325775532087	0.298799572921112\\
6.76100495969072	0.329307538421759\\
6.85815010907082	0.358208087119834\\
6.95508332417278	0.385554925829709\\
7.04623420013854	0.409852747124532\\
7.13276951032304	0.431607497907834\\
7.21597668044805	0.451291928738158\\
7.29753740795188	0.469394413937181\\
7.37776475159279	0.486037618808737\\
7.45733172646904	0.501398453612552\\
7.5364443639768	0.515539851874777\\
7.61543356756344	0.528536990023795\\
7.69451693376031	0.54043406803082\\
7.7738512275597	0.551258993679667\\
7.85369762470461	0.561046685758583\\
7.93410393369774	0.569799513742765\\
8.01539273717871	0.577545149703582\\
8.09753435050548	0.584271032470681\\
8.18093266464041	0.589996784682796\\
8.26549857609669	0.594700455853356\\
8.35174866961142	0.598390699551274\\
8.4395483846034	0.601038536935519\\
8.52958462101553	0.602635608077676\\
8.62169693519691	0.603144968789664\\
8.71686381541069	0.602528343659685\\
8.81495066912719	0.600733432044166\\
8.91301525325404	0.597819628631427\\
9.0220704529075	0.593333339197727\\
9.12551802987898	0.58793041691824\\
9.23670240182379	0.580948276209369\\
9.35905128593696	0.571940530567488\\
9.39746317844501	0.568397470318757\\
9.43604864136863	0.563482811802956\\
9.47406815339865	0.557333184622962\\
9.51509649621067	0.549288451773521\\
9.55142196660637	0.540988533274515\\
9.59107258149581	0.530716360161847\\
9.63322916079397	0.518467567634374\\
9.67759839402767	0.504173862531374\\
9.72346192794375	0.487976708924904\\
9.76943842764653	0.470385895508613\\
9.81520043711545	0.451631773934082\\
9.86988916912536	0.42773321303622\\
9.92432868370531	0.402503704816894\\
9.97852609116256	0.376136345078393\\
10.0414499032163	0.344192800920183\\
10.1129587918965	0.306489864858392\\
10.2017067104422	0.258203050224196\\
10.3418693377002	0.180269629465116\\
10.4871471277943	0.0998973277432071\\
10.5776528435519	0.0512441493689693\\
10.6545334218665	0.011308070074378\\
10.7196133165868	-0.0212451170925618\\
10.782379090753	-0.0513803774854011\\
10.8476926152944	-0.0812755015334687\\
10.9020063694175	-0.104903205728093\\
10.9546312088909	-0.126660585220762\\
11.0137654385872	-0.1497103769604\\
11.0703148210941	-0.170313657759081\\
11.1238867390572	-0.188495142510774\\
11.1733814464683	-0.204112495593051\\
11.2231678371909	-0.218662308150201\\
11.2756286638053	-0.232725102900602\\
11.3278586404017	-0.245428945012961\\
11.3792716961565	-0.256672701538037\\
11.4271313422553	-0.26602214483853\\
11.481090153285	-0.275285757325335\\
11.5368123729578	-0.283454173809588\\
11.5907782118686	-0.290040244277435\\
11.6406578406196	-0.294997846134549\\
11.6947618482527	-0.299184621559547\\
11.7521075758246	-0.302316707053677\\
11.8123005285442	-0.304219755756449\\
11.8656820110602	-0.304776716573429\\
11.9206264426233	-0.304297096670894\\
11.9865977954781	-0.30239316391787\\
12.0544496461219	-0.299029399262565\\
12.1241582758332	-0.294213349849599\\
12.1959068761564	-0.287953950743685\\
12.2695530741398	-0.280304490965563\\
12.3553061442997	-0.270038060644026\\
12.4408104007927	-0.258571650100791\\
12.5456813459881	-0.24316736800138\\
12.6109714463748	-0.233776733345493\\
12.6819544972851	-0.224815353330357\\
12.7553339297544	-0.216872964691277\\
12.8274029795661	-0.210323850004471\\
12.9056726498952	-0.204540298566116\\
12.9871135261745	-0.199894231431827\\
13.0655336499304	-0.196637725729786\\
13.149376438113	-0.194352679994211\\
13.2443292080464	-0.193086642536144\\
13.3442381592602	-0.193050851445662\\
13.4506065010168	-0.194202933745885\\
13.5803781959724	-0.196851077734294\\
13.7563429019607	-0.201747206799439\\
14.0951524776331	-0.211346696720611\\
14.2413732233566	-0.214095600326559\\
14.365775135841	-0.215249397378926\\
14.4800143117225	-0.215141926312043\\
14.5961027438873	-0.213739089807277\\
14.7031624756276	-0.211192781186242\\
14.8017184206126	-0.207730268562415\\
14.9054985265473	-0.202890370917586\\
15.0115350190942	-0.196664019362899\\
15.1175353908777	-0.189148512338477\\
15.2252013288836	-0.180216116511062\\
15.3301944722071	-0.170284503258799\\
15.442162502515	-0.158426406901299\\
15.5511301048346	-0.145707947257449\\
15.6725533299434	-0.130278224585297\\
15.8034805289419	-0.112313687405958\\
15.9409320920321	-0.0921817384084562\\
16.0978635603347	-0.0679105875726833\\
16.285813050736	-0.0375485786971979\\
16.6321862222215	0.019987681232692\\
16.8878335258621	0.0618486301890613\\
17.0743278436357	0.0911348865194554\\
17.2252688052145	0.113696657739538\\
17.3652254609559	0.133504285009096\\
17.5113655995257	0.152888622610888\\
17.6406155031227	0.16882473928041\\
17.7747113311295	0.184082243313753\\
17.9003691079173	0.19715057722631\\
18.0219573445941	0.208634227310668\\
18.141451691162	0.218792976636749\\
18.2600396058077	0.227767552049755\\
18.3918321728007	0.236456809647109\\
18.5245368612916	0.24386221743454\\
18.6589513223987	0.250023649789583\\
18.795847410858	0.254960721979472\\
18.9360149884968	0.258675872236708\\
19.0803636542602	0.261156095774552\\
19.2299982450193	0.262369497063816\\
19.3703030981654	0.262331649082171\\
19.5174448731892	0.261153461120855\\
19.6768133631088	0.258659165741776\\
19.8383753831237	0.254940093805416\\
20.0110216959389	0.249753318979074\\
20.0885481903752	0.246647392267761\\
20.1431079123918	0.243148111864883\\
20.1972610123801	0.238553730981067\\
20.2584779831064	0.232084245569482\\
20.317057288148	0.22470275978462\\
20.3803408743255	0.215515605714927\\
20.444753472238	0.204981099443103\\
20.5179134284612	0.191720168427139\\
20.5906186669511	0.177350191053677\\
20.6718985096436	0.160106037406564\\
20.7704959283407	0.137876906876855\\
20.8947491406523	0.108488086389727\\
21.2691393682124	0.0189510823439178\\
21.3770037117648	-0.00506106788790817\\
21.4677590753421	-0.0240922014272655\\
21.5595378182802	-0.042048966590329\\
21.6473125762766	-0.0578674542338149\\
21.7300658353491	-0.0714704588285855\\
21.8068399207282	-0.0828938341209806\\
21.887593472907	-0.0936255777777433\\
21.9651688359857	-0.102673788675247\\
22.0414036101035	-0.110355793334108\\
22.1202666550655	-0.117048943576563\\
22.1973863569934	-0.122382157300308\\
22.2757420290512	-0.126606811837192\\
22.357217942502	-0.129771022285894\\
22.4390292933114	-0.131748417256375\\
22.5262014002314	-0.132612162572471\\
22.6178394274122	-0.132240000267885\\
22.7135040526488	-0.130579325934292\\
22.8130493180188	-0.127623961857438\\
22.9271436417452	-0.122908855432954\\
23.0449594468631	-0.116806528858149\\
23.1847268075255	-0.108309180204628\\
23.3372474157256	-0.0989668431724624\\
23.4422778399907	-0.0939533522896312\\
23.5540176591524	-0.0898662325809312\\
23.665475301178	-0.0869626288795473\\
23.7910632949697	-0.0849379655793285\\
23.9288508715344	-0.0840111232266594\\
24.0808681446559	-0.0842397250620053\\
24.2691601587189	-0.0857979383167446\\
24.6387962834613	-0.0904744498990482\\
24.8974235835278	-0.0930405574871713\\
25.0899155400434	-0.0937203820478913\\
25.2630335384202	-0.0930952520635202\\
25.4183772858315	-0.0913814944083384\\
25.5711668209927	-0.0885514812581363\\
25.7302126877704	-0.0843652674954427\\
25.8927653014384	-0.0787885546687761\\
26.0562432903666	-0.0719097942626661\\
26.2218914516126	-0.0637356151724333\\
26.4025721360302	-0.0535917455499515\\
26.6010182427576	-0.0412148828385739\\
26.834900908596	-0.0253800225247289\\
27.1894296522552	2.87226999802215e-05\\
27.6417834300874	0.0322551730021132\\
27.892806338864	0.0488881359283511\\
28.1132761227389	0.0622854135215931\\
28.3228033313736	0.0737464520076543\\
28.5131127428747	0.0829629247497365\\
28.7036557100763	0.0909854688733134\\
28.8889644846693	0.097602850116445\\
29.0734262021985	0.103030603118995\\
29.2597465200695	0.107363568283347\\
29.4499963288441	0.110641946419335\\
29.6465390418633	0.112879338083275\\
29.8515910302601	0.114055822167401\\
30.0691145115597	0.114121705508353\\
30.2973060954076	0.113012931591342\\
30.5481793487169	0.110560153879632\\
30.7685011772181	0.10726874205524\\
30.8506049360819	0.104745884186073\\
30.9399758204008	0.10076778214848\\
31.0347005480137	0.0952752522009206\\
31.1359533282695	0.0881231287305582\\
31.2453860866673	0.0791347076077642\\
31.3718110195532	0.0674678153381336\\
31.5323920050515	0.0513040044214534\\
32.0451116262233	-0.00137577811322842\\
32.1860749994687	-0.0140874474936084\\
32.3112600199949	-0.0242079403265336\\
32.4346057383245	-0.0329539715235683\\
32.5552579864341	-0.0402423263022058\\
32.6716766204718	-0.0460449882206433\\
32.7854564664938	-0.0505427389607718\\
32.9048929473053	-0.054040833552591\\
33.022239743137	-0.0563089550543481\\
33.1457263870845	-0.0575264478938635\\
33.280531030625	-0.0576113360603756\\
33.4245364730455	-0.0564447670728327\\
33.5878896795701	-0.0538137656242839\\
33.7701268694231	-0.0496152384184256\\
34.1439114717901	-0.0404477223647035\\
34.3099038215517	-0.0381645514423639\\
34.4935290372025	-0.0368488155654347\\
34.7103346663759	-0.0365781622026304\\
34.9878931607592	-0.0375454234458061\\
35.7324322288179	-0.0407770358952035\\
35.9868005853807	-0.0403421398845296\\
36.221577175688	-0.0387872222631671\\
36.4599440669743	-0.0359886073310633\\
36.7024942218671	-0.0318965906204625\\
36.9621974146215	-0.0262474834107849\\
37.2507114479125	-0.018703292815303\\
37.6172108612146	-0.00780601187578611\\
38.5946936560424	0.0219647196462773\\
38.9155997509455	0.0301622257271532\\
39.2143855387007	0.0365821976271548\\
39.496568583075	0.041459939711828\\
39.7866444614412	0.0452457259908385\\
40.0822233105873	0.0478604835592193\\
40.3918268598655	0.0493511059735852\\
40.7258547327524	0.0496840365889284\\
41.0776324556333	0.0488072820238443\\
41.4423356882147	0.0467145639331363\\
41.5677388787001	0.0448721021110998\\
41.7091866009263	0.0415307248115084\\
41.8653180568733	0.0365393005005359\\
42.0463661146561	0.0294597087098794\\
42.3038945333545	0.01801344020393\\
42.7267307991725	-0.000793468441742107\\
42.9345189497094	-0.00871237661320379\\
43.1181910454773	-0.0145324771201842\\
43.2991919394676	-0.0190327231066618\\
43.472074727866	-0.0221405574604319\\
43.6556626218119	-0.024197066987611\\
43.8420906917855	-0.0250743694802082\\
44.0451355849125	-0.024824900946129\\
44.2778823685668	-0.0233055739627019\\
44.5779837267004	-0.0200225767483033\\
44.8163688935869	-0.0176362597686506\\
45.0813391256829	-0.0162518437063284\\
45.394493878265	-0.0159231583586319\\
45.8711835134827	-0.0168398891508019\\
46.4694141570432	-0.0177481059332933\\
46.8501028250585	-0.0170615129008098\\
47.2040120597802	-0.015222311023706\\
47.5845098179746	-0.0119826499496725\\
48.0219093819256	-0.00696744032680385\\
48.7732639001067	0.00315235997153707\\
49.4042024737804	0.0110985746211583\\
49.8687199870182	0.015728998757055\\
50.3102655327249	0.0188876176213029\\
50.7506803437289	0.020803727988266\\
51.2225142140259	0.0216111853433887\\
51.7454262943509	0.0212569354630716\\
52.2001062352518	0.0198687705206311\\
52.4033741382177	0.0178604181321589\\
52.6499944190632	0.0141265126757872\\
53.0059015426716	0.00733699774750107\\
53.5195245129298	-0.00233657561605582\\
53.8100673283526	-0.00650764600898412\\
54.0828043249033	-0.00917774022630624\\
54.3580597126508	-0.0106129090365243\\
54.6499656238737	-0.0108987692914866\\
55.0020042607546	-0.00995308689145702\\
55.8498007726863	-0.00697718522258128\\
56.3978399380724	-0.00716728952022549\\
57.3893838065791	-0.00759766865334655\\
57.9358007114798	-0.00645487759229013\\
58.5255508025735	-0.00395452241690464\\
59.4890263766692	0.00160055711927498\\
60.361971035948	0.0061147985557426\\
61.0318899069145	0.00832869491974009\\
61.7134733570799	0.00932831680655966\\
62.4872349399364	0.00919842263641613\\
63.0079355757129	0.00813984196653905\\
63.3618953397541	0.00587999763467906\\
64.6646087798003	-0.00365142479128622\\
65.0744369706292	-0.00466578787572303\\
65.5367219912967	-0.00455709750509925\\
67.0447268678195	-0.00310594712109946\\
68.3767667172373	-0.00308929857783369\\
69.2428262657499	-0.00162997846052804\\
71.8836523712134	0.00378665239699671\\
72.9691203130259	0.00406289881222222\\
73.8095382460459	0.0032520033765735\\
74.5206053258132	0.000782155164898768\\
75.2921474087845	-0.00150021820174118\\
75.9237103753273	-0.00207066859385918\\
76.8311244692334	-0.00146384216419904\\
77.8887724917966	-0.00138536319063576\\
79.4671917731997	-0.00110339350074184\\
84.082346724609	0.00169456624563225\\
84.8439176401009	0.000917610984089379\\
86.3492673548187	-0.000863955618186196\\
87.5788993360255	-0.000618691883545353\\
95.2715892310185	0.000558085338738579\\
97.8665019578172	-0.000330796463501315\\
100.009601964004	-0.000275077257299472\\
};
\addlegendentry{Full}

\addplot [color=red, dashed, line width=2.0pt]
  table[row sep=crcr]{%
0	0\\
0.0446302474218356	-0.0509337672303474\\
0.0955841316008303	-0.107307268225725\\
0.142903776700322	-0.157822346570498\\
0.191100404098748	-0.20732641584695\\
0.237796328419378	-0.25330684328928\\
0.281962119217482	-0.294914269014996\\
0.325213646006176	-0.333817918781804\\
0.368869830114107	-0.371178137937179\\
0.413452927065876	-0.407303647301617\\
0.447424961204348	-0.43342989867439\\
0.481701663465131	-0.458545912402386\\
0.515952549677294	-0.482382746543919\\
0.54999955490004	-0.504820396292388\\
0.583849117399566	-0.525879251149945\\
0.617507675442056	-0.545581353986819\\
0.651064406523631	-0.563994230984008\\
0.684697380425035	-0.581218846559437\\
0.718464391925735	-0.597277838355041\\
0.752423235805196	-0.612188208489584\\
0.786762428903188	-0.626011629814045\\
0.821353311270101	-0.638674834314727\\
0.856168081806629	-0.650157598578147\\
0.891178939413521	-0.660446449589401\\
0.926367972523778	-0.669536836336761\\
0.961618177130688	-0.677408748721263\\
0.996925138217961	-0.684079866694674\\
1.03228444076927	-0.689570380098075\\
1.06779164269406	-0.69391358086726\\
1.10348773131504	-0.697128484310667\\
1.13942999033246	-0.699232145277605\\
1.1756757034465	-0.700237336068611\\
1.21254542574549	-0.700148481776822\\
1.24990798792177	-0.698956891999771\\
1.28774249296212	-0.696665135846388\\
1.32602804385337	-0.693281995831896\\
1.3649568441909	-0.688794750837971\\
1.40410519255916	-0.683269211571584\\
1.44343471570426	-0.676748239833685\\
1.4962084885681	-0.666564382542262\\
1.54977133677889	-0.654670456019346\\
1.60491596332731	-0.640929621743922\\
1.66280663861606	-0.625034746230384\\
1.72482943138638	-0.606527113381119\\
1.78606533569415	-0.586949090821918\\
1.86114069636166	-0.561452715342895\\
1.9164479298819	-0.542771499853103\\
1.96530123716208	-0.527772741498239\\
2.01384755236685	-0.514252390893162\\
2.05893856862861	-0.50290362063339\\
2.10660068778584	-0.492138421966445\\
2.1549327731757	-0.482468985087081\\
2.20361637143876	-0.473949123438331\\
2.25393455300787	-0.46637710267909\\
2.30740390851008	-0.459644554956597\\
2.36478692394287	-0.453845457888889\\
2.42492712443499	-0.449247192968002\\
2.48553501895042	-0.446022420529076\\
2.5472142617888	-0.444062549364006\\
2.61186299096154	-0.443299576423541\\
2.68304253905811	-0.443817381144015\\
2.75856636399381	-0.445711359082921\\
2.83428748215401	-0.448762690632577\\
2.92885586728592	-0.453817404840152\\
3.04757284862571	-0.461435419088033\\
3.3698466060347	-0.482894234536147\\
3.47142346098568	-0.488210758630672\\
3.56871625033851	-0.49202997835016\\
3.64426545741838	-0.493945716593487\\
3.71849362029097	-0.494801708100354\\
3.79195574038961	-0.494542923767611\\
3.86505297801587	-0.493103552111108\\
3.93805705661393	-0.490413330755018\\
4.01112707588611	-0.486404935150446\\
4.0843506887182	-0.481019574980991\\
4.15778411814975	-0.474209148816755\\
4.23148526182774	-0.465934867099847\\
4.305538153848	-0.456163027056249\\
4.38006404021149	-0.444859802719961\\
4.45521791634525	-0.431987277162591\\
4.53117225791136	-0.417502620395098\\
4.60809038557358	-0.401361810047504\\
4.68609165371554	-0.383528818000144\\
4.76521207577933	-0.363990288011919\\
4.84536913671406	-0.3427735225536\\
4.92634691161057	-0.319961880570176\\
5.00781933197831	-0.29569874882327\\
5.08941569898268	-0.270173685973631\\
5.17080707124526	-0.243593574709251\\
5.27195867777402	-0.209168901152879\\
5.37237720591025	-0.173677798033111\\
5.49218362678839	-0.129946818411355\\
5.63246887730223	-0.0773470501285374\\
5.86243618152717	0.0104791653917005\\
6.08803775699765	0.0963019072918598\\
6.23475803925867	0.150843048866875\\
6.34412140879071	0.190443165227293\\
6.45411550514569	0.229123707556596\\
6.56087765426503	0.265376599083936\\
6.66442065550002	0.2991593409632\\
6.7650430727594	0.330555818819093\\
6.86321591747331	0.359717405506274\\
6.95950223908451	0.386821675721677\\
7.05449246632934	0.412042660254841\\
7.14876831616444	0.435534809737959\\
7.24288354346282	0.457425016915266\\
7.33734240139238	0.477806694560712\\
7.43255167371412	0.496731456380488\\
7.50461380986746	0.509974455730315\\
7.57710582461127	0.522359438018611\\
7.64996087854217	0.533863693863779\\
7.72275285311235	0.544419696191241\\
7.81890134098657	0.556937670902883\\
7.91325629856742	0.567665371983992\\
8.00573939728113	0.576712436682897\\
8.09730516911007	0.584275222018206\\
8.18999864776684	0.590564010347165\\
8.28728005070765	0.595741056508331\\
8.36653526616422	0.598921259970794\\
8.44837353548169	0.601259316255593\\
8.55291368245732	0.602891485986476\\
8.64225398210023	0.60312596052789\\
8.74898314699388	0.602077136511795\\
8.86277499482226	0.599461426112939\\
8.96993365459215	0.595655626120916\\
9.07859189146671	0.590527748336442\\
9.19717603121302	0.583579023837189\\
9.30817202130807	0.575865732174293\\
9.39208189245646	0.56898576541829\\
9.42904365106931	0.564488830160528\\
9.46679943505735	0.558622228537274\\
9.50124840501783	0.552179332335442\\
9.5383193131798	0.544122236894566\\
9.57818525854138	0.534205513392607\\
9.62015913879151	0.522422116280325\\
9.66402921576289	0.508707454435509\\
9.70977906251262	0.492968737641789\\
9.75807474902193	0.474868103641725\\
9.8089638430302	0.454270957946392\\
9.86251942404526	0.431056470247029\\
9.91853349040633	0.405266195935212\\
9.97802516695756	0.376397685507854\\
10.0393633128856	0.345284284781357\\
10.1092348968962	0.308496410785679\\
10.2009036516557	0.258655851995499\\
10.351619255103	0.174841086849085\\
10.4818917074021	0.102772391483654\\
10.5720119234454	0.0542364781997691\\
10.6478149111526	0.0147431388646595\\
10.7139216395611	-0.0184456889074909\\
10.7742329860252	-0.0475412707601492\\
10.8342752424407	-0.0752592479503704\\
10.8938407318637	-0.10142284139539\\
10.9500857782172	-0.124826035055094\\
11.0041332028476	-0.146060203022117\\
11.0589311623044	-0.166282789148383\\
11.1151731213353	-0.185628433004268\\
11.1722421002876	-0.20376611961602\\
11.2291429010691	-0.220330846617429\\
11.2853627236041	-0.235193707881322\\
11.3413343904832	-0.248501538077008\\
11.3975427662865	-0.260374208482475\\
11.4543631475168	-0.270870619331788\\
11.5120641241249	-0.280004059116507\\
11.5703711760048	-0.287704566475\\
11.6290322218858	-0.293941930221365\\
11.687874351641	-0.298724895449155\\
11.7469914443864	-0.302098129717024\\
11.8066464985648	-0.304108137004988\\
11.867228274362	-0.304785779748713\\
11.9291061940623	-0.304138582524629\\
11.9926264080915	-0.302156250126941\\
12.0574167284219	-0.298857916903827\\
12.1228804387418	-0.294321117883015\\
12.1888797643436	-0.28863051802999\\
12.2704333256364	-0.280214813000569\\
12.3577627662906	-0.269731575292226\\
12.4500506891163	-0.25727601730668\\
12.5625620054729	-0.240640895544644\\
12.6299240354254	-0.231265725146983\\
12.7025862612799	-0.222454215980889\\
12.777194379806	-0.214764973202591\\
12.8545973516301	-0.208169636763628\\
12.9344290945555	-0.202750627166793\\
13.0108401121523	-0.19879162194438\\
13.093931765544	-0.195738710414361\\
13.1788468473984	-0.193827437494718\\
13.2672224606871	-0.19297466127405\\
13.3650018316126	-0.193187292533267\\
13.4857077147609	-0.194806102970958\\
13.6276202461512	-0.198076485566503\\
13.8653246428047	-0.20506082105112\\
14.0944447147703	-0.211307098366078\\
14.2328801475447	-0.213950197305877\\
14.3646566116714	-0.215240679646058\\
14.4749092479001	-0.215193302369585\\
14.5844865446979	-0.213984436244303\\
14.6940745338083	-0.211512098060794\\
14.8039804725061	-0.207689508299211\\
14.9144935536213	-0.202446750237073\\
15.0257273554722	-0.19574068428193\\
15.1383661448892	-0.187504269586356\\
15.252634568704	-0.177702273653537\\
15.3694935111248	-0.166233600567267\\
15.4885868622828	-0.153125535748686\\
15.6099727735279	-0.138393147531417\\
15.7322874771613	-0.12227021848372\\
15.8544230253597	-0.10503283527639\\
15.9957445720628	-0.0838788944585644\\
16.1557316708782	-0.0586964554040748\\
16.376868411904	-0.0224887666410467\\
16.9024438166317	0.0641801412248242\\
17.0934329355191	0.0940296166486689\\
17.2544398470738	0.11790393098336\\
17.4075530325189	0.13926474027933\\
17.5554188727403	0.158471760190196\\
17.6751042741394	0.172901229364115\\
17.7932752467945	0.186111683921382\\
17.9109238637378	0.198207073691165\\
18.028959490765	0.209263843552037\\
18.1483708225169	0.219343212336767\\
18.2691631122362	0.228408729872839\\
18.3905452317469	0.236381094295893\\
18.5108263861796	0.243169274044533\\
18.6509612488824	0.249708501875858\\
18.7884317743642	0.254736294983658\\
18.9302234310003	0.258551566466167\\
19.0602759862308	0.260894478673364\\
19.1976504324158	0.262229456851784\\
19.3531732602433	0.262401829982551\\
19.5083943933864	0.26126335273149\\
19.6594827370244	0.258997280691162\\
19.819689682277	0.255433326911685\\
19.9933221927835	0.250346689040555\\
20.0890384662928	0.246626805940508\\
20.1433284034471	0.24313818994186\\
20.1998163961305	0.23831674960141\\
20.254750284582	0.232521572783298\\
20.3139477098374	0.22512837886687\\
20.3769089507444	0.216049883898407\\
20.4444516833427	0.205038124280918\\
20.5172132570857	0.191858116189792\\
20.5951008214455	0.176434245497546\\
20.6793529607564	0.158474520335318\\
20.7745625238794	0.136939066506827\\
20.9020836769215	0.106725102266338\\
21.2588301231279	0.0213142162786681\\
21.3688614823364	-0.00329810918134399\\
21.4592022406777	-0.0223500253718356\\
21.5434004320001	-0.0389932848366499\\
21.6241279227099	-0.0538260346657751\\
21.7106510621409	-0.0683995234542323\\
21.7889324115093	-0.0803365377399956\\
21.8687979138054	-0.0912483910984463\\
21.947956712786	-0.100775801026984\\
22.0263254092472	-0.108934723644381\\
22.105333367798	-0.115882072491019\\
22.1857834059109	-0.121657983957547\\
22.2674676260168	-0.126218656415702\\
22.3497459560298	-0.129534451543506\\
22.4324839557429	-0.131636747428061\\
22.5164270102567	-0.132580391831027\\
22.6026429812106	-0.132390861067236\\
22.6916917676336	-0.131067758324576\\
22.7961371100243	-0.128210630748072\\
22.9022930814913	-0.124048580491404\\
23.0148908439305	-0.118470999930665\\
23.1534950894438	-0.110308895241687\\
23.3614058001489	-0.097710845275131\\
23.4726704135174	-0.0927213105273523\\
23.579972820105	-0.0890923202462517\\
23.6984239017026	-0.0863121369667113\\
23.8304153800699	-0.0845499242552421\\
23.9696767274254	-0.0839594000422323\\
24.1308441056637	-0.084540613648727\\
24.3361218577525	-0.0865790292291564\\
24.9269749664698	-0.0932171683836032\\
25.1133150011215	-0.0937146530762476\\
25.2777826468724	-0.0930008431153055\\
25.4423487109477	-0.0910374255637549\\
25.5895115542849	-0.0881482420967643\\
25.738138331616	-0.0841231704185219\\
25.8893289325162	-0.0789067325458319\\
26.0443659794796	-0.0724354024558238\\
26.2239885197803	-0.0636208531242204\\
26.4071184407954	-0.0533332511786995\\
26.6102277691997	-0.0406310382383452\\
26.850570005066	-0.0242857666615919\\
27.2223717430043	0.0024342338353307\\
27.6241087075828	0.0310331740479057\\
27.8695147316479	0.0473931155280383\\
28.1009665022495	0.0615748022998446\\
28.2961216619482	0.0723722739861472\\
28.4856612724857	0.081716681881133\\
28.6741247804381	0.0898257119719261\\
28.8651490348159	0.0968164145880053\\
29.0590984765426	0.102651218260505\\
29.250246470263	0.107176005604785\\
29.4576169848088	0.110754402610567\\
29.6486189789184	0.11289839298739\\
29.8622381665572	0.114090595518093\\
30.0847362482707	0.114085116538547\\
30.323560686946	0.112817784361326\\
30.5702824613294	0.110290486319499\\
30.7753576054951	0.107103968108547\\
30.8596984267603	0.104400655119775\\
30.9457482683032	0.100471435132818\\
31.0421044401353	0.09479673196428\\
31.1384893646051	0.0879301873589498\\
31.2460497146257	0.079078088187174\\
31.3771199336121	0.0669544113495277\\
31.5391821803452	0.0505992176475729\\
32.0356444753794	-0.000479447383767706\\
32.1810084266344	-0.0136543781909495\\
32.3087271599167	-0.0240165855364154\\
32.428723807272	-0.0325676629800995\\
32.5424074493394	-0.0395285914863308\\
32.6552220749422	-0.0453003864914194\\
32.7675891005616	-0.0499145059171582\\
32.882547109764	-0.0534809631165274\\
32.9997083412786	-0.0559614518851248\\
33.1296993783089	-0.0574344804692259\\
33.2634453987291	-0.0576687644272909\\
33.4038677401857	-0.0566853965332257\\
33.5620144349074	-0.0543124014942151\\
33.7342274173319	-0.0505289289316266\\
34.1844515097978	-0.0397852203636688\\
34.3590476583624	-0.037698056306823\\
34.5513780890823	-0.0366591372509504\\
34.7719270831098	-0.0366938426578827\\
35.0935574904227	-0.0381210457339307\\
35.6345361042679	-0.04064023925649\\
35.9111955917057	-0.0406065922799428\\
36.1489048812231	-0.039403984763652\\
36.3891376592877	-0.0369499858396125\\
36.6357836556953	-0.0331355613169961\\
36.893148451884	-0.0278643142883226\\
37.1781241292474	-0.0207144767242795\\
37.5198231316437	-0.0108016634356289\\
38.7243070726218	0.0254242238766551\\
39.0411185188035	0.0330165442841661\\
39.324014660445	0.0386190378118698\\
39.6111556098881	0.0431029635199138\\
39.9002840657055	0.0463980828347701\\
40.2006447277078	0.0485760355040838\\
40.5100060834523	0.0496126954023453\\
40.8496496994747	0.0495091435993373\\
41.2301667800764	0.0480899018218395\\
41.4929457027007	0.0461080087212054\\
41.6283507029	0.0435940859608337\\
41.7726347072047	0.039651000413798\\
41.9305762424357	0.0341225111528871\\
42.1323557208498	0.0257580793891776\\
42.8974955982961	-0.00739819371135297\\
43.0868938729487	-0.013628240152542\\
43.2682638509885	-0.0183546100041667\\
43.4480793520547	-0.0217793861314277\\
43.6333115299335	-0.0240131406443282\\
43.8223691569981	-0.025035830376595\\
44.0201881211321	-0.0249160036790528\\
44.2564539874745	-0.0234911792611996\\
44.5501518184556	-0.0203703882365858\\
44.8181466283446	-0.0176231698493154\\
45.0774511282367	-0.0162641220613864\\
45.3820229010563	-0.0159160189496532\\
45.8371822481861	-0.0167560994186147\\
46.4696994737925	-0.0177490688199242\\
46.8537747392827	-0.0170525103794859\\
47.2076905698156	-0.0151948610374575\\
47.5806429665929	-0.0120205955559243\\
48.0071498185502	-0.00715632440179093\\
48.7380435899167	0.00267615834903268\\
49.3835792259741	0.0108649360107194\\
49.8436055585908	0.0155138427402619\\
50.271567062401	0.0186608852715437\\
50.7238079215374	0.0207218363199928\\
51.1954701506165	0.0215969984069488\\
51.7173091673527	0.0213051023460764\\
52.1943116673857	0.0199090489802956\\
52.3983493313273	0.0179233558270511\\
52.6372826224454	0.0143455198416831\\
52.9808808322524	0.00783694001958679\\
53.5241986702861	-0.00241387644493329\\
53.8173055142455	-0.00659518692316396\\
54.0884665637067	-0.00922031440656212\\
54.3662930263143	-0.0106376869542828\\
54.6559530713508	-0.0108933538265319\\
54.9992624243772	-0.00996467786586663\\
55.8536114614698	-0.0069745811454851\\
56.4076380835752	-0.00717822947568436\\
57.3887090286118	-0.0076002361043237\\
57.9274471730787	-0.00648030343408834\\
58.5273228098573	-0.00394663289847585\\
59.4982160752715	0.00165485374736818\\
60.3686578026156	0.00614404700723981\\
61.0352482440635	0.00833668892232708\\
61.7250241432829	0.00933548660836436\\
62.5022243380487	0.00918549790802103\\
63.0040962025862	0.00815773859773117\\
63.3616673327063	0.00588162928875136\\
64.6556894220456	-0.0036160772960443\\
65.0692490281577	-0.00466034182433361\\
65.5302518984542	-0.00456573051120301\\
67.0702891677214	-0.00311756532894947\\
68.3632396807833	-0.00310262689627905\\
69.2282973924737	-0.00166343888595577\\
71.9051748645092	0.0038043038900355\\
72.9879986810895	0.00405844312486181\\
73.81372848674	0.00324100425778795\\
74.5359313106466	0.000724909417030517\\
75.2980282881084	-0.00151134247569473\\
75.9245435102537	-0.00207074338246116\\
76.8368319408805	-0.00146016229979296\\
77.8939064807672	-0.00138652941163286\\
79.4575346821229	-0.00111018772223304\\
81.6499618270031	0.00110468706340328\\
83.186802647934	0.00177784456714392\\
84.605422747483	0.00126817253318734\\
86.6690534497882	-0.000899169001016276\\
91.3041411339565	2.66640295905063e-05\\
94.0547188706705	0.000776120611391207\\
95.8613853103512	0.000157282687837323\\
97.485081076356	-0.000383558689975416\\
100.002997367862	-0.000275281092982027\\
};
\addlegendentry{Graph-style ROM}

\end{axis}

\begin{axis}[%
width=0.7\textwidth,
height=0.7\textwidth,
at={(0in,0in)},
scale only axis,
xmin=0,
xmax=1,
ymin=0,
ymax=1,
axis line style={draw=none},
ticks=none,
axis x line*=bottom,
axis y line*=left
]
\end{axis}
\end{tikzpicture}%

%% file: Highnonlin_timehist_x3_new.tikz
% This file was created by matlab2tikz.
%
%The latest updates can be retrieved from
%  http://www.mathworks.com/matlabcentral/fileexchange/22022-matlab2tikz-matlab2tikz
%where you can also make suggestions and rate matlab2tikz.
%
\definecolor{mycolor1}{rgb}{0.12941,0.12941,0.12941}%
\begin{tikzpicture}

\begin{axis}[%
width=0.7\textwidth,
height=0.7\textwidth,
at={(0in,0in)},
scale only axis,
xmin=0,
xmax=10.473272861437,
xlabel style={font=\color{mycolor1}},
xlabel={$t [s]$},
ymin=-2,
ymax=0.6,
ylabel style={font=\color{mycolor1}},
ylabel={$x_3$},
axis background/.style={fill=white},
xmajorgrids,
ymajorgrids,
legend style={legend cell align=left, align=left, font=\tiny}
]
\addplot [color=black, line width=2.5pt]
  table[row sep=crcr]{%
0	0.0134347315302978\\
0.101784258059871	0.0174288167088061\\
0.2113112154903	0.0261186768075063\\
0.344980810968421	0.0412301076674257\\
0.681154242015284	0.0813104414090873\\
0.788641065132149	0.0886200213772668\\
0.887465572088505	0.0911393565192213\\
0.982913291064493	0.0893617616728335\\
1.07985297940446	0.083345298306984\\
1.18389752003872	0.072637765864723\\
1.30757961982297	0.0554769578191383\\
1.53469435905829	0.018504926495158\\
1.69785113561043	-0.00568961139988211\\
1.85761104513677	-0.0246426817296008\\
2.26262767116553	-0.0704730438889722\\
2.38254286098408	-0.0897062165792128\\
2.48591336881523	-0.110561100515534\\
2.5816260163096	-0.134272602545032\\
2.66962485576837	-0.160417567565281\\
2.75268637135405	-0.189355518388396\\
2.83542901150573	-0.22264042332484\\
2.91559759216094	-0.259389585801758\\
2.99471091848866	-0.300174528837221\\
3.07721879491716	-0.347605409110519\\
3.15973231314804	-0.400073228205477\\
3.2448142456807	-0.459388291062364\\
3.33241180231861	-0.525818501106855\\
3.42202514538967	-0.599102928206289\\
3.52058937692104	-0.685402076204303\\
3.62831853430909	-0.785672096330758\\
3.75026552046934	-0.905155461773678\\
3.91718030768905	-1.07521591976288\\
4.19843323868226	-1.36215838776723\\
4.31646708546171	-1.47597944260661\\
4.41249991321974	-1.56307286727256\\
4.49859475531624	-1.63576926267777\\
4.58039143720523	-1.69925826388413\\
4.65416525244631	-1.75124147129163\\
4.72548943832589	-1.79625870824296\\
4.79075207847153	-1.83258991497246\\
4.85367309125663	-1.8629552210439\\
4.91570899014699	-1.88820111724533\\
4.9768786004842	-1.90836239910874\\
5.03574331699883	-1.92320030509232\\
5.09470591510422	-1.93348140302312\\
5.15137531594388	-1.9389812313878\\
5.20728241033235	-1.94016696253124\\
5.2640091684329	-1.93705870589588\\
5.32038170832249	-1.92968552494465\\
5.37880854892663	-1.91758189443196\\
5.43704514975092	-1.90107156630298\\
5.49735255265785	-1.87940339046973\\
5.56016372062816	-1.85204240061733\\
5.62163177983024	-1.82071424816635\\
5.68895383896352	-1.78148671093679\\
5.75949354446141	-1.73520380132866\\
5.83316298888992	-1.68163833082632\\
5.91257243379082	-1.61848725511743\\
5.99779796111237	-1.54522263025596\\
6.09446399042804	-1.45635452788249\\
6.2080644691799	-1.34586499476262\\
6.36126686291756	-1.19035094338751\\
6.68333092577531	-0.861905929475004\\
6.80299939896308	-0.746624271544029\\
6.90787028424949	-0.651130537998581\\
7.00707086897079	-0.566558459301454\\
7.09975226141746	-0.493239690156091\\
7.18824520479153	-0.428770271113718\\
7.27206439535189	-0.372918363030292\\
7.352759315527	-0.32404060094006\\
7.43393733108152	-0.279722002932028\\
7.5140068078569	-0.24070462257864\\
7.59300588738966	-0.206638949015947\\
7.67851461353291	-0.174488234458549\\
7.76425426362703	-0.146850977716886\\
7.85396985057914	-0.12242089495386\\
7.95070628942989	-0.100613005579691\\
8.052455080779	-0.0819520618708385\\
8.17539715212243	-0.0639452434104602\\
8.33288616114397	-0.0455576394732624\\
8.7126431143018	-0.00295347797504952\\
8.85990204683543	0.0192353416234425\\
9.19082391640304	0.0707349348410649\\
9.30210431802885	0.0826587712951561\\
9.40249128506676	0.0891672239185031\\
9.49889590860885	0.0911486490594822\\
9.59705879826074	0.0888217064018821\\
9.70141090655326	0.0819580591881728\\
9.82264465127991	0.0695215029752063\\
10.0474361002874	0.0409377036002141\\
10.1958798515026	0.0244995108318733\\
10.306030230322	0.0164297356928351\\
10.4075589988049	0.0131865583368338\\
10.473272861437	0.0134347318154635\\
};
\addlegendentry{Full System (6D)}

\addplot [color=red, dashed, line width=2.5pt]
  table[row sep=crcr]{%
0	0.0134347315302978\\
0.10157373425335	0.0174190033843757\\
0.209930036101056	0.0259862712425214\\
0.343900042327759	0.0410939111380095\\
0.681757224382991	0.0813664932997415\\
0.783064979859958	0.0883629354602942\\
0.887952655498262	0.0911519148518227\\
0.98279610478334	0.0893578395241761\\
1.08293074310509	0.0830865222454769\\
1.18565533265894	0.0724274735237511\\
1.31268938104854	0.054693565422177\\
1.55945126431606	0.0145435292879732\\
1.70756653494276	-0.00696310045971416\\
1.87458846053213	-0.0264708902281114\\
2.24715115307739	-0.0683043137726962\\
2.36588763220469	-0.0867481789739983\\
2.47176060352449	-0.107434883049601\\
2.56203313220449	-0.129040973511565\\
2.6545572520202	-0.15562619063201\\
2.73670377737647	-0.183449157969706\\
2.82299717156272	-0.217343267973858\\
2.90291325743969	-0.253277364339198\\
2.98176494765446	-0.293189641309782\\
3.06435050343628	-0.339878003074249\\
3.14630673943197	-0.39119675387791\\
3.23040713481459	-0.448979978287177\\
3.31884249626872	-0.515183936057904\\
3.40874521783801	-0.587920142477865\\
3.50886990069918	-0.674853453232531\\
3.61308427348864	-0.77116327397367\\
3.73023965432354	-0.885176899909919\\
3.89146761314767	-1.0487365700291\\
4.21174201050049	-1.37528606561322\\
4.3281421889042	-1.48685959460114\\
4.42342585923563	-1.57258822822805\\
4.5067951564667	-1.64237747923038\\
4.58396508411644	-1.70188020993263\\
4.66463249546148	-1.75817381537606\\
4.73302916012482	-1.80068697534929\\
4.79854850456723	-1.83659224111494\\
4.85743854491984	-1.86460201573475\\
4.91947036509446	-1.88955739685258\\
4.9678136115469	-1.90563257305698\\
5.03058312768231	-1.92202450074712\\
5.08309299386082	-1.93176447605559\\
5.14012126615974	-1.9381529645333\\
5.21347285139425	-1.93994418515732\\
5.28211109189715	-1.93524874236933\\
5.34707639466788	-1.92470566825969\\
5.40810765598114	-1.90980177394506\\
5.46879631037782	-1.89021350267567\\
5.53232242437541	-1.86473374796773\\
5.59282853700294	-1.83591449690509\\
5.65330887356054	-1.8028572665392\\
5.72492965054723	-1.75849697370936\\
5.79918809141929	-1.70694923825161\\
5.87767273472601	-1.64686391237463\\
5.95796790686701	-1.58009045440061\\
6.0536967264686	-1.49446538755534\\
6.16397079601925	-1.38932583880241\\
6.2872393117764	-1.26610049499637\\
6.7540074694274	-0.79306491487835\\
6.86757772513069	-0.68710293157419\\
6.96942584446708	-0.597916570222571\\
7.0641063488643	-0.520732499187911\\
7.15279695084794	-0.453906428203215\\
7.23743477966241	-0.395347977492397\\
7.3170968629866	-0.345029528908741\\
7.40055554459355	-0.297346667692059\\
7.48446677649821	-0.254553686277598\\
7.56767182456147	-0.217080745985344\\
7.65158822890475	-0.184090477411235\\
7.73673333233853	-0.155236460207234\\
7.82562853404856	-0.129666666467287\\
7.92293777477424	-0.106427063133209\\
8.02583137045078	-0.0864566025489335\\
8.14184264181505	-0.0684502005084475\\
8.28419327372082	-0.0508669025753719\\
8.76992613968088	0.00521084905334135\\
8.94865699251435	0.0338672223856467\\
9.15398151792009	0.0658872500871333\\
9.26399593423156	0.0790870653679843\\
9.36807183660915	0.0874286023519204\\
9.4598610996064	0.0908173787320266\\
9.56168706840561	0.0901723193118915\\
9.66594909766913	0.0847522101982516\\
9.78452751289892	0.0738377525952671\\
9.94292611809418	0.0544136634591013\\
10.1661422635043	0.0273633290450181\\
10.2836263128163	0.0177106087099119\\
10.3861917192994	0.0134953252665344\\
10.4736065884256	0.0134231282184558\\
};
\addlegendentry{Arc-length ROM (2D)}

\end{axis}

\begin{axis}[%
width=0.7\textwidth,
height=0.7\textwidth,
at={(0in,0in)},
scale only axis,
xmin=0,
xmax=1,
ymin=0,
ymax=1,
axis line style={draw=none},
ticks=none,
axis x line*=bottom,
axis y line*=left
]
\end{axis}
\end{tikzpicture}%

%% file: discont_x6_timehist.tikz
% This file was created by matlab2tikz.
%
%The latest updates can be retrieved from
%  http://www.mathworks.com/matlabcentral/fileexchange/22022-matlab2tikz-matlab2tikz
%where you can also make suggestions and rate matlab2tikz.
%
\definecolor{mycolor1}{rgb}{0.12941,0.12941,0.12941}%
\begin{tikzpicture}

\begin{axis}[%
width=0.7\textwidth,
height=0.7\textwidth,
at={(0in,0in)},
scale only axis,
xmin=0,
xmax=50,
xlabel style={font=\color{mycolor1}},
xlabel={$t [s]$},
ymin=-0.1,
ymax=0.1,
ylabel style={font=\color{mycolor1}},
ylabel={$\dot{x}_3=y_3$},
axis background/.style={fill=white},
xmajorgrids,
ymajorgrids,
legend style={legend cell align=left, align=left, font=\tiny}
]
\addplot [color=black, line width=2.0pt]
  table[row sep=crcr]{%
0	0.00928155248219298\\
0.0950651822927213	0.0109400421844725\\
0.180480230247561	0.0121134335640534\\
0.259749262967851	0.0128739606843737\\
0.332490794888692	0.0132563948907887\\
0.403859251465519	0.0133136004390764\\
0.474280444949429	0.0130457056945801\\
0.549283950590819	0.0123968454277659\\
0.622610847902244	0.0113986288546997\\
0.696672198739854	0.0100317596099373\\
0.772023986781228	0.00828574417255368\\
0.848293481280827	0.00617627582540337\\
0.929641343366356	0.00358108795084178\\
1.01836669104085	0.000395865787304217\\
1.11632085693062	-0.00347164985317505\\
1.23721836309868	-0.00860685457256949\\
1.64523709503006	-0.0262634378347499\\
1.74587409365765	-0.030085142159848\\
1.8320155798129	-0.0330280697512393\\
1.91774333548015	-0.0356022116440329\\
2.12288460450836	-0.041413031937644\\
2.20707858187951	-0.0443980994731419\\
2.29729732261389	-0.0479652759179743\\
2.39684788947579	-0.0522724358332667\\
2.50559740247516	-0.0573183958048631\\
2.66715221778206	-0.065198809418618\\
2.89476189079388	-0.0762840903705566\\
3.0060426747828	-0.0813486223272619\\
3.09811175395154	-0.0852174052553423\\
3.18084419509606	-0.0883755369783259\\
3.26388666300765	-0.0911864173711123\\
3.33845745583752	-0.0933628947795668\\
3.4138088342	-0.0951933719273796\\
3.48322240771924	-0.0965257873086998\\
3.55072943745972	-0.0974761493269654\\
3.6187486843596	-0.0980728684653513\\
3.68562255371455	-0.0982935226231021\\
3.74825986326944	-0.0981622629744479\\
3.80784634011955	-0.0977286128496999\\
3.87004791431805	-0.0969514216751648\\
3.93370506315414	-0.0958115643826645\\
3.99461914085374	-0.0943954650649985\\
4.0599720625351	-0.0925257088257254\\
4.1230638079557	-0.0903820130313946\\
4.18461529752655	-0.0879773794933953\\
4.25189525300723	-0.0850056623571191\\
4.31866241305083	-0.0817156250701316\\
4.38538920573455	-0.0781042250261095\\
4.45917671779308	-0.0737561112337417\\
4.53380265599889	-0.0690066315243385\\
4.60959390930881	-0.0638520168654466\\
4.6939731395539	-0.0577627047692317\\
4.78051891307924	-0.0511840502945802\\
4.8772195148285	-0.0435034572030446\\
4.99355169880984	-0.0339129277365089\\
5.1508464297549	-0.0205769005333707\\
5.44003927900115	0.00398176530080718\\
5.5630717740948	0.0140523707775912\\
5.66660389233325	0.0222060855414483\\
5.76405280131652	0.0295472290406309\\
5.8529794915391	0.035920650348622\\
5.93700500263783	0.0416264161665296\\
6.01987282626501	0.0469288881523369\\
6.10063137079564	0.0517679830717057\\
6.18003343031624	0.0561958337612509\\
6.26045980061801	0.060336499291104\\
6.33818803106653	0.0640016021327838\\
6.41534820808814	0.0673080508869575\\
6.49009169424806	0.0701935904956343\\
6.56920755076732	0.0729076362566445\\
6.64441033855562	0.075164856424955\\
6.72378930362392	0.0772102097776397\\
6.79857216813777	0.0788254594441682\\
6.87835816289428	0.0802225529908824\\
6.96201084469604	0.0813354413341827\\
7.0465735486093	0.0821064537483736\\
7.13207656128477	0.0825384350464873\\
7.22265382506575	0.0826324397193119\\
7.31269680578468	0.0823749655688033\\
7.4024816953246	0.0817912533172631\\
7.49641034506653	0.080855683953672\\
7.59678927061029	0.0795189311655946\\
7.70434494579761	0.0777390402458025\\
7.8161164898827	0.0755555946813047\\
7.93989059697918	0.0728038322478852\\
8.08200495573207	0.0693060685955871\\
8.27194564699779	0.0642666516255588\\
8.66256350392378	0.0537999102389364\\
8.81151237942331	0.0502170620051814\\
8.85000660322876	0.0492476795996311\\
8.98848229077763	0.0391913344398489\\
9.12042392068908	0.029743392169955\\
9.20619420825987	0.0239491597511901\\
9.27892933846184	0.0193673189798886\\
9.34572906304602	0.0154879447802401\\
9.40636583034518	0.0122750598618993\\
9.46748007608863	0.00935981561646315\\
9.52909513988818	0.00676773171502276\\
9.58940046123928	0.0045802414087035\\
9.64755065718568	0.00280391375647326\\
9.70689013952573	0.00132853071851002\\
9.7678933124416	0.000161863445107713\\
9.82823053689831	-0.000653102377860648\\
9.8913117823653	-0.00116087587400671\\
9.95612695424673	-0.00133884581280341\\
10.0247870255544	-0.00117909600166399\\
10.097886525384	-0.000658746537993693\\
10.1764837579173	0.000244223245509545\\
10.2673898828403	0.00164081749906586\\
10.3793207505528	0.00371441240447723\\
10.7211439149637	0.0102895659584519\\
10.8141430714852	0.0116221639288696\\
10.900737666887	0.0125300332374678\\
10.9796773462175	0.0130297494691831\\
11.0564222801438	0.0131841947375406\\
11.1303669313857	0.0130044604443569\\
11.2031820440242	0.0125017719658089\\
11.278545815148	0.0116378298999322\\
11.356043969751	0.0103894057054674\\
11.4345347473234	0.00876601253354181\\
11.513905141664	0.00677891187878288\\
11.5988693734941	0.00430064027494836\\
11.6894897901975	0.001307386141157\\
11.7922435909563	-0.00244376918743683\\
11.915302533272	-0.0072929059936655\\
12.3380283222621	-0.0242737894653047\\
12.439727813862	-0.0278397273788613\\
12.5348906655005	-0.0308181826919167\\
12.6182538556417	-0.033092688105981\\
12.8195494902941	-0.0382831881033283\\
12.9036970844576	-0.0410027118026264\\
12.9940040169072	-0.044264484532718\\
13.0930777100711	-0.0481870694963575\\
13.2088288452374	-0.053118991511468\\
13.3795592676585	-0.0607862015897709\\
13.5963207773677	-0.0704882225066328\\
13.7066968473193	-0.0751073204529433\\
13.8069163800822	-0.0789631802176203\\
13.8977111086425	-0.082096517294616\\
13.9796139971979	-0.0845731907780447\\
14.0530265935079	-0.0864737088355909\\
14.1276472899795	-0.0880655412607609\\
14.1983459849396	-0.0892336260560285\\
14.2661274640835	-0.0900243792278133\\
14.3345680411122	-0.0904809815219565\\
14.4026904297258	-0.0905825101507958\\
14.4661394908303	-0.0903524604171224\\
14.5337800267286	-0.0897567813180729\\
14.5965115983362	-0.0888784419654414\\
14.6590924949326	-0.087689401788289\\
14.7253534661336	-0.086091968664654\\
14.7886910634529	-0.0842437991810385\\
14.8567838335666	-0.0819138860285875\\
14.9232914264026	-0.0793047823403299\\
14.9956615619039	-0.0761054315650895\\
15.0678457450827	-0.072558022567506\\
15.1403983123853	-0.0686550576523928\\
15.2137365486404	-0.0643901784253629\\
15.2950427023609	-0.0593187678197395\\
15.378090109314	-0.053805980460865\\
15.4704392515105	-0.0473349624902824\\
15.5732343007478	-0.0397860582861682\\
15.6959245916239	-0.0304226276120403\\
15.8771637457294	-0.0162048255889857\\
16.1177830351562	0.00262265204987244\\
16.2440885128232	0.0121600948449867\\
16.3555172265703	0.0202341560912487\\
16.4559538922175	0.027169108993597\\
16.5507253863666	0.033367717630874\\
16.6426919097562	0.039027576596105\\
16.7273263914545	0.0439031704421211\\
16.8110556937746	0.0483951278356685\\
16.894768443344	0.0525436358470301\\
16.9738173709265	0.0561381188408916\\
17.0573822510036	0.0595911316405733\\
17.1415924692124	0.0627072983104071\\
17.2242835788244	0.0654120980422874\\
17.3039343843996	0.0676870568502181\\
17.3882159754258	0.0697461063632616\\
17.4684788285391	0.0713807196477987\\
17.5541121434112	0.072782907509243\\
17.6355754058676	0.0737995776683888\\
17.7197889687676	0.0745375539957251\\
17.8047912007633	0.0749735213788085\\
17.8986951181431	0.0751125953977052\\
17.9924061233812	0.0749140019538785\\
18.0925350278002	0.0743552759304933\\
18.1946598452794	0.0734462854172051\\
18.3051907940027	0.072113608816295\\
18.422638407711	0.0703462669678672\\
18.5515413893421	0.0680511469695375\\
18.6944718348033	0.0651560475435531\\
18.8667277080444	0.0613177664765416\\
19.4778033888558	0.0473656753317258\\
19.5664147768146	0.0454453802240593\\
19.7001863850755	0.0366424505869176\\
19.8461771613929	0.0271346374640444\\
19.937442311516	0.0215481258582884\\
20.0109687949658	0.0173666130023165\\
20.0784361267634	0.0138369720812506\\
20.1465004079043	0.0106137918521441\\
20.2082458971253	0.00800838636652657\\
20.270495777912	0.00570536044751435\\
20.3357609701557	0.00365076829613287\\
20.3965816184844	0.00207189921483319\\
20.4611269731347	0.000748823314481228\\
20.5237398938176	-0.000195325499070975\\
20.5884750543446	-0.000833657566751356\\
20.6562881712397	-0.00115480563093229\\
20.7249833603976	-0.00114556848473057\\
20.7994999133056	-0.000794295310988957\\
20.8796916984811	-7.46008471210757e-05\\
20.9708711113077	0.00109036708857957\\
21.0816095854991	0.00285518807685037\\
21.4823480027097	0.00957169667589852\\
21.5762428708604	0.010627680816178\\
21.6632328019266	0.011275827503141\\
21.7468726693069	0.0115559701423251\\
21.8260768224998	0.0114843126706248\\
21.9054832287496	0.0110675469572925\\
21.9812863734494	0.0103414754068254\\
22.060258786298	0.00924620224564876\\
22.1393947633787	0.00781217588144756\\
22.2255101693494	0.00589056468306381\\
22.3157348236216	0.00350907265561062\\
22.408919336699	0.000704862171474474\\
22.5125222613708	-0.00274726582634344\\
22.6452422749301	-0.00753073519963721\\
23.0417303283489	-0.0220909235764069\\
23.1510340662809	-0.0256167299858987\\
23.2437066654231	-0.0282832829052211\\
23.3322925764397	-0.030502728229358\\
23.5015891026805	-0.0345228348186879\\
23.5865651463452	-0.0370143801834146\\
23.6791414372839	-0.0400694959925332\\
23.7787370299291	-0.0436918656190244\\
23.8973981490484	-0.0483569770731691\\
24.0659833241504	-0.0553674600039287\\
24.2995667694583	-0.0650613252517189\\
24.4178168163229	-0.0696152075188863\\
24.5169175522648	-0.0730980754989403\\
24.6067002591987	-0.0759154254150474\\
24.6877135641367	-0.0781318122259265\\
24.7694912005601	-0.0800147871901586\\
24.8436055009371	-0.0813840924510174\\
24.9134171081312	-0.0823590991841527\\
24.986200831437	-0.083032260683332\\
25.0534819357377	-0.0833295876302103\\
25.1228214447728	-0.0832990220498999\\
25.1939106941538	-0.0829048492259759\\
25.2618219773336	-0.0821808133451754\\
25.3314829621529	-0.0810838822853341\\
25.3975543058924	-0.0797133281200502\\
25.4673261069476	-0.0779213177327236\\
25.5344513207435	-0.0758697251175988\\
25.6066214155467	-0.0733158994765546\\
25.6779968087409	-0.0704488621307959\\
25.7492794360363	-0.0672629734192611\\
25.8275567511942	-0.0634152924408511\\
25.9068851822997	-0.0591706687190339\\
25.9876917386046	-0.0545223651078715\\
26.0773002502623	-0.0490284791409579\\
26.1767819661158	-0.0425733434160946\\
26.2873608150649	-0.0350475405487245\\
26.4269647691983	-0.0251794194774959\\
26.9030594769996	0.00881763902837918\\
27.021811202169	0.016772174621309\\
27.1280873683345	0.0235468944067705\\
27.225722503734	0.029434424692127\\
27.3171785465289	0.0346238510659944\\
27.4101447662409	0.0395503933094972\\
27.4964277840875	0.0437905460311043\\
27.5853276536978	0.0478113854742688\\
27.6683425320399	0.0512390350921166\\
27.7511696695896	0.0543395271475902\\
27.8345991394645	0.0571382945709402\\
27.9164590754587	0.0595688970337491\\
28.0034427198926	0.0618125821271889\\
28.0874910282266	0.0636542722333857\\
28.1765902207824	0.0652648196561572\\
28.2625371703016	0.0664950154273285\\
28.3548131465994	0.0674753667840164\\
28.4493076372345	0.0681299323352746\\
28.5448733670303	0.0684509079308597\\
28.6440633655462	0.0684429171065517\\
28.7492787727161	0.0680808147007923\\
28.8573355955713	0.0673603938961733\\
28.9703746712281	0.0662643454771228\\
29.0924080631866	0.0647322995982265\\
29.223744835791	0.0627353123916663\\
29.3715903221537	0.060134745417038\\
29.550008783347	0.0566319441493519\\
29.8425911095113	0.0504740610692807\\
30.0931571562743	0.0453389511132869\\
30.2573447277347	0.0423207657472204\\
30.2647482875259	0.0420651810807584\\
30.3983503332359	0.0339570171548189\\
30.5499039283393	0.024833873438368\\
30.6417850030205	0.0196385567083439\\
30.7221729642479	0.015432497324106\\
30.7965859629501	0.0118881011619223\\
30.8647958761611	0.00897419354959084\\
30.9335721168957	0.00638819625740439\\
30.9959946751836	0.00436268032087384\\
31.0630631544141	0.00253656088534626\\
31.1288768295989	0.00109939644161727\\
31.1950283887034	4.65473382149639e-06\\
31.2615590413122	-0.000753340518052426\\
31.3310075390661	-0.00119563698495284\\
31.4021092508219	-0.00130505207371101\\
31.4764930796373	-0.00108305686119792\\
31.5586812282717	-0.000489346097197085\\
31.6511867626338	0.000537344844453003\\
31.7600397337288	0.00210001655187142\\
31.9376913258969	0.00504940475517657\\
32.1088187460872	0.00777312640721561\\
32.21754808502	0.00917402263777944\\
32.3150965244945	0.0100900167889506\\
32.4041804866251	0.0105824047344925\\
32.4901433686027	0.0107068931359748\\
32.5735908875838	0.0104747863919243\\
32.6560533107613	0.00989211302443493\\
32.737852913589	0.00896558040537343\\
32.8193171481296	0.00770544309889232\\
32.9075320478285	0.005979950348042\\
32.9942570703626	0.00394678394469139\\
33.0921061298439	0.00129988286797555\\
33.2010846046834	-0.00201045996284677\\
33.3293275531821	-0.00625700689873554\\
33.7745564322871	-0.0213320686046217\\
33.8851799594932	-0.0245393528832594\\
33.9838923654762	-0.0270467222959354\\
34.0762527165863	-0.0290375658441633\\
34.1552363501689	-0.0307479361325278\\
34.2422482959584	-0.0329562115080506\\
34.339021506389	-0.0357711395880358\\
34.4451330422502	-0.0392288210335821\\
34.5588268211247	-0.0432722683634807\\
34.7116296553883	-0.0490606366063062\\
35.0096800561681	-0.0604360785423737\\
35.1270152381257	-0.0645381442784583\\
35.2254342096107	-0.0676593151910438\\
35.3147117167717	-0.0701709977205454\\
35.3954081647203	-0.072135072744409\\
35.4770898656098	-0.0737905232736082\\
35.5510553374738	-0.0749744938196741\\
35.6283060601019	-0.0758685121010672\\
35.7055494071359	-0.0763934422728028\\
35.7765726597902	-0.0765372319296986\\
35.8474210006063	-0.0763480986849601\\
35.9186793067548	-0.075817008172784\\
35.9865641054393	-0.0749905512248787\\
36.0540376494186	-0.0738593749478582\\
36.1248176221567	-0.0723438457331511\\
36.1992058840641	-0.0703944874839166\\
36.2715604786911	-0.0681573064064551\\
36.3494552874265	-0.0653869702803505\\
36.4271784128783	-0.0622670164397192\\
36.5054234695354	-0.0587893163827076\\
36.5913851743684	-0.0546091934997506\\
36.6790999196738	-0.0499933846337797\\
36.7761084385443	-0.0445277258433521\\
36.8834884157797	-0.0381083080004032\\
37.0028553084658	-0.0306202272675549\\
37.1611386162774	-0.0203235575040708\\
37.5357677864473	0.00418469217822803\\
37.6630541999939	0.0121027002316723\\
37.7776356209792	0.0188957962317318\\
37.8802029034455	0.0246508648081658\\
37.982637940495	0.0300498366401172\\
38.0773529685952	0.0347036789819981\\
38.1721419455664	0.039015072969832\\
38.2621744422739	0.0427759066689362\\
38.3544572113654	0.0462835658549849\\
38.4434492234937	0.049328388468858\\
38.5341066649942	0.0520877764654344\\
38.6241909307942	0.0544891298122394\\
38.7122604061226	0.056512889430465\\
38.8059983642302	0.0583224249028333\\
38.89691685363	0.059747592596004\\
38.9954078481625	0.0609378523489852\\
39.0889062012411	0.0617417429668308\\
39.1940603108954	0.0622855389034243\\
39.2975198235733	0.0624697834487478\\
39.4071389805769	0.062311141065841\\
39.5193904050684	0.0618005787454621\\
39.6402805105967	0.0608926683540716\\
39.7654634773051	0.0596033401110176\\
39.8997117449049	0.0578779491978736\\
40.0450356525127	0.0556774474167412\\
40.2210135303132	0.0526595544260573\\
40.4702531218354	0.0479979769501568\\
40.8077166795558	0.0417057027181187\\
40.9690882002801	0.0387138128961908\\
41.1115575539225	0.030697905512028\\
41.2649199664651	0.0221734205059576\\
41.3577615879891	0.0173529761826003\\
41.4383429343932	0.0134988240743326\\
41.5128700574608	0.0102676126462526\\
41.5811834125036	0.00762324719868701\\
41.6500265170372	0.00528987965847705\\
41.717793633433	0.00333322370680378\\
41.7864173632411	0.00170359232526351\\
41.8553962893788	0.000422646874028487\\
41.9241492838754	-0.000504506790562687\\
41.9950002207524	-0.00110887878115307\\
42.0674603084029	-0.00138033122631498\\
42.1446293766016	-0.00131657446559785\\
42.2271463134317	-0.000892644835872147\\
42.3172731261727	-7.75051026380424e-05\\
42.4227806455496	0.00123362465459564\\
42.5679150770274	0.00341450508038577\\
42.8147456786366	0.00715591968440066\\
42.9300929440439	0.00852789021224964\\
43.0303493506078	0.00937252487040041\\
43.1197833269487	0.00979315069218245\\
43.2082616899718	0.00986328955306703\\
43.2930901067457	0.00958612509491275\\
43.376100580628	0.00897908417055504\\
43.4640427571791	0.00797427538285689\\
43.5502383171915	0.00663928839178141\\
43.6420509050573	0.00485949046920808\\
43.736642363355	0.00267707920943394\\
43.8394141765971	-3.72655019234003e-05\\
43.9589995111644	-0.00354650555942726\\
44.1251155158198	-0.00880205891638752\\
44.3962337691169	-0.0173941565923244\\
44.5179973663785	-0.0208884311023141\\
44.6267017385642	-0.0236611351789193\\
44.7217592640943	-0.0257525763048818\\
44.9181739358821	-0.0298455158068762\\
45.0166688040389	-0.0324204686144256\\
45.1220667536661	-0.0355248224394913\\
45.2439514506033	-0.0394690904212212\\
45.4044978535808	-0.0450346202835092\\
45.7104584339658	-0.0557051931646555\\
45.8274233631034	-0.0594176092809278\\
45.9345902066067	-0.0624752359909948\\
46.0239121923839	-0.0647069567170249\\
46.1139175327125	-0.0666143366232248\\
46.1965064769803	-0.068027008357916\\
46.2790453555772	-0.0690868960159179\\
46.3590807902781	-0.0697561419509043\\
46.434532583498	-0.0700479092270783\\
46.50724279597	-0.0700070032078415\\
46.5804813137283	-0.0696392212273125\\
46.652946788393	-0.068949000984361\\
46.7279046071706	-0.0678928401137497\\
46.7992515547698	-0.0665666839044832\\
46.8744031832353	-0.06483705544602\\
46.9539673319485	-0.0626438344726381\\
47.0324343814235	-0.060130078285404\\
47.1107362589249	-0.0572916960611494\\
47.196168169631	-0.0538428347994753\\
47.2828601389487	-0.0499976889600404\\
47.3782452049748	-0.0454065250573947\\
47.4762790486284	-0.040345994086735\\
47.5849102154025	-0.0344038035320651\\
47.7135292599859	-0.0270237422598782\\
47.8918309212051	-0.0164242443569904\\
48.202114558949	0.00205114595185307\\
48.3443433392152	0.0101413434113198\\
48.4662989319041	0.0167364464179727\\
48.5791684836012	0.0224922609378169\\
48.6811310040488	0.027363857023488\\
48.7799071455074	0.031758981691155\\
48.8803760343261	0.0358817517774952\\
48.9754815716302	0.0394480421637198\\
49.0687061819615	0.0426182910632775\\
49.1613339864836	0.0454452486392256\\
49.2548917611754	0.0479733112539975\\
49.3475772291205	0.0501559988395215\\
49.4466877102766	0.0521414259130779\\
49.5441078747867	0.0537504489832088\\
49.6393619056156	0.0550057053901014\\
49.7420495485496	0.0560203259802279\\
49.8464272278967	0.0567080958971999\\
49.9513729306303	0.0570691071769431\\
50	0.0571295515852484\\
};
\addlegendentry{Full System (6D)}

\addplot [color=red, dashed, line width=2.0pt]
  table[row sep=crcr]{%
0	0.00928064636582349\\
0.0906430398257925	0.0109331321164419\\
0.166795065555966	0.0120518101453015\\
0.249582863066173	0.0129143668181229\\
0.313290509793362	0.0132971834224236\\
0.379243276658507	0.0134222716198309\\
0.453253946082015	0.0132282371646752\\
0.522902171614518	0.0127215184388803\\
0.588160686656295	0.0119614698225305\\
0.669502229898434	0.0106269682676015\\
0.731002245770853	0.00933285136566298\\
0.798631042000075	0.00762977691351807\\
0.882178615639667	0.00513677290970804\\
0.971437349039853	0.00204153480052582\\
1.05619962455088	-0.00124284907976602\\
1.17415937450286	-0.00621789194241984\\
1.36151166282321	-0.0146114978470067\\
1.54795249440706	-0.0229461998837834\\
1.65379682629627	-0.0273559294269887\\
1.75185902493698	-0.031039275699861\\
1.83219273153973	-0.0337137142989761\\
1.92844151828168	-0.0365573198102283\\
2.11245556081636	-0.0417716795073488\\
2.18681772664734	-0.0442059897569962\\
2.27077886382673	-0.0473211047482707\\
2.35193490948029	-0.050683528325024\\
2.47646514379051	-0.0562965306738903\\
2.93565019202323	-0.0775135744787434\\
3.04106861057019	-0.0820139136201377\\
3.12756310626091	-0.0854122249378904\\
3.21316077399331	-0.0884216266678095\\
3.28128943702183	-0.0905166705115263\\
3.34951552469956	-0.0923165459137039\\
3.41830848529536	-0.0938072250910054\\
3.48816211091668	-0.0949712625563492\\
3.55956853289975	-0.0957847549883368\\
3.63300001162341	-0.0962152950025512\\
3.68958366229526	-0.0962620292972431\\
3.74749687559547	-0.0960503870751452\\
3.80764517573294	-0.0955507718296147\\
3.8693888249486	-0.0947398038726348\\
3.93275418732057	-0.0935924650814002\\
3.99776762692697	-0.0920832487261904\\
4.06530305633931	-0.0901602176493483\\
4.13460031820098	-0.0878133657129254\\
4.20576786150078	-0.0850145550237187\\
4.2789141352275	-0.0817365061399755\\
4.35571975587292	-0.0778708142930569\\
4.43583666588459	-0.0733975965358979\\
4.52009941193464	-0.0682386194396472\\
4.60934254069515	-0.0623114112371326\\
4.67757590420072	-0.0574918410284155\\
4.74580926770628	-0.0524529333982358\\
4.83737270572966	-0.0453973233191149\\
4.92893614375305	-0.0380715732399608\\
5.08099257146892	-0.0255075427255136\\
5.43228641209186	0.00381573888373765\\
5.55571061839984	0.0137113297174309\\
5.67515063131687	0.0228846171737729\\
5.76490132902109	0.0294579063234295\\
5.86074335623279	0.0361311580801598\\
5.93117670203044	0.0407872985949211\\
6.00814870881796	0.0456198961641192\\
6.10856026214543	0.0515020616793009\\
6.19700639814744	0.0562706740529961\\
6.27656815635233	0.0602194211997684\\
6.35130004347626	0.0636281660460156\\
6.42444293106269	0.0666788126272166\\
6.49637921789913	0.0694011462629121\\
6.56749130277311	0.0718201405808543\\
6.63802435887438	0.0739519398250508\\
6.70838800425017	0.0758137671389036\\
6.77856535441974	0.077408460008229\\
6.84853952490233	0.078739701404011\\
6.9407387919077	0.0801034670424343\\
7.03065247018762	0.0810129798886763\\
7.1171915643431	0.0815060082572856\\
7.19994654371135	0.0816392079652033\\
7.279218966319	0.0814724651093428\\
7.37431022594444	0.0809219700495376\\
7.46674492635191	0.0800625418321701\\
7.57584331518808	0.0787079978027592\\
7.70251447952494	0.0767792251896182\\
7.83036931399931	0.0745232992022693\\
7.96180636893421	0.0718618245297122\\
8.0822253519075	0.0690432837902293\\
8.21880872806237	0.0654353541264427\\
8.38790858706601	0.0609207117748412\\
8.49629026777594	0.0583982422544267\\
8.7035824706548	0.0538782390578731\\
8.76261114325707	0.0522245884320682\\
8.82774134188213	0.0499819772882475\\
8.84882240105447	0.0490761429192403\\
8.88843757529139	0.0468831194576822\\
8.94469306022055	0.0433677811686408\\
9.00129118138996	0.0394484934889121\\
9.10155810744479	0.0320585602040566\\
9.18787393156445	0.0258445583567735\\
9.25220053210235	0.0215681014461921\\
9.30793407522643	0.0181691139876108\\
9.36374763176637	0.0150508531769589\\
9.41960119501424	0.0121981176926056\\
9.49450188468358	0.00876105769035718\\
9.55164326394799	0.00642933112067112\\
9.61986597372865	0.00398802116352215\\
9.67423855222913	0.00232997623441378\\
9.73750190270964	0.000747548092533634\\
9.80988399290529	-0.000582588622570768\\
9.87100961533839	-0.00130434803482871\\
9.93213523777148	-0.00167486015912743\\
9.99161015455255	-0.00172726597701001\\
10.0721089041241	-0.00138121119487522\\
10.1564671884685	-0.000603539066403869\\
10.2466882609765	0.000581551579657003\\
10.3663347871251	0.0025502592908424\\
10.5443795876626	0.00594464433552133\\
10.6986760458755	0.00882865660766896\\
10.7995558210155	0.0103809441856484\\
10.8958802402834	0.0114503652741149\\
10.9677867485402	0.0119315727104947\\
11.0408699377544	0.0121202511428393\\
11.1145752205032	0.0119958084230234\\
11.1890755162352	0.0115470408189395\\
11.2645437443991	0.0107612545540263\\
11.3450602171171	0.00955467813636091\\
11.4055853722699	0.00839747860793949\\
11.4729580340968	0.00685999069566634\\
11.5697141254369	0.00421259768467053\\
11.6611763190601	0.00128286336411776\\
11.7463536186501	-0.00174968222660254\\
11.8657499873375	-0.00634715719557732\\
12.08524593074	-0.0152688782048003\\
12.2373585467579	-0.0213573294388851\\
12.3486719042134	-0.0254784416691081\\
12.4559027329633	-0.0289954994677259\\
12.549404567983	-0.0316617232452998\\
12.6992262625033	-0.0354452097414963\\
12.8061015731418	-0.0383063861516533\\
12.8823936011276	-0.040640816116742\\
12.9710643889167	-0.0437218936731156\\
13.0603494656669	-0.0471721145109427\\
13.2178424464326	-0.0537191859844768\\
13.5462143363425	-0.0674800502494719\\
13.6878508448487	-0.073081276782176\\
13.7919720209274	-0.0768416849450801\\
13.8776990249869	-0.079598850683908\\
13.9628283180622	-0.0819625335122183\\
14.0310837494221	-0.083553789196344\\
14.0999656840051	-0.0848642012984584\\
14.1699696851312	-0.0858774530846631\\
14.2415868656704	-0.0865711062842536\\
14.3152862227024	-0.0869148541696632\\
14.3914932348449	-0.0868695131775539\\
14.4705614436633	-0.0863870289350075\\
14.5315909045167	-0.085709177284059\\
14.595050701671	-0.0847208978820859\\
14.6601703732958	-0.0834059337762483\\
14.7269707134972	-0.0817409681133512\\
14.7954725163809	-0.0797032222867031\\
14.8667069393153	-0.0772337197744477\\
14.9400570721484	-0.0743250511439157\\
15.0158172834819	-0.0709431302476276\\
15.0942819419176	-0.0670536697435296\\
15.1784566727309	-0.062470312655563\\
15.2686330946083	-0.0571279000015963\\
15.3676122608386	-0.0508086723615335\\
15.4781952247106	-0.0432755142368535\\
15.6023255514806	-0.0343608559507373\\
15.7616886252859	-0.022461150409292\\
16.1060327383757	0.00349424032160783\\
16.2295636098084	0.012436832844017\\
16.3489868156238	0.0207175731399261\\
16.4398765138555	0.0267252833608467\\
16.5359834735656	0.0327614394510434\\
16.6407387159879	0.0389339573216532\\
16.7248031591236	0.0435584116778145\\
16.8088676022592	0.0478748296149334\\
16.8909753783536	0.0517827071838397\\
16.9691698203658	0.055213760341033\\
17.043315870575	0.0582004124167739\\
17.1400726718859	0.0617020723710624\\
17.2350470504094	0.0646999605626917\\
17.3291697463226	0.0672406545487334\\
17.4228972566758	0.0693463903693825\\
17.5161808081477	0.0710250382058533\\
17.6083657323983	0.0722800465461901\\
17.6983311381155	0.0731243637963601\\
17.7850284496019	0.0735910365937755\\
17.8680169379007	0.0737298651316394\\
17.9672038037862	0.0735224480627323\\
18.061772240421	0.0729783278047549\\
18.1721467711884	0.0719741859781209\\
18.2990009275843	0.0704293445534816\\
18.4442004838168	0.0682802917350145\\
18.5734873260213	0.0660659218427355\\
18.6890667514842	0.0637973673641525\\
18.8159583369415	0.0609444696455625\\
19.1208066667358	0.0537197357817689\\
19.2550176837592	0.0511191884707358\\
19.3746999185796	0.0487031511202929\\
19.4340602933373	0.0472099472528598\\
19.4998360236459	0.0451734563720834\\
19.5273871894228	0.0440606885021424\\
19.5733687917187	0.0417060515155683\\
19.6295354055562	0.0384553358565825\\
19.7013064662311	0.033860396608425\\
19.8738011695862	0.0225439969388646\\
19.9402432701566	0.018629091410638\\
19.995258342155	0.0156665346923148\\
20.0700054235291	0.0120347599876212\\
20.1447525049033	0.00881742089566728\\
20.2205773508046	0.00595182615560219\\
20.2878213796114	0.00375531194552536\\
20.3601725108325	0.00178621836118253\\
20.4350771932609	0.000216345824760822\\
20.5099818756893	-0.00085451566658179\\
20.5732001263761	-0.00137538346818644\\
20.6339735607003	-0.00156775385732288\\
20.7139136876561	-0.00141546149868077\\
20.7965234644883	-0.000859535071612072\\
20.8843572360102	7.18199315556944e-05\\
21.0019936355725	0.00171184506181277\\
21.1506457063297	0.00418873956245847\\
21.3855984282265	0.00818455814091834\\
21.4864128832058	0.0095378141660305\\
21.5827454851138	0.0104440277134472\\
21.6549457819548	0.0108298727239244\\
21.7279024144193	0.0109466096818238\\
21.8018304458817	0.0107790258516687\\
21.8766166875648	0.0103156485177678\\
21.9533094157509	0.00953337635657192\\
22.0329569067223	0.00839113522492596\\
22.1295053119846	0.0065579285043782\\
22.2270407918501	0.00422431733124284\\
22.3208401932201	0.00156793668519839\\
22.4357498171131	-0.00212878807945316\\
22.5562516553533	-0.00635851120004816\\
22.9744874555628	-0.0214742720854417\\
23.0763844601953	-0.0246839389279216\\
23.1687346796524	-0.0272482569598651\\
23.280556234266	-0.0299512886795128\\
23.4810748393173	-0.0346437104965887\\
23.5583816542102	-0.0367882367294996\\
23.6497538483136	-0.0396725973751373\\
23.7479757008469	-0.0431311119088704\\
23.9161265468915	-0.0494773948932163\\
24.240925065025	-0.0617357721993841\\
24.3819420564987	-0.0667125477505053\\
24.4857594156203	-0.070025632041677\\
24.5712932747007	-0.0724331245709138\\
24.6564112589982	-0.0744790724880318\\
24.7420235783863	-0.0761423339109797\\
24.8114777303689	-0.0771802230824008\\
24.8823259264628	-0.0779392418937306\\
24.9550488948098	-0.0783953583060963\\
25.0300941210652	-0.078516849391292\\
25.1078502880403	-0.0782638553055008\\
25.1886173950815	-0.0775890623699311\\
25.251186046652	-0.0767760811101965\\
25.3153900941766	-0.0756777779450033\\
25.3820978777681	-0.0742535400305329\\
25.4505309660615	-0.0724941151099898\\
25.5207702875601	-0.0703773143374917\\
25.592896770767	-0.0678817107115464\\
25.6684054952908	-0.0649292454403891\\
25.7468809241466	-0.0615072193149118\\
25.8289924266958	-0.0575625664765766\\
25.9154093722998	-0.0530398065405109\\
26.0124819422292	-0.0475530853586648\\
26.1249739856229	-0.0407370597767951\\
26.2051085057068	-0.0356373614441736\\
26.3174341400931	-0.02823228275048\\
26.5317467606214	-0.013640985578256\\
26.7658337104675	0.00227016221930398\\
26.8902189584481	0.010432884622972\\
27.0099505139284	0.0179678720720915\\
27.1312213823043	0.0251946193896657\\
27.2295199837862	0.0307060845202756\\
27.3389680244313	0.0364388306783709\\
27.4215347290557	0.0404626998269322\\
27.5041014336801	0.0442152355006087\\
27.5848654616771	0.0476147106789071\\
27.6865639010893	0.0515039553621506\\
27.7841049328263	0.0548165348469283\\
27.879479408619	0.0576554880234141\\
27.9737785254944	0.0600719164479599\\
28.0675894566111	0.0620912048180813\\
28.1610274788264	0.063723594772533\\
28.2536347212092	0.0649723899784362\\
28.3444070739331	0.0658453520989468\\
28.432200980917	0.0663662186863121\\
28.5163385539012	0.0665762579041242\\
28.616386870777	0.0664758956807674\\
28.7121337594127	0.0660512397547777\\
28.8231274035008	0.0652076303966922\\
28.9501530523911	0.0638690448805619\\
29.0950986154655	0.0619815634672989\\
29.2423895358253	0.0597338225291182\\
29.3779443530226	0.0573138396967181\\
29.507018439079	0.0546502241527236\\
29.7895286434552	0.0485971327569246\\
29.9253677994466	0.046212350345975\\
30.0463411086446	0.0440140921168393\\
30.1055187180693	0.0426763255529181\\
30.1708892020368	0.0408596663380862\\
30.2001375898773	0.0398062838486339\\
30.250265619408	0.0374772697143868\\
30.3065048391875	0.0345147318751771\\
30.3780081581534	0.030363369379991\\
30.5353156835297	0.021007889994948\\
30.6006648669372	0.0174728901843793\\
30.6729899439591	0.0139333207843038\\
30.7507163558311	0.0105559777186031\\
30.8277951015818	0.00760291119454592\\
30.9058875165155	0.0049923000506098\\
30.9711093296952	0.00311772393998666\\
31.0496175534022	0.00126582984966461\\
31.1281257771091	-0.000106549512032927\\
31.2016730357385	-0.000942247124847029\\
31.2652983642128	-0.00132456422586813\\
31.3447305860121	-0.0013989637269276\\
31.4251370718755	-0.00108677825284786\\
31.5097371762345	-0.000425958812179772\\
31.6245335049938	0.000868579825926474\\
31.7705735504953	0.00295155754570686\\
32.1098507098734	0.00805225516369745\\
32.2068169679008	0.00906339851486848\\
32.3025686454549	0.00968486423749937\\
32.3998248389409	0.00988625746327187\\
32.4737430914789	0.00973907201857571\\
32.5485085651817	0.00932495316152426\\
32.6251535445562	0.00862340005976137\\
32.7047103141092	0.00759793298981037\\
32.8008266530353	0.00595680773161433\\
32.8985172753796	0.00385328283486785\\
32.9926264234879	0.00145163500503287\\
33.1075417802468	-0.00188279589568907\\
33.2597341344328	-0.00673774459463772\\
33.645498000998	-0.0193310098987425\\
33.7465471793041	-0.0222137824523614\\
33.8382082579201	-0.0245229974367049\\
33.9578035058091	-0.0271343194682387\\
34.1536895605569	-0.0312757470641358\\
34.2519912763943	-0.0337791299583188\\
34.3492721915874	-0.0366408491666093\\
34.4603220708009	-0.0402639509846381\\
35.0354908143977	-0.0596501939534519\\
35.1396748130629	-0.0627144600497189\\
35.2423330557675	-0.0653662310851857\\
35.3274693989327	-0.0672224006665232\\
35.4129640486091	-0.0687315253852816\\
35.5000069937776	-0.0698749184449738\\
35.5712124965235	-0.0705022517700229\\
35.6443957892778	-0.0708504787686692\\
35.7199945889742	-0.0708893333207499\\
35.798381635104	-0.0705812214948338\\
35.8798350390471	-0.0698820956528792\\
35.9645236771459	-0.068743500789779\\
36.029899634777	-0.0675771235233782\\
36.0978078081196	-0.0661012355184738\\
36.1674985355874	-0.0643088414286623\\
36.2390943107654	-0.0621788192897696\\
36.3127176272386	-0.0596906386274725\\
36.3901435933074	-0.0567595913382277\\
36.4710589604269	-0.0533692230143714\\
36.5563890287032	-0.0494564943095881\\
36.6470590982424	-0.0449547160107286\\
36.7522599344196	-0.0393444948892778\\
36.8375500211543	-0.0345407469804044\\
36.9882194451812	-0.025636240294773\\
37.1930667264858	-0.0130487618881148\\
37.4552658462951	0.00303303884491157\\
37.6076379212652	0.0119788485027712\\
37.7270141202963	0.0186394042529017\\
37.8508730049088	0.0251479984518923\\
37.9505468895691	0.0300474335780621\\
38.0640415934091	0.0352240946108537\\
38.1491626212892	0.0388085808863394\\
38.2576398976192	0.0429895857731353\\
38.3594487696099	0.046506829484521\\
38.4570310622262	0.0495004186538992\\
38.5524247829103	0.0520654002827499\\
38.646737876169	0.0542485352212196\\
38.7405652912445	0.0560729463208389\\
38.8340246256884	0.057547888195522\\
38.9266488143067	0.0586761967640825\\
39.0174178483002	0.0594648176676458\\
39.1051915295281	0.0599354011760624\\
39.2101059073193	0.0601342136237903\\
39.3089273202908	0.0599794157608358\\
39.4224779779291	0.0594308430243302\\
39.5507079787731	0.0584094482072928\\
39.6955275283158	0.0568677846950791\\
39.8597570988777	0.0547510122485164\\
40.011158028336	0.0524399308042476\\
40.1350108356035	0.0502281855397229\\
40.5276377869169	0.0428234032494856\\
40.7337086260481	0.0394781835880593\\
40.8117302979196	0.0377548856134098\\
40.866067759899	0.0362261391673684\\
40.9160802101078	0.0341688788387344\\
40.9722791261653	0.0315350798632963\\
41.0440540468617	0.0278020319164796\\
41.2310864826828	0.0178332347335157\\
41.2992157657049	0.014613690382248\\
41.3786452139075	0.0112589159196403\\
41.4620999835911	0.00816568746215296\\
41.5373706484766	0.0057214407987729\\
41.6219748706407	0.0033629853208339\\
41.6980089252243	0.00162412714028193\\
41.7641197900201	0.000433842746502933\\
41.8302306548159	-0.000440751108754966\\
41.8948454834742	-0.00098982108191592\\
41.9762909811909	-0.00127555208398178\\
42.0560372211947	-0.0011688542509205\\
42.1383416580271	-0.000731999020089802\\
42.2488447422511	0.0002456582743946\\
42.3676066768769	0.00164370189839502\\
42.544130781297	0.00409696873661858\\
42.7454225461489	0.00684100509418784\\
42.8430144675883	0.00787526568402797\\
42.9391701363962	0.00857580225767407\\
43.0358307797559	0.0089064941937238\\
43.133897763801	0.00883568505620502\\
43.2335552179755	0.00833928617421975\\
43.3107083984601	0.00766084749651696\\
43.391229045457	0.00667847291845902\\
43.4923115883094	0.00505147606462941\\
43.5890665946562	0.00310119824524691\\
43.6805197881102	0.00094321789312346\\
43.7955474510662	-0.00211405187903324\\
43.9488796270872	-0.00656282515095796\\
44.3011835689739	-0.016981401813112\\
44.4346036098636	-0.02042047680618\\
44.548722373743	-0.0229230871696302\\
44.6863702401772	-0.0255167765407052\\
44.8288804056372	-0.0282899880254419\\
44.9279482960659	-0.0305774722900622\\
45.0253276310794	-0.0331741693200769\\
45.1669103500095	-0.0373895552679571\\
45.6584605913859	-0.0524481618089183\\
45.7806506241448	-0.0558164400486021\\
45.883748250991	-0.0583343116097197\\
45.9689428726964	-0.060122016238374\\
46.0542571847783	-0.0616042467208189\\
46.1406152216616	-0.0627622527013685\\
46.2288842946448	-0.063569992935129\\
46.3016067365339	-0.0639406114872898\\
46.3766562997402	-0.0640382442082341\\
46.4544222838072	-0.0638303164444878\\
46.5352045561585	-0.0632781655589767\\
46.6191932181277	-0.0623387070619827\\
46.7064859497827	-0.0609666911166045\\
46.797167515824	-0.0591160160496358\\
46.8676510454884	-0.0573814534792589\\
46.9401209523462	-0.0553331644709658\\
47.0147307694759	-0.0529516300758814\\
47.0934661948386	-0.0501503820651195\\
47.1760949900171	-0.046910488684361\\
47.2637642020522	-0.0431629400175311\\
47.3576208779844	-0.0388342383366975\\
47.4697580686737	-0.0332965176808315\\
47.5662119092535	-0.0282762835095056\\
47.7434362403036	-0.0186304070827674\\
48.1887766613293	0.00598858744399422\\
48.3391214945616	0.013799767212042\\
48.4597311700048	0.0197127316833416\\
48.5558985082006	0.024160011902751\\
48.660780606829	0.0287098562021484\\
48.7728083515352	0.0331966139295758\\
48.8848360962414	0.037275276802994\\
48.989050278212	0.0406878359424567\\
49.0878889282281	0.0435753552576941\\
49.18393224456	0.0460506352807784\\
49.2785405044633	0.0481680628893173\\
49.372496976495	0.0499555100932199\\
49.4661455686742	0.0514257467851493\\
49.5593020840503	0.0525826031769512\\
49.6511502037693	0.0534287061952909\\
49.7619144762097	0.0540681408061587\\
49.8672729320056	0.0543027030975125\\
49.967277622339	0.0542067327253477\\
50	0.0541121109684539\\
};
\addlegendentry{Arc-length ROM (2D)}

\end{axis}

\begin{axis}[%
width=0.7\textwidth,
height=0.7\textwidth,
at={(0in,0in)},
scale only axis,
xmin=0,
xmax=1,
ymin=0,
ymax=1,
axis line style={draw=none},
ticks=none,
axis x line*=bottom,
axis y line*=left
]
\end{axis}
\end{tikzpicture}%